\documentclass[a4paper,10pt]{amsart}
 
\usepackage[english]{babel}
\usepackage[utf8]{inputenc}

\usepackage{bbding}

\usepackage{amssymb} 
\usepackage{ifthen}
\usepackage{amsmath} 
\usepackage{amsthm}  
\usepackage{xspace}
\usepackage{enumerate} 
\usepackage{mathrsfs}
\usepackage{pdfpages}
\usepackage[all]{xy}
\usepackage{wasysym}       
\usepackage{bbm} 
\usepackage{bm} 
\usepackage{color} 
\usepackage{braket} 
\usepackage{rotating}
\usepackage{hyperref}

\newcommand{\mathgrc}[1]{\mathbf{#1}}

\newtheorem{theorem}{Theorem}[section] 
\newtheorem{proposition}[theorem]{Proposition}  
\newtheorem{corollary}[theorem]{Corollary}    
\newtheorem{lemma}[theorem]{Lemma}        

\theoremstyle{definition}           

\newtheorem{remark}[theorem]{Remark}   
\newtheorem{notation}[theorem]{Notation}   
    
\newtheorem{definition}[theorem]{Definition} 
\newtheorem{example}[theorem]{Example} 
\newtheorem{assumption}[theorem]{Assumption}

\newcommand{\protectreference}{\ref{thm-diagrams}}

\usepackage{color}

\newcommand{\tildeh}{ h}
\newcommand{\Io}[1][n]{\ensuremath{\mathbb{I}_{#1,0}}\xspace}
\newcommand{\Ii}[1][n]{\ensuremath{\mathbb{I}_{#1,1}}\xspace}
\newcommand{\Iii}[1][n]{\ensuremath{\mathbb{I}_{#1,2}}\xspace}
\newcommand{\e}[2]{\ensuremath{\mathfrak{e}_{#1,#2}}\xspace}
\newcommand{\eo}[1][n]{\e{#1}{0}}
\newcommand{\ei}[1][n]{\e{#1}{1}}
\newcommand{\eii}[1][n]{\e{#1}{2}}
\newcommand{\FFFF}[2]{\ensuremath{\mathfrak{f}_{#1,#2}}\xspace}
\newcommand{\Fo}[1][n]{\FFFF{#1}{0}}

\newcommand{\Fii}[1][n]{\FFFF{#1}{2}}
\newcommand{\ksix}{\ensuremath{K_{\mathrm{six}}}\xspace}
        
        \newcommand{\kk}{\ensuremath{\mathit{KK}}\xspace}

        \newcommand{\kkE}{\ensuremath{\mathit{KK}_{\mathcal{E}}}\xspace}
        \newcommand{\catkkE}{\ensuremath{\mathbf{KK_{\mathcal{E}}}}\xspace}
        
        \newcommand{\kE}{\ensuremath{\underline{\mathit{K}}_{\mathcal{E}}}\xspace}
        
        \newcommand{\Homsix}{\ensuremath{\operatorname{Hom}_{\mathrm{six}}}\xspace}
        \newcommand{\Extsix}{\ensuremath{\operatorname{Ext}^1_{\mathrm{six}}}\xspace}

\newcommand{\scstes}{standard cyclic six term exact sequence\xspace}

\newcommand{\defemph}[1]{\textbf{#1}}

\newcommand{\flip}{\ensuremath{\operatorname{flip}}\xspace}
\newcommand{\Sflip}{\ensuremath{\operatorname{S-flip}}\xspace}
\newcommand{\Cflip}{\ensuremath{\operatorname{C-flip}}\xspace}
\newcommand{\sets}[2]{\ensuremath{\Set{#1 | #2}}\xspace}
\newcommand{\Matn}[1][n]{\ensuremath{\textsf{M}_{#1}}\xspace}

\newcommand{\Z}{\ensuremath{\mathbb{Z}}\xspace}

\newcommand{\N}{\ensuremath{\mathbb{N}}\xspace}
\newcommand{\C}{\ensuremath{\mathbb{C}}\xspace}

\newcommand{\uct}{the bootstrap category $\mathcal{N}$\xspace}

\newcommand{\ev}{\ensuremath{\operatorname{ev}}\xspace}

\newcommand{\mc}{\ensuremath{\mathfrak{mc}}\xspace}

\newcommand{\sfS}{\ensuremath{\mathsf{S}}\xspace}
\newcommand{\sfC}{\ensuremath{\mathsf{C}}\xspace}

\newcommand{\A}{\ensuremath{\mathfrak{A}}\xspace}
\newcommand{\B}{\ensuremath{\mathfrak{B}}\xspace}
\newcommand{\CC}{\ensuremath{\mathfrak{C}}\xspace}
\newcommand{\D}{\ensuremath{\mathfrak{D}}\xspace}
\newcommand{\E}{\ensuremath{\mathfrak{E}}\xspace}

\newcommand{\ca}{\mbox{$C\sp *$-}alge\-bra\xspace}
\newcommand{\cas}{\mbox{$C\sp *$-}alge\-bras\xspace}

\newcommand{\starhomo}{\mbox{$*$-}homo\-mor\-phism\xspace}
\newcommand{\starhomos}{\mbox{$*$-}homo\-mor\-phisms\xspace}

\newcommand{\starisos}{\mbox{$*$-}iso\-mor\-phisms\xspace}

\newcommand{\Hom}{\operatorname{Hom}}

\newcommand{\ident}{\mathrm{id}}

\newcommand{\simplecomplex}[2]%
{\ensuremath{\complex{Compl}^{#1,#2}}\xspace}
\newcommand{\proj}[2]%
{\ensuremath{\complex{Proj}^{#1,#2}}\xspace}
\newcommand{\free}[2]%
{\ensuremath{\complex{Free}^{#1,#2}}\xspace}
\newcommand{\cofree}[2]%
{\ensuremath{\complex{Cofree}^{#1,#2}}\xspace}
\newcommand{\complex}[1]%
{\ensuremath{\mathbf{#1}_\bullet}\xspace}

\newcommand{\mor}[1]%
{\ensuremath{\bm{#1}_\bullet}\xspace}

\newcommand{\obj}[2]%
{\ensuremath{\mathsf{#1}_{#2}}\xspace}

\newcommand{\map}[2]%
{\ensuremath{\mathsf{#1}_{#2}}\xspace}

\newcommand{\rep}[2]%
{\ensuremath{(\obj{#1}{},\map{#2}{})}\xspace}

\newcommand{\functorF}{\ensuremath{\mathsf{F}}\xspace}

\newcommand{\ie}{\textit{i.e.}\xspace}
\newcommand{\eg}{\textit{e.g.}\xspace}
\newcommand{\cf}{\textit{cf.}\xspace}

\newcommand{\bonkat}[1][271076]{%
  \ifthenelse{\equal{#1}{271076}}%
  {\cite{bonkat}}{\cite[#1]{bonkat}}}

\newcommand{\hooklongrightarrow}{\lhook\joinrel\relbar\joinrel\rightarrow}
\newcommand{\twoheadlongrightarrow}{\relbar\joinrel\twoheadrightarrow}

\newcommand{\extwithoutmaps}[3]{\ensuremath{#1\hookrightarrow#2\twoheadrightarrow#3}\xspace}
\newcommand{\extwithmaps}[5]{\ensuremath{#1\overset{#4}{\hooklongrightarrow}#2\overset{#5}{\twoheadlongrightarrow}#3}\xspace}

\newcommand{\ext}%
{\ensuremath{e\colon\B\hookrightarrow\E\twoheadrightarrow\A}\xspace}
\newcommand{\extprime}%
{\ensuremath{e'\colon\B'\hookrightarrow\E'\twoheadrightarrow\A'}\xspace}
\newcommand{\EXT}[4]%
{\ensuremath{#1\colon#2\hookrightarrow#3\twoheadrightarrow#4}\xspace}
\newcommand{\Ext}%
{\operatorname{Ext}}

        \advance \topmargin by -\headheight
        \advance \topmargin by -\headsep

        \advance \evensidemargin by -\oddsidemargin
        \advance \evensidemargin by -\textwidth 
        \advance \evensidemargin by -1in        
\begin{document}

\title{Ideal related $K$-theory with coefficients}

\author{S{\o}ren Eilers}


\address{
  Department of Mathematical Sciences \\
  University of Copenhagen \\
  Universitetsparken 5 \\ 
  DK-2100 Copenhagen \O \\ 
  Denmark
}

\email{eilers@math.ku.dk}

\author{Gunnar Restorff}

\address{
  Department of Science and Technology \\
  University of the Faroe Islands \\
  Nóatún 3 \\
  FO-100 Tórshavn \\ 
  Faroe Islands
}

\email{gunnarr@setur.fo}

\author{Efren Ruiz}

\address{
  Department of Mathematics \\
  University of Hawaii, Hilo \\
  200 W.\ Kawili St. \\ 
  Hilo \\
  Hawaii, 96720-4091 \\
  USA
}

\email{ruize@hawaii.edu}


\maketitle

\begin{abstract}
  In this paper, we define an invariant, which we believe should be
  the substitute for total $K$-theory in the case when
  there is one distinguished ideal.  Moreover, some diagrams relating the new groups to the
  ordinary $K$-groups with coefficients are constructed.  
  These diagrams will in most cases help to determine the new
  groups, and will in a companion paper be used to prove a universal
  multi-coefficient theorem for the one distinguished ideal case for a large class
  of algebras. 
\end{abstract}

\tableofcontents

\section{Introduction}
\noindent 
To characterize the automorphism groups of purely infinite \cas up to,
say, approximate unitary equivalence, one naturally looks at the work
of Dadarlat and Loring, which gave such a characterization of the
automorphism groups of certain stably finite \cas of real rank zero as
a corollary to their Universal Multi-Coefficient Theorem (UMCT),
\cf\ \cite{MR1404333}.
But even for nuclear, separable, purely infinite \cas with real
rank zero, finitely generated $K$-theory and only one non-trivial
ideal, there are substantial 
problems in doing so. 
The work of Rørdam (\cf\ \cite{MR99b:46108}) 
clearly indicates that the right invariant contains the associated six
term exact sequence in $K$-theory, and the work of Dadarlat and Loring
indicates that one should consider $K$-theory with coefficients in a
similar way. 

In the paper \cite{cmb}, the authors gave a series of examples showing
that the naïve approach --- of combining the six term exact sequence with
total $K$-theory --- does not work. 
There are several obstructions given in the paper, and they can even
be obtained using Cuntz-Krieger algebras of type (II) with 
exactly one non-trivial ideal. 

With this as motivation, a new invariant ---  \emph{ideal-related
  $K$-theory with coefficients} --- is defined, and we will argue that it 
should be thought of as the substitute for total $K$-theory,
when working with \cas with one specified ideal. 
It is easy to show that all the obstructions from the paper~\cite{cmb}
vanish when using this invariant.  
Furthermore, a lot of diagrams, which are part of the new
invariant, are exhibited. 
These diagrams can --- in many cases --- be very useful when computing
the new groups which go into the invariant. 
Also these diagrams are used in another paper (\cite{automorphism-se-gr-er}) by the three
authors, where {they} show a UMCT for $\kkE$ for a class of \cas including
all Cuntz-Krieger algebras of type (II) with one specified ideal --- in this case, the
invariant can actually be reduced quite a lot.

The paper is organized as follows. 
The first section contains an introduction. 
The second section sets up notation and conventions, and proves some
results related to suspensions, cones, and mapping cone sequences. 
The third section recalls definitions and results about homology and
cohomology theories for \cas. 
The fourth section exhibits relations between homology (and
cohomology) theories and mapping cone sequences. 
The fifth section gives some concrete examples of the results
developed in Section~\ref{sec-homologyandmappingcones}. 
In the sixth section, \emph{Ideal-related $K$-theory with coefficients} is defined. 
In the seventh section, some important diagrams involved in the
invariant are constructed. 
The eighth section contains the proof of Theorem~\ref{thm-diagrams} from
Section~\ref{sec-somediagrams}. 
The ninth section contains some examples of the theory developed in this paper. 

Parts of this paper have appeared in the second named author's
PhD-thesis \cite{grthesis}.

\section{Suspensions, cones, and mapping cones}

This section is devoted to setting up notation, recalling basic concepts, and proving some results that will be needed later.  Throughout the paper, $\N_0$, $\N$, and $\N_{\geq 2}$ will denote the set of non-negative integers, the set of positive integers, and the set of integers greater or equal to $2$, respectively. 
Moreover, $\Matn$ will denote the algebra of $n\times n$ matrices with complex entries.

\label{sec-mappingcones}

\begin{definition}
  Let \A be a \ca.  Define the \defemph{suspension} and the \defemph{cone}\footnote{Note that some authors place the algebra at $0$ rather than $1$ --- \eg\ Blackadar in \cite{MR99g:46104}} 
  of \A as 
  \begin{align*}
    \sfS\A&=\sets{f\in C([0,1],\A)}{f(0)=0,f(1)=0}, \\
    \sfC\A&=\sets{f\in C([0,1],\A)}{f(0)=0},
  \end{align*}
 respectively.
\end{definition}

\begin{remark}
  For each \ca \A, there is a canonical short exact sequence:
  $$\extwithoutmaps{\sfS\A}{\sfC\A}{\A}.$$
  It is
  well-known, that $\sfS$ and $\sfC$ are exact functors.  
\end{remark}

\begin{notation}
  Whenever convenient 
  $\sfC\sfC\A$, $\sfS\sfC\A$, $\sfC\sfS\A$, and $\sfS\sfS\A$ will be identified with subalgebras of $C([0,1]^2,\A)$
  by writing $f(x,y)$ for $(f(x))(y)$. 
  In this way $\ev_1 (f)$ will become $f(1,-)$ while $(\sfS\ev_1)(f)$ or
  $(\sfC\ev_1)(f)$ will be $f(-,1)$. 

The operation on $C([0,1]^2,\A)$ that \emph{flips} a
  function on $[0,1]^2$ along the diagonal will be denoted by $\flip$, \ie, 
  $\flip(f)(x,y)=f(y,x)$.  
\end{notation}

\begin{definition}
  Let \A and \B be \cas, and let $\phi\colon\A\rightarrow\B$ be a
  \starhomo.  
  The \defemph{mapping cone of $\phi$}, $\sfC_\phi$, is the
  pullback of the maps $\A\overset{\phi}{\longrightarrow}\B$ and
  $\sfC\B\overset{\ev_1}{\longrightarrow}\B$. 
As usual, the pullback can be realized as the restricted direct sum:
  $$\sfC_\phi=\A\oplus_{\phi,\ev_1}\sfC\B=\sets{(x,y)\in\A\oplus\sfC\B}{\phi(x)=\ev_1(y)=y(1)}.$$
\end{definition}

\begin{remark}
  Let $\phi\colon\A\rightarrow\B$ be a \starhomo between \cas.  Then
  there is a canonical short exact sequence
  $$\extwithoutmaps{\sfS\B}{\sfC_\phi}{\A}$$
  called the \defemph{mapping cone sequence}.  This sequence is
  natural in \A and \B, \ie, if there exists a commuting diagram
  $$\xymatrix{\A_1\ar[r]^{\phi_1}\ar[d]_f&\B_1\ar[d]^g\\
    \A_2\ar[r]^{\phi_2}&\B_2}$$
  then there is a (canonical) \starhomo
  $\omega\colon \sfC_{\phi_1}\rightarrow \sfC_{\phi_2}$ making the
  diagram
  $$\xymatrix{0\ar[r]&\sfS\B_1\ar[d]_{\sfS g}\ar[r]&\sfC_{\phi_1}\ar[d]^\omega\ar[r]&\A_1\ar[d]^f\ar[r]&0\\
    0\ar[r]&\sfS\B_2\ar[r]&\sfC_{\phi_2}\ar[r]&\A_2\ar[r]&0}$$
  commutative (\cf\ \cite[Section~19.4]{MR99g:46104}).  Actually, a concrete description of $\omega$ is as follows:
  $\omega(a,h)=(f(a),g\circ h)$ for all
  $(a,h)\in\A_1\oplus_{\phi_1,\ev_1}\sfC\B_1=\sfC_{\phi_1}$.
\end{remark}

\begin{remark}
  \label{rem-example19.4.2(a)inBlackadar}
  The mapping cone sequence of the
  identity homomorphism $\ident_\A$ is the canonical sequence 
  $\extwithoutmaps{\sfS\A}{\sfC\A}{\A}$.
  For each \starhomo $\phi\colon\A\rightarrow\B$ between \cas, 
  there exists a canonical \starisos $\Sflip$ from $\sfS\sfC_\phi$ to $\sfC_{\sfS\phi}$,
  and $\Cflip$ from $\sfC\sfC_\phi$ to $\sfC_{\sfC\phi}$, given by  
  $$\sfS\sfC_\phi
  =\sfS(\A\oplus_{\phi,\ev_1}\sfC\B)\ni(x,y)\mapsto
  (x,\flip(y))\in
  \sfS\A\oplus_{\sfS\phi,\ev_1}\sfC\sfS\B=\sfC_{\sfS\phi},$$
  $$\sfC\sfC_\phi
  =\sfC(\A\oplus_{\phi,\ev_1}\sfC\B)\ni(x,y)\mapsto
  (x,\flip(y))\in
  \sfC\A\oplus_{\sfC\phi,\ev_1}\sfC\sfC\B=\sfC_{\sfC\phi},$$
  respectively. 
  See Definition~\ref{def-commutatorofMCwithSandC} and
  Lemma~\ref{lem-commutatorofMCwithSandC} for more on these
  isomorphisms.  
\end{remark}

\begin{definition}\label{def-functorsonmathcalE}
  Define functors $\mc$, $\sfS$ and $\sfC$ on
  the category of all
  extensions of \cas (with the morphisms being triples of \starhomos
  making the obvious diagram commutative) 
  as follows.   
  For an extension 
  $e\colon\extwithmaps{\A_0}{\A_1}{\A_2}{\iota}{\pi}$
   set 
  $$\mc(e)\colon\extwithmaps{\sfS\A_2}{\sfC_\pi}{\A_1}{\iota_\mc}{\pi_\mc},$$
  $$\sfS(e)=\sfS e\colon
  \extwithmaps{\sfS\A_0}{\sfS\A_1}{\sfS\A_2}{\sfS\iota}{\sfS\pi},$$ 
  $$\sfC(e)=\sfC e\colon
  \extwithmaps{\sfC\A_0}{\sfC\A_1}{\sfC\A_2}{\sfC\iota}{\sfC\pi}.$$ 
  For a morphism $\phi=(\phi_0,\phi_1,\phi_2)$ from $e$ to $e'$, 
  let $\mc(\phi)$ be the morphism $(\sfS\phi_2,\omega,\phi_1)$ defined
  using the naturality of the mapping cone construction (see
  above), 
  let $\sfS(\phi)=\sfS\phi$ be the morphism
  $(\sfS\phi_0,\sfS\phi_1,\sfS\phi_2)$, and 
  we let $\sfC(\phi)=\sfC\phi$ be the morphism
  $(\sfC\phi_0,\sfC\phi_1,\sfC\phi_2)$.   

  It is easy to verify that these are functors. 
  Moreover, one easily verifies, that they are exact 
  (\ie, they send short exact 
  sequences of extensions to short exact sequences of extensions). 
\end{definition}

\begin{definition}
  Let $e\colon\extwithmaps{\A_0}{\A_1}{\A_2}{\iota}{\pi}$ be an extension of \cas. 
  Then construct two new extensions, 
  $\mathfrak{i}(e)$ and $\mathfrak{q}(e)$, from $e$ as follows. 
  Let $\mathfrak{i}(e)$ denote the extension 
  $\A_0=\A_0\twoheadrightarrow 0$, and let 
  $\mathfrak{q}(e)$ denote the extension 
  $0\hookrightarrow\A_2=\A_2$. 
  Then there exists a canonical short exact sequence of extensions: 
  $\extwithmaps{\mathfrak{i}(e)}{e}{\mathfrak{q}(e)}{\mathfrak{i}_e}{\mathfrak{q}_e}$. 
\end{definition}

\begin{remark}
  Note that if $e\colon\extwithmaps{\A_0}{\A_1}{\A_2}{\iota}{\pi}$
  is an extension of \cas, then
  there exists a commuting diagram
  $$\xymatrix{
    0\ar@{^(->}[r]\ar@{^(->}[d]&\sfS\A_2\ar@{=}[r]\ar@{^(->}[d]^{\iota_{\mc}}&\sfS\A_2\ar@{^(->}[d]\\     \A_0\ar@{^(->}[r]^{f_e}\ar@{=}[d]&\sfC_\pi\ar@{->>}[r]\ar@{->>}[d]^{\pi_\mc}&\sfC\A_2\ar@{->>}[d]^{\ev_1}\\ 
    \A_0\ar@{^(->}[r]^\iota&\A_1\ar[r]^\pi&\A_2    }$$
  with short exact rows and columns.  The map
  $f_e\colon\A_0\rightarrow \sfC_\pi$ induces isomorphism on the level
  of $K$-theory (actually, this holds more generally for additive,
  homotopy-invariant, half-exact functors, 
  \cf\ \cite[Proposition~21.4.1]{MR99g:46104}).  Actually, this diagram is nothing but the short exact sequence
  $\extwithmaps{\mc(\mathfrak{i}(e))}{\mc(e)}{\mc(\mathfrak{q}(e))}{\mc(\mathfrak{i}_e)}{\mc(\mathfrak{q}_e)}$
  induced by applying the functor $\mc$ to the short exact sequence
  $\extwithmaps{\mathfrak{i}(e)}{e}{\mathfrak{q}(e)}{\mathfrak{i}_e}{\mathfrak{q}_e}$.
\end{remark}

\begin{definition}\label{def-commutatorofMCwithSandC}
  Let there be given an extension 
  $e\colon\extwithmaps{\A_0}{\A_1}{\A_2}{\iota}{\pi}$ of \cas. 
  Form the extensions $\sfS(\mc(e))$, $\mc(\sfS(e))$, $\sfC(\mc(e))$,
  and $\mc(\sfC(e))$ as above. 
  Then define morphisms $\theta_e$ 
  from $\sfS(\mc(e))$ to $\mc(\sfS e)$  
  and $\eta_e$ 
  from $\sfC(\mc(e))$ to $\mc(\sfC e)$ as follows:  
  $$\xymatrix{
    \sfS(\mc(e))\colon\ar[d]_\cong^{\theta_e} & 
    0\ar[r]&
    \sfS\sfS\A_2\ar[r]^{\sfS(\iota_\mc)}\ar[d]^{\flip}_\cong&
    \sfS\sfC_\pi\ar[r]^{\sfS(\iota_\mc)}\ar[d]^{\Sflip}_\cong&
    \sfS\A_1\ar[r]\ar@{=}[d]&0\\
    \mc(\sfS e)\colon &
    0\ar[r]&
    \sfS\sfS\A_2\ar[r]^{(\sfS\iota)_\mc}&
    \sfC_{\sfS\pi}\ar[r]^{(\sfS\pi)_\mc}&
    \sfS\A_1\ar[r]&0,
    }$$
  $$\xymatrix{
    \sfC(\mc(e))\colon\ar[d]_\cong^{\eta_e} &
    0\ar[r]&
    \sfC\sfS\A_2\ar[r]^{\sfC(\iota_\mc)}\ar[d]^{\flip}_\cong&
    \sfC\sfC_\pi\ar[r]^{\sfC(\iota_\mc)}\ar[d]^{\Cflip}_\cong&
    \sfC\A_1\ar[r]\ar@{=}[d]&0\\
    \mc(\sfC e)\colon &
    0\ar[r]&
    \sfS\sfC\A_2\ar[r]^{(\sfC\iota)_\mc}&
    \sfC_{\sfC\pi}\ar[r]^{(\sfC\pi)_\mc}&
    \sfC\A_1\ar[r]&0,
    }$$
  where the \starhomos 
  $\sfS\sfC_\pi\rightarrow\sfC_{\sfS\pi}$ and 
  $\sfC\sfC_\pi\rightarrow\sfC_{\sfC\pi}$ are the
  canonical isomorphisms from Remark~\ref{rem-example19.4.2(a)inBlackadar}. 
\end{definition}

\begin{lemma}\label{lem-commutatorofMCwithSandC}
  The above morphisms, $\theta_e$ and $\eta_e$, are functorial, \ie, they implement 
  isomorphisms from the functor $\sfS\circ\mc$ to the functor
  $\mc\circ\sfS$ and from the functor $\sfC\circ\mc$ to the functor
  $\mc\circ\sfC$, respectively. 
\end{lemma}
\begin{proof}
  This is a long, straightforward verification. 
\end{proof}

\begin{lemma}
  Let $e$ be an extension of \cas. 
  Then there exists an isomorphism of short exact sequences of extensions as follows: 
  $$\xymatrix{
    0\ar[r] &
    \sfS\mc(e)\ar[d]^{\theta_e}\ar[r] &
    \sfC\mc(e)\ar[d]^{\eta_e}\ar[r] &
    \mc(e)\ar@{=}[d]\ar[r] &
    0 \\
    0\ar[r] &
    \mc(\sfS e)\ar[r] &
    \mc(\sfC e)\ar[r] &
    \mc(e)\ar[r] &
    0. 
    }$$
\end{lemma}
\begin{proof}
  This is a straightforward verification. 
\end{proof}

\begin{lemma}\label{lem-helpforthelemmawithbigbigbigdiagram}
  Let there be given a commutative diagram 
  $$\xymatrix{
    \mathfrak{X}\ar[r]^{\phi_1}\ar[d]_{\phi_2} &
    \mathfrak{Y}_1\ar[d]^{\psi_1} \\
    \mathfrak{Y}_2\ar[r]_{\psi_2} &
    \mathfrak{Z}
    }$$
  of \cas and \starhomos.  Then there are canonically induced \starhomos 
  $\sfC_{\phi_1}\rightarrow\sfC_{\psi_2}$ and 
  $\sfC_{\phi_2}\rightarrow\sfC_{\psi_1}$.  
  The mapping cones $\sfC_{\sfC_{\phi_1}\rightarrow\sfC_{\psi_2}}$ and 
  $\sfC_{\sfC_{\phi_2}\rightarrow\sfC_{\psi_1}}$ are canonically
  isomorphic to 
  $$\sets{(x,f_1,f_2,h)\in\mathfrak{X}\oplus\sfC\mathfrak{Y}_1\oplus\sfC\mathfrak{Y}_2\oplus\sfC\sfC\mathfrak{Z}}%
  {\begin{matrix}
    \phi_1(x)=f_1(1),&\psi_1\circ f_1(-)=h(1,-),\\
    \phi_2(x)=f_2(1),&\psi_2\circ f_2(-)=h(-,1)
    \end{matrix}
    }$$
  $$\sets{(x,f_2,f_1,h)\in\mathfrak{X}\oplus\sfC\mathfrak{Y}_2\oplus\sfC\mathfrak{Y}_1\oplus\sfC\sfC\mathfrak{Z}}%
  {\begin{matrix}
    \phi_1(x)=f_1(1),&\psi_1\circ f_1(-)=h(-,1),\\
    \phi_2(x)=f_2(1),&\psi_2\circ f_2(-)=h(1,-)
    \end{matrix}
    }$$
  respectively. 
  So $(x,f_1,f_2,h)\mapsto(x,f_2,f_1,\flip(h))$ is an isomorphism from
  $\sfC_{\sfC_{\phi_1}\rightarrow\sfC_{\psi_2}}$ to 
  $\sfC_{\sfC_{\phi_2}\rightarrow\sfC_{\psi_1}}$. 
\end{lemma}
\begin{proof}
  This is straightforward to check by writing out the mapping cones as
  restricted direct sums. 
  Note that we only need to check the first statement, since the
  second follows by symmetry (by interchanging $1$ and $2$). 
\end{proof}

For a morphism $\phi=(\phi_0,\phi_1,\phi_2)$ between extensions of
\cas, let $\sfC_\phi$ denote the object 
$\extwithoutmaps{\sfC_{\phi_0}}{\sfC_{\phi_1}}{\sfC_{\phi_2}}$ 
(\cf\ also \bonkat[Definition~3.4.1]). 

\begin{lemma}\label{lem-bigbigbig-commutingdiagram}
  Let there be given a commuting diagram
  $$\xymatrix{
    \A_0\ar@{^(->}[d]^{\iota_\A}\ar@{^(->}[r]^{x_0}&
    \B_0\ar@{^(->}[d]^{\iota_\B}\ar@{->>}[r]^{y_0}&
    \CC_0\ar@{^(->}[d]^{\iota_\CC}\\ 
    \A_1\ar@{->>}[d]^{\pi_\A}\ar@{^(->}[r]^{x_1}&
    \B_1\ar@{->>}[d]^{\pi_\B}\ar@{->>}[r]^{y_1}&
    \CC_1\ar@{->>}[d]^{\pi_\CC}\\ 
    \A_2\ar@{^(->}[r]^{x_2}&
    \B_2\ar@{->>}[r]^{y_2}&
    \CC_2 
    }$$
  with the rows and columns being short exact sequences
  of \cas{} 
  --- we will write this short as
  $\extwithmaps{e_\A}{e_\B}{e_\CC}{\mathbf{x}}{\mathbf{y}}$.   
  Then there exists an isomorphism $\xi_\mathbf{y}$ from 
  $\sfC_{\mc(\mathbf{y})}$ to $\mc(\sfC_\mathbf{y})$ given as follows: 
  $$\xymatrix{
    \sfC_{\mc(\mathbf{y})}\colon\ar[d]_\cong^{\xi_{\mathbf{y}}} &
    0\ar[r] &
    \sfC_{\sfS y_2}\ar_\cong^{\Sflip}[d]\ar[r] &
    \sfC_{\sfC_{\pi_\B}\rightarrow\sfC_{\pi_\CC}}\ar[d]_\cong\ar[r] &
    \sfC_{y_1}\ar@{=}[d]\ar[r] & 0 \\
    \mc(\sfC_\mathbf{y})\colon &
    0\ar[r] &
    \sfS\sfC_{y_2}\ar[r] &
    \sfC_{\sfC_{y_1}\rightarrow\sfC_{y_2}}\ar[r] &
    \sfC_{y_1}\ar[r] & 0 
    }$$
  where the isomorphism from  
  $\sfC_{\sfC_{\pi_\B}\rightarrow\sfC_{\pi_\CC}}$ to 
  $\sfC_{\sfC_{y_1}\rightarrow\sfC_{y_2}}$ is given as in the above
  lemma.
  Moreover, the map given by the matrix  
  $$\begin{pmatrix}
    0      & \theta_{e_\CC}   & \theta_{e_\CC} \\
    \ident & \xi_{\mathbf{y}} & \eta_{e_\CC}\\
    \ident & \ident           & \ident 
  \end{pmatrix}$$
  between the standard diagrams 
  $$\vcenter{\xymatrix{
      0\ar@{^(->}[d]\ar@{^(->}[r] &
      \sfS\mc(e_\CC)\ar@{^(->}[d]\ar@{=}[r] &
      \sfS\mc(e_\CC)\ar@{^(->}[d] \\
      \mc(e_\A)\ar@{=}[d]\ar@{^(->}[r] &
      \sfC_{\mc(\mathbf{y})}\ar@{->>}[d]\ar@{->>}[r] &
      \sfC\mc(e_\CC)\ar@{->>}[d] \\
      \mc(e_\A)\ar@{^(->}[r] &
      \mc(e_\B)\ar@{->>}[r] &
      \mc(e_\CC)
      }}
  \quad and\quad 
  \vcenter{\xymatrix{
      0\ar@{^(->}[d]\ar@{^(->}[r] &
      \mc(\sfS e_\CC)\ar@{^(->}[d]\ar@{=}[r] &
      \mc(\sfS e_\CC)\ar@{^(->}[d] \\
      \mc(e_\A)\ar@{=}[d]\ar@{^(->}[r] &
      \mc(\sfC_{\mathbf{y}})\ar@{->>}[d]\ar@{->>}[r] &
      \mc(\sfC e_\CC)\ar@{->>}[d] \\
      \mc(e_\A)\ar@{^(->}[r] &
      \mc(e_\B)\ar@{->>}[r] &
      \mc(e_\CC)
      }}$$
  makes everything commutative. 
\end{lemma}
\begin{proof}
  Using the above, we have that
  $\sfC_{\sfC_{\pi_\B}\rightarrow\sfC_{\pi_\CC}}$ is isomorphic to 
  $$\sets{(x,f_2,f_1,h)\in\B_1\oplus\sfC\B_2\oplus\sfC\CC_1\oplus\sfC\sfC\CC_2}%
  {\begin{matrix}
    y_1(x)=f_1(1),&\pi_\CC\circ f_1(-)=h(-,1),\\
    \pi_\B(x)=f_2(1),&y_2\circ f_2(-)=h(1,-)
    \end{matrix}
    }$$
  and $\sfC_{\sfC_{y_1}\rightarrow\sfC_{y_2}}$ is isomorphic to 
  $$\sets{(x,f_1,f_2,h)\in\B_1\oplus\sfC\CC_1\oplus\sfC\B_2\oplus\sfC\sfC\CC_2}%
  {\begin{matrix}
    y_1(x)=f_1(1),&\pi_\CC\circ f_1(-)=h(1,-),\\
    \pi_\B(x)=f_2(1),&y_2\circ f_2(-)=h(-,1)
    \end{matrix}
    }$$
  and, moreover, 
  \begin{align*}
    \sfC_{\sfS y_2}
    &=\sets{(f_2,h)\in\sfS\B_2\oplus\sfC\sfS\CC_2}{y_2\circ f_2(-)=h(1,-)}, \\
    \sfS\sfC_{y_2}
    &=\sets{(f_2,h)\in\sfS\B_2\oplus\sfS\sfC\CC_2}{y_2\circ f_2(-)=h(-,1)}, \\
    \sfC_{y_1}
    &=\sets{(x,f_1)\in\B_1\oplus\sfC\CC_1}{y_1(x)=f_1(1)}. 
  \end{align*}
  Using these identifications, we compute the extensions: 
  $$\xymatrix@R=12pt{
    \sfC_{\mc(\mathbf{y})}\colon&
    0\ar[r] &
    \sfC_{\sfS y_2}\ar[rr]^-{(f_2,h)\mapsto(0,f_2,0,h)} &&
    \sfC_{\sfC_{\pi_\B}\rightarrow\sfC_{\pi_\CC}}\ar[rr]^-{(x,f_2,f_1,h)\mapsto(x,f_1)} &&
    \sfC_{y_1}\ar[r] & 0, \\
    \mc(\sfC_{\mathbf{y}})\colon&
    0\ar[r] &
    \sfS\sfC_{y_2}\ar[rr]^-{(f_2,h)\mapsto(0,0,f_2,h)} &&
    \sfC_{\sfC_{y_1}\rightarrow\sfC_{y_2}}\ar[rr]^-{(x,f_1,f_2,h)\mapsto(x,f_1)} &&
    \sfC_{y_1}\ar[r] & 0. 
    }$$
  Now it is routine to check that the given diagram commutes. 

  Second part: 
  The above results show that every square which does not involve
  $\sfC_{\mc(\mathbf{y})}$ and $\mc(\sfC_{\mathbf{y}})$ commutes. 
  The long and straightforward proof of the commutativity of the
  remaining four squares of morphisms of extension is left to the
  reader. 
\end{proof}

\section{Homology and cohomology theories for \cas}
\label{sec-homologyandcohomologytheories}

We recall some of the definitions and results about homology and
cohomology theories on \cas. 
The reader is referred to \cite[Chapters~21 and~22]{MR99g:46104} 
(these two chapters are
primarily due to Cuntz, Higson, Rosenberg, and Schochet --- see the
monograph for further references). 

\begin{definition}
  Let $\mathgrc{S}$ be a subcategory of the category of all \cas,
  which is closed under quotients, extensions, and closed under
  suspension in the sense that if $\A$ is an object of $\mathgrc{S}$,
  then the suspension $\sfS\A$ of $\A$ is also an object of
  $\mathgrc{S}$, $\sfS\phi$ is a morphism in $\mathgrc{S}$ whenever
  $\phi$ is, $\sfS\C$ is an object of $\mathgrc{S}$ and every
  \starhomo from $\sfS\C$ to every object of $\mathgrc{S}$ is a
  morphism in $\mathgrc{S}$.
  
  Let $\mathgrc{Ab}$ denote the category of abelian groups.  We will
  consider functors $\functorF$ from $\mathgrc{S}$ to $\mathgrc{Ab}$.
  Such functors may or may not satisfy each of the axioms 
  \defemph{Homotopy-invariance} (H), 
  \defemph{Stability} (S), 
  \defemph{$\sigma$-additivity} (A),
  \defemph{completely additive}, 
  \defemph{additive},
  \defemph{Half-exactness} (HX)
  defined in 
  \cite[Chapters~21 and~22]{MR99g:46104}. 

  In \cite[Definition~21.1.1]{MR99g:46104} 
  a \defemph{homology theory} on $\mathgrc{S}$ is defined as a
  sequence $(h_n)$ of
  covariant functors from $\mathgrc{S}$ to $\mathgrc{Ab}$ satisfying
  (H) and having natural long exact sequences associated with short
  exact sequences.  A \defemph{cohomology theory} is defined analogous. 
\end{definition}

\begin{definition}\label{def-longexactsequence}
  Let $\functorF$ be an additive functor from $\mathgrc{S}$ to
  $\mathgrc{Ab}$ satisfying (H) and (HX). 

  Let $e\colon\extwithmaps{\A_0}{\A_1}{\A_2}{\iota}{\pi}$ be a given
  extension. 
  Then for each $n\in\N_0$, 
  set 
  $$\functorF_n=\functorF\circ\sfS^n, \quad\text{and}\quad
  \partial_{n+1}=\functorF_{n}(f_e)^{-1}\circ\functorF_{n}(\iota_\mc)
  \colon\functorF_{n+1}(\A_2)\rightarrow\functorF_{n}(\A_0), 
  \text{ if }\functorF\text{ is covariant,}$$
  $$\functorF^n=\functorF\circ\sfS^n,\quad\text{and}\quad
  \partial^{n}=\functorF^{n}(\iota_\mc)\circ\functorF^{n}(f_e)^{-1}
  \colon\functorF^{n}(\A_0)\rightarrow\functorF^{n+1}(\A_2),
  \text{ if }\functorF\text{ is contravariant,}$$
  where 
  $\iota_\mc\colon\sfS\A_2\rightarrow\sfC_\pi$ 
  and $f_e\colon\A_0\rightarrow\sfC_\pi$ are the canonical
  \starhomos.\footnote{Note that $\sfS^n$ denotes the composition of
    $\sfS$ with itself $n$ times, while the superscript in $\functorF^n$
    indicates that this is some kind of $n$'th cohomology.}  
\end{definition}

From \cite[Theorem~21.4.3]{MR99g:46104} we have the following
theorem. 

\begin{theorem}
  Let $\functorF$ be an additive functor from $\mathgrc{S}$ to
  $\mathgrc{Ab}$ satisfying (H) and (HX). 
  If $\functorF$ is covariant, then
  $(\functorF_n)_{n=0}^\infty$ is a homology theory. 
  If $\functorF$ is contravariant, then
  $(\functorF^n)_{n=0}^\infty$ is a cohomology theory. 
\end{theorem}

\begin{corollary}\label{cor-splitexactness}
  If $\functorF$ is an additive functor from $\mathgrc{S}$ to
  $\mathgrc{Ab}$ satisfying (H) and (HX), then $\functorF$ is
  split-exact, \ie, $\functorF$ sends split-exact sequences from
  $\mathgrc{S}$ to split-exact sequences of abelian groups.  
\end{corollary}

\begin{proof}
  Let $\extwithmaps{\A_0}{\A_1}{\A_2}{\iota}{\pi}$ be a split-exact
  sequence of \cas, and assume that $\functorF$ is covariant. 
  It is clear that $\functorF(\pi)$ and $\functorF(\sfS\pi)$ are
  surjective (since $\functorF$ and $\functorF\circ\sfS$ are
  functors). 
  From preceding theorem it follows that $\partial_1=0$, so
  $\functorF\iota$ is injective. 
  The proof in the contravariant case is dual. 
\end{proof}

The following theorem is taken from
\cite[Corollary~22.3.2]{MR99g:46104}. 

\begin{theorem}
  Let $\functorF$ be an additive functor from $\mathgrc{S}$ 
  to $\mathgrc{Ab}$ satisfying 
  (H), (S), and (HX). 
  Then $\functorF$ is naturally isomorphic to $\functorF\circ\sfS^2$. 
\end{theorem}

\begin{definition}\label{def-standardcyclicsixtermexactsequence}
  Let $\functorF$ be an additive functor from $\mathgrc{S}$ to
  $\mathgrc{Ab}$ satisfying (H), (S), and (HX), and let 
  $\beta_\A\colon\functorF(\A)\rightarrow\functorF(\sfS^2\A)$ 
  denote the natural isomorphism.  
  Then for each short exact sequence
  $$e\colon\extwithmaps{\A_0}{\A_1}{\A_2}{\iota}{\pi}$$
  of \cas we make the following definition. 
  If $\functorF$ is covariant, then 
  define  
  $\partial_0\colon\functorF(\A_2)\rightarrow\functorF(\sfS\A_0)$ as
  the composition of the homomorphisms
  $$\xymatrix{ \functorF(\A_2)\ar[r]^{\beta_{\A_2}} &
    \functorF(\sfS^2\A_2)\ar[r]^{\partial_2} & \functorF(\sfS\A_0).
    }$$
  If $\functorF$ is contravariant, then 
  define $\widetilde{\partial^{1}}\colon\functorF(\sfS\A_0)\rightarrow\functorF(\A_2)$
  as the composition of the homomorphisms
  $$\xymatrix{ \functorF(\sfS\A_0)\ar[r]^{\partial^1} &
    \functorF(\sfS^2\A_2)\ar[r]^{\beta_{\A_2}^{-1}} & \functorF(\A_2).
    }$$
  So with each such short exact sequence we have associated a cyclic
  six term exact sequence
  $$\vcenter{\xymatrix@R16pt{ \functorF(\A_0)\ar[r]^{\functorF(\iota)} &
      \functorF(\A_1)\ar[r]^{\functorF(\pi)}
      & \functorF(\A_2)\ar[d]^{\partial_0} \\
      \functorF(\sfS\A_2)\ar[u]^{\partial_1} &
      \functorF(\sfS\A_1)\ar[l]^{\functorF(\sfS\pi)} &
      \functorF(\sfS\A_0)\ar[l]^{\functorF(\sfS\iota)}
      }}\quad\text{respectively}\quad 
  \vcenter{\xymatrix@R16pt{
      \functorF(\A_2)\ar[r]^{\functorF(\pi)} &
      \functorF(\A_1)\ar[r]^{\functorF(\iota)}
      & \functorF(\A_0)\ar[d]^{\partial^0} \\
      \functorF(\sfS\A_0)\ar[u]^{\widetilde{\partial^{1}}} &
      \functorF(\sfS\A_1)\ar[l]^{\functorF(\sfS\iota)} &
      \functorF(\sfS\A_2)\ar[l]^{\functorF(\sfS\pi)} }}$$
  which is
  natural with respect to morphisms of short exact sequences of \cas.
  We will occasionally misuse the notation and write $\partial^1$
  instead of $\widetilde{\partial^{1}}$ (which should not cause any
  confusions). 
\end{definition}

\begin{remark}
  While it is obvious how to generalize 
  homotopy-invariance, stability, additivity, and
  split-exactness for a functor from $\mathgrc{S}$ to an additive
  category $\mathgrc{A}$, it is not obvious how to generalize
  half-exactness.  
  
  In Section~21.4 in \cite{MR99g:46104}, Blackadar defines
  half-exactness for such functors 
  (\ie, $\operatorname{Hom}_{\mathgrc{A}}(X,\functorF(-))$ and
  $\operatorname{Hom}_{\mathgrc{A}}(\functorF(-),X)$ should be
  half-exact for all objects $X$). 
  It is natural to ask whether this \emph{extends} the original
  definition, and the answer is no. 
  This is seen by applying $\operatorname{Hom}_\Z(\Z_3,K_1(-))$ to the
  short exact sequence $\extwithoutmaps{\sfS\Matn[3]}{\Io[3]}{\C}$ 
  (\cf\ Definition~\ref{def-dimensiondropinterval}). 
  On the other hand, for the category of modules over a unital ring,
  $\operatorname{Hom}_R(R,M)$ is naturally isomorphic to $M$ ---
  so this property is stronger than the ordinary
  half-exactness. 
  To avoid confusions, we will not use this terminology. 
\end{remark}

\section{(Co-)Homology theories and mapping cone sequences}
\label{sec-homologyandmappingcones}

In this section we show exactly how the cyclic six term exact sequence
of the mapping cone sequence for an extension of \cas is related to the
cyclic six term exact sequence of the original extension.  First we will need the following lemma, which Bonkat uses a version of
in the proof of \bonkat[Lemma~7.3.1].  The proof given here is much more elementary. 

\begin{lemma}
  \label{lem-anticommute-kk}
  Let $\functorF_0$ and $\functorF_1$ be covariant additive functors from the category
  $\mathgrc{S}$ to the category $\mathgrc{Ab}$, which have the
  properties (H), (S), and (HX). 
  Assume that $\partial_0^-$ and $\partial_1^-$ are boundary maps
  making $(\functorF_i,\partial_i)_{i=0}^1$ into a cyclic homology
  theory on $\mathgrc{S}$.  Let there also be given a commuting diagram
  $$\xymatrix{
    \A_0\ar@{^(->}[d]\ar@{^(->}[r]&\A_1\ar@{^(->}[d]\ar@{->>}[r]&\A_2\ar@{^(->}[d]\\ 
    \B_0\ar@{->>}[d]\ar@{^(->}[r]&\B_1\ar@{->>}[d]\ar@{->>}[r]&\B_2\ar@{->>}[d]\\ 
    \CC_0\ar@{^(->}[r]&\CC_1\ar@{->>}[r]&\CC_2 
    }$$
  with the rows and columns being short exact sequences
  of \cas.  
  Let $e_\A,e_\B$ and $e_\CC$ denote the three horizontal
  extensions, while $e_0,e_1$ and $e_2$ denote the three vertical
  extensions.  
  Then there exists a diagram
  $$\def\objectstyle{\scriptstyle} 
  \def\labelstyle{\scriptstyle}
  \xymatrix@R=18pt{
    &\ar[d]^-{\partial_1^{e_0}}&\ar[d]^-{\partial_1^{e_1}}&\ar[d]^-{\partial_1^{e_2}}&
    \ar[d]^-{\partial_0^{e_0}}&\ar[d]^-{\partial_0^{e_1}}&\ar[d]^-{\partial_0^{e_2}}\\
    \ar[r]^-{\partial_1^{e_{\A}}}&
    \functorF_0(\A_0)\ar[d]\ar[r]&
    \functorF_0(\A_1)\ar[d]\ar[r]&
    \functorF_0(\A_2)\ar[d]\ar[r]^-{\partial_0^{e_{\A}}}&
    \functorF_1(\A_0)\ar[d]\ar[r]&
    \functorF_1(\A_1)\ar[d]\ar[r]&
    \functorF_1(\A_2)\ar[d]\ar[r]^-{\partial_1^{e_{\A}}}& \\
    \ar[r]^-{\partial_1^{e_{\B}}}&
    \functorF_0(\B_0)\ar[d]\ar[r]&
    \functorF_0(\B_1)\ar[d]\ar[r]&
    \functorF_0(\B_2)\ar[d]\ar[r]^-{\partial_0^{e_{\B}}}&
    \functorF_1(\B_0)\ar[d]\ar[r]&
    \functorF_1(\B_1)\ar[d]\ar[r]&
    \functorF_1(\B_2)\ar[d]\ar[r]^-{\partial_1^{e_{\B}}}&\\
    \ar[r]^-{\partial_1^{e_{\CC}}}&
    \functorF_0(\CC_0)\ar[d]^-{\partial_0^{e_0}}\ar[r]&
    \functorF_0(\CC_1)\ar[d]^-{\partial_0^{e_1}}\ar[r]&
    \functorF_0(\CC_2)\ar[d]^-{\partial_0^{e_2}}\ar[r]^-{\partial_0^{e_{\CC}}}&
    \functorF_1(\CC_0)\ar[d]^-{\partial_1^{e_0}}\ar[r]&
    \functorF_1(\CC_1)\ar[d]^-{\partial_1^{e_1}}\ar[r]&
    \functorF_1(\CC_2)\ar[d]^-{\partial_1^{e_2}}\ar[r]^-{\partial_1^{e_{\CC}}}&\\
    \ar[r]^-{\partial_0^{e_{\A}}}&
    \functorF_1(\A_0)\ar[d]\ar[r]&
    \functorF_1(\A_1)\ar[d]\ar[r]&
    \functorF_1(\A_2)\ar[d]\ar[r]^-{\partial_1^{e_{\A}}}&
    \functorF_0(\A_0)\ar[d]\ar[r]&
    \functorF_0(\A_1)\ar[d]\ar[r]&
    \functorF_0(\A_2)\ar[d]\ar[r]^-{\partial_0^{e_{\A}}}& \\
    \ar[r]^-{\partial_0^{e_{\B}}}&
    \functorF_1(\B_0)\ar[d]\ar[r]&
    \functorF_1(\B_1)\ar[d]\ar[r]&
    \functorF_1(\B_2)\ar[d]\ar[r]^-{\partial_1^{e_{\B}}}&
    \functorF_0(\B_0)\ar[d]\ar[r]&
    \functorF_0(\B_1)\ar[d]\ar[r]&
    \functorF_0(\B_2)\ar[d]\ar[r]^-{\partial_0^{e_{\B}}}&\\
    \ar[r]^-{\partial_0^{e_{\CC}}}&
    \functorF_1(\CC_0)\ar[d]^-{\partial_1^{e_0}}\ar[r]&
    \functorF_1(\CC_1)\ar[d]^-{\partial_1^{e_1}}\ar[r]&
    \functorF_1(\CC_2)\ar[d]^-{\partial_1^{e_2}}\ar[r]^-{\partial_1^{e_{\CC}}}&
    \functorF_0(\CC_0)\ar[d]^-{\partial_0^{e_0}}\ar[r]&
    \functorF_0(\CC_1)\ar[d]^-{\partial_0^{e_1}}\ar[r]&
    \functorF_0(\CC_2)\ar[d]^-{\partial_0^{e_2}}\ar[r]^-{\partial_0^{e_{\CC}}}&\\
    &&&&&&}$$
  with the cyclic six term exact sequence both horizontally and
  vertically.  
  The two squares
  $$\xymatrix{
    \functorF_0(\CC_2)\ar[d]^-{\partial_0^{e_2}}\ar[r]^-{\partial_0^{e_{\CC}}}&
    \functorF_1(\CC_0)\ar[d]^-{\partial_1^{e_0}}&
    \functorF_1(\CC_2)\ar[d]^-{\partial_1^{e_2}}\ar[r]^-{\partial_1^{e_{\CC}}}&
    \functorF_0(\CC_0)\ar[d]^-{\partial_0^{e_0}}\\
    \functorF_1(\A_2)\ar[r]^-{\partial_1^{e_{\A}}}& \functorF_0(\A_0)&
    \functorF_0(\A_2)\ar[r]^-{\partial_0^{e_{\A}}}& \functorF_1(\A_0)& }
  $$
  anticommute, while all the other squares (in the big diagram)
  commute.

  If $\functorF$ is contravariant instead, the dual statement holds. 
\end{lemma}
\begin{proof}
  That all the other squares commute, is evident (using that $\functorF_0$ and
  $\functorF_1$ are functors and that the maps $\partial_0$ and $\partial_1$ are
  natural).  Let $\D$ denote the pullback of $\CC_2$ along $\B_2\rightarrow\CC_2$
  and $\CC_1\rightarrow\CC_2$. 
  Then there exist short exact sequences 
  $$\xymatrix@R=12pt{e_{\textrm{sum}}:&\A_0\ar@{^(->}[r]&\A_1+\B_0\ar@{->>}[r]&\A_2\oplus\CC_0\\
    e_{\textrm{pullback}}:&\A_2\oplus\CC_0\ar@{^(->}[r]&\D\ar@{->>}[r]&\CC_2,}$$
  where we identify $\A_1$ and $\B_0$ with their images inside
  $\B_1$. 
  Split-exactness of $\functorF_0$ and $\functorF_1$, 
  \cf\ Corollary~\ref{cor-splitexactness}, and naturality of
  $\partial_0$ and $\partial_1$ together with the morphisms of extensions 
  $$\def\objectstyle{\scriptstyle} \def\labelstyle{\scriptstyle}
  \xymatrix{
    \A_0\ar@{=}[d]\ar@{^(->}[r]&\B_0\ar[d]\ar@{->>}[r]&\CC_0\ar[d]&
    \A_0\ar@{=}[d]\ar@{^(->}[r]&\A_1\ar[d]\ar@{->>}[r]&\A_2\ar[d]\\
    \A_0\ar@{^(->}[r]&\A_1+\B_0\ar@{->>}[r]&\A_2\oplus\CC_0&
    \A_0\ar@{^(->}[r]&\A_1+\B_0\ar@{->>}[r]&\A_2\oplus\CC_0\\
    \A_2\oplus\CC_0\ar[d]\ar@{^(->}[r]&\D\ar@{->>}[r]\ar[d]&\CC_2\ar@{=}[d]&
    \A_2\oplus\CC_0\ar[d]\ar@{^(->}[r]&\D\ar@{->>}[r]\ar[d]&\CC_2\ar@{=}[d]\\
    \A_2\ar@{^(->}[r]&\B_2\ar@{->>}[r]&\CC_2&
    \CC_0\ar@{^(->}[r]&\CC_1\ar@{->>}[r]&\CC_2
    }$$
  give that the map  
  $\partial_{1-j}^{e_{\textrm{sum}}}\partial_{j}^{e_{\textrm{pullback}}}\colon 
  \functorF_j(\CC_2)\rightarrow \functorF_j(\A_0)$ 
  is exactly
  $\partial_{1-j}^{e_\A}\partial_j^{e_2}+\partial_{1-j}^{e_0}\partial_j^{e_\CC}$, for $j=0,1$.  
  But it turns out that
  $\partial_{1-j}^{e_{\textrm{sum}}}\partial_{j}^{e_{\textrm{pullback}}}=0$ proving anticommutativity. 
  For we have the following commuting diagram with short exact rows and columns
  $$\xymatrix{
    \A_0\ar@{=}[d]\ar@{^(->}[r]&\A_1+\B_0\ar@{^(->}[d]\ar@{->>}[r]&\A_2\oplus\CC_0\ar@{^(->}[d]\\
    \A_0\ar@{->>}[d]\ar@{^(->}[r]&\B_1\ar@{->>}[d]\ar@{->>}[r]&\D\ar@{->>}[d]\\
    0\ar@{^(->}[r]&\CC_2\ar@{=}[r]&\CC_2,
    }$$
  so the map
  $\partial_j^{e_{\textrm{pullback}}}$ factors through 
  $\functorF_{1-j}(\A_1+\B_0)\rightarrow
  \functorF_{1-j}(\A_2\oplus\CC_0)$.

  The proof in the case that $\functorF$ is contravariant is dual. 
\end{proof}

\begin{lemma}\label{lem-mappingconesequence-new}
  Let $\functorF$ be an additive functor from the category
  $\mathgrc{S}$ to 
  the category $\mathgrc{Ab}$, which has the
  properties (H), (S), and (HX). 
  Let \A be an arbitrary \ca. 
  The \scstes{}\footnote{as defined in 
    Definitions~\ref{def-longexactsequence} 
    and~\ref{def-standardcyclicsixtermexactsequence}} 
  associated with 
  $\extwithoutmaps{\sfS\A}{\sfC\A}{\A}$ is the
  sequence 
  $$\xymatrix{\functorF(\sfS\A)\ar[r]&0\ar[r]&\functorF(\A)\ar[d]_\cong^{-\beta_\A}\\
    \functorF(\sfS\A)\ar[u]_{-\ident}&0\ar[l]&\functorF(\sfS^2\A),\ar[l]}$$
  in the covariant case, and the sequence 
  $$\xymatrix{\functorF(\A)\ar[r]&0\ar[r]&\functorF(\sfS\A)\ar[d]^{-\ident}\\
    \functorF(\sfS^2\A)\ar[u]_{-\beta_\A^{-1}}^\cong&0\ar[l]&\functorF(\sfS\A),\ar[l]}$$
  in the contravariant case. 
\end{lemma}
\begin{proof}
  Assume that $\functorF$ is covariant. 
  Since the cone, $\sfC\A$, of $\A$ is homotopy equivalent to the zero
  \ca, $\functorF(\sfC\A)\cong\functorF(\sfS\sfC\A)\cong 0$ 
  (\cf\ \cite[Example~4.1.5]{MR2001g:46001}).  We have the commutative diagram
  $$\xymatrix{
    0\ar@{^(->}[r]\ar@{^(->}[d]&\sfS\A\ar@{=}[r]\ar@{^(->}[d]^{\iota_\mc}&\sfS\A\ar@{^(->}[d]\\
    \sfS\A\ar@{^(->}[r]^{f}\ar@{=}[d]&\sfC_\pi\ar@{->>}[r]\ar@{->>}[d]^{\pi_\mc}&\sfC\A\ar@{->>}[d]\\
    \sfS\A\ar@{^(->}[r]^\iota&\sfC\A\ar@{->>}[r]^\pi&\A
    }$$
  with short exact rows and columns. 
  Note that $\sfC_\pi$ is realized as
  $\{(x,y)\in\sfC\A\oplus\sfC\A\,|\,x(1)=y(1)\}$. 
  Using this picture there exists a 
  \starhomo $\varphi\colon\sfC\A\ni x\mapsto (x,x)\in\sfC_\pi$. 
  Note that the composed \starhomo $\varphi\circ\iota$ is just
  $f+\iota_\mc$. 
  Since $\functorF(\sfC\A)=0$, we must have
  $\functorF(\varphi\circ\iota)=0$. 
  Using the split-exactness of $\functorF$ 
  (\cf\ Corollary~\ref{cor-splitexactness}), 
  we get a canonical
  identification of $\functorF(\sfS\A\oplus\sfS\A)$ with
  $\functorF(\sfS\A)\oplus\functorF(\sfS\A)$. 
  Under this identification, we get
  $$\xymatrix{
    \functorF(\sfS\A)\ar[rr]\ar[drr]|-{x\mapsto(x,x)}
    &&\functorF(\sfS\A\oplus\sfS\A)\ar[d]_\cong\ar[rr]
    &&\functorF(\sfC_\pi)\\
    &&\functorF(\sfS\A)\oplus\functorF(\sfS\A)\ar[urr]|-{(x,y)\mapsto\functorF(f)(x)+\functorF(\iota_\mc)(y)}}$$
  Consequently,
  $$\functorF(f)+\functorF(\iota_\mc)
  =\functorF(f+\iota_\mc)=\functorF(\varphi\circ\iota)=0,$$ 
  and hence 
  $\functorF(f)=-\functorF(\iota_\mc)$.
  Therefore, we have
  $\partial_1=\functorF(f)^{-1}\circ\functorF(\iota_\mc)=-\ident$.

  The map
  $\partial_0\colon\functorF(\A)\rightarrow\functorF(\sfS^2\A)$ 
  is the composition of the maps 
  $$\xymatrix{
    \functorF(\A)\ar[r]^-{\beta_{\A}}
    & \functorF(\sfS^2\A)\ar[r]^-{\partial_2}
    & \functorF(\sfS^2\A),
    }$$
  where 
  $\partial_2=\functorF(\sfS f)^{-1}\circ\functorF(\sfS\iota_{\mc})$. 
  It is easy to see that the matrix 
  $$\begin{pmatrix}
    0 & \flip & \flip \\
    \flip & (\flip,\flip) & \flip \\
    \flip & \flip & \ident
  \end{pmatrix}$$
  implements a map between the diagrams 
  $$\vcenter{\xymatrix{
      0\ar@{^(->}[r]\ar@{^(->}[d] &
      \sfS\sfS\A\ar@{=}[r]\ar@{^(->}[d]^{\sfS\iota_\mc} &
      \sfS\sfS\A\ar@{^(->}[d] \\
      \sfS\sfS\A\ar@{=}[d]\ar@{^(->}[r]^{\sfS f} &
      \sfS\sfC_{\pi}\ar@{->>}[r]\ar@{->>}[d]^{\sfS\pi_\mc} &
      \sfS\sfC\A\ar@{->>}[d] \\
      \sfS\sfS\A\ar@{^(->}[r]^{\sfS\iota} &
      \sfS\sfC\A\ar@{->>}[r]^{\sfS\pi} &
      \sfS\A
      }}\quad
  \text{and}\quad
  \vcenter{\xymatrix{
      0\ar@{^(->}[r]\ar@{^(->}[d] &
      \sfS(\sfS\A)\ar@{=}[r]\ar@{^(->}[d] &
      \sfS(\sfS\A)\ar@{^(->}[d] \\
      \sfS(\sfS\A)\ar@{=}[d]\ar@{^(->}[r] &
      \sfC_{\rho}\ar@{->>}[r]\ar@{->>}[d] &
      \sfC(\sfS\A)\ar@{->>}[d] \\
      \sfS(\sfS\A)\ar@{^(->}[r] &
      \sfC(\sfS\A)\ar@{->>}[r]^{\rho} &
      \sfS\A
      }}$$
  such that everything commutes. 
  So by the above, we have $\partial_0=-\beta_\A$. 
  
 The proof when $\functorF$ is contravariant is dual. 
\end{proof}

\begin{proposition}\label{prop-K-theoryofmappingconesequence}
  Let $\functorF$ be an additive functor from $\mathgrc{S}$ 
  to the category $\mathgrc{Ab}$, which has the
  properties (H), (S), and (HX). 
  Let there be given an extension 
  $$e:\extwithmaps{\A_0}{\A_1}{\A_2}{\iota}{\pi}.$$ 
Then there exist isomorphisms of cyclic six term exact sequences as follows: 
  $$\def\objectstyle{\scriptstyle}
  \def\labelstyle{\scriptstyle}
  \xymatrix{
    \ar[r]^{-\functorF(\sfS\pi)}&
    \functorF(\sfS\A_2)\ar@{=}[d]\ar[r]^-{\partial_1^e}&
    \functorF(\A_0)\ar[d]^{\functorF(f_e)}_\cong\ar[r]^-{\functorF(\iota)}&
    \functorF(\A_1)\ar@{=}[d]\ar[r]^-{-\functorF(\pi)}&
    \functorF(\A_2)\ar[d]^{\beta_{\A_2}}_\cong\ar[r]^-{\partial_0^e}&
    \functorF(\sfS\A_0)\ar[d]^{\functorF(\sfS f_e)}_\cong\ar[r]^-{\functorF(\sfS\iota)}&
    \functorF(\sfS\A_1)\ar@{=}[d]\ar[r]^-{-\functorF(\sfS\pi)}&\\
    \ar[r]^{\partial_1^{\mc(e)}}
    &\functorF(\sfS\A_2)\ar[r]^{\functorF(\iota_\mc)}
    &\functorF(\sfC_\pi)\ar[r]^{\functorF(\pi_\mc)}
    &\functorF(\A_1)\ar[r]^{\partial_0^{\mc(e)}}
    &\functorF(\sfS\sfS\A_2)\ar[r]^{\functorF(\sfS\iota_\mc)}
    &\functorF(\sfS\sfC_\pi)\ar[r]^{\functorF(\sfS\pi_\mc)}
    &\functorF(\sfS\A_1)\ar[r]^-{\partial_1^{\mc(e)}}&
    }$$
  in the covariant case, and 
  $$\def\objectstyle{\scriptstyle}
  \def\labelstyle{\scriptstyle}
  \xymatrix{
    \ar[r]^-{-\functorF(\pi)}&
    \functorF(\A_1)\ar[r]^-{\functorF(\iota)}&
    \functorF(\A_0)\ar[r]^-{\partial^0_e}&
    \functorF(\sfS\A_2)\ar[r]^-{-\functorF(\sfS\pi)}&
    \functorF(\sfS\A_1)\ar[r]^-{\functorF(\sfS\iota)}&
    \functorF(\sfS\A_0)\ar[r]^-{\partial^1_e}&
    \functorF(\A_2)\ar[r]^-{-\functorF(\pi)}&
    \\
    \ar[r]^{\partial^1_{\mc(e)}}
    &\functorF(\A_1)\ar@{=}[u]\ar[r]^{\functorF(\pi_\mc)}
    &\functorF(\sfC_\pi)\ar[u]_{\functorF(f_e)}^\cong\ar[r]^{\functorF(\iota_\mc)}
    &\functorF(\sfS\A_2)\ar@{=}[u]\ar[r]^{\partial^0_{\mc(e)}}
    &\functorF(\sfS\A_1)\ar@{=}[u]\ar[r]^{\functorF(\sfS\pi_\mc)}
    &\functorF(\sfS\sfC_\pi)\ar[u]_{\functorF(\sfS f_e)}^\cong\ar[r]^{\functorF(\sfS\iota_\mc)}
    &\functorF(\sfS\sfS\A_2)\ar[u]_{\beta_{\A_2}^{-1}}^\cong\ar[r]^-{\partial^1_{\mc(e)}}&
    }$$
  in the contravariant case. 
\end{proposition}
\begin{proof}
  Assume that $\functorF$ is covariant. 
  The diagram 
  $$\xymatrix{\A_1\ar[r]^\pi\ar[d]^\pi&\A_2\ar@{=}[d]\\
    \A_2\ar@{=}[r]&\A_2}$$
  induces the morphism of extensions
  $$\xymatrix{0\ar[r]&\sfS\A_2\ar@{=}[d]\ar[r]&\sfC_\pi\ar[d]^\omega\ar[r]&\A_1\ar[d]^\pi\ar[r]&0\\
    0\ar[r]&\sfS\A_2\ar[r]&\sfC_{\ident_{\A_2}}\ar[r]&\A_2\ar[r]&0.}$$
  Note that $\sfC_{\ident_{\A_2}}$ is canonically isomorphic to $\sfC\A_2$.
  According to Lemma~\ref{lem-mappingconesequence-new}, 
  this induces a morphism between cyclic six term exact sequences:
  $$\def\objectstyle{\scriptstyle}
    \def\labelstyle{\scriptstyle}
    \xymatrix{
      \ar[r]^-{\partial_1^{e_\mc}}
      &\functorF(\sfS\A_2)\ar@{=}[d]\ar[r]^{\functorF(\iota_\mc)}
      &\functorF(\sfC_\pi)\ar[d]\ar[r]^{\functorF(\pi_\mc)}
      &\functorF(\A_1)\ar[d]^{\functorF(\pi)}\ar[r]^{\partial_0^{e_\mc}}
      &\functorF(\sfS\sfS\A_2)\ar@{=}[d]\ar[r]^{\functorF(\sfS\iota_\mc)}
      &\functorF(\sfS\sfC_\pi)\ar[d]\ar[r]^{\functorF(\sfS\iota_\mc)}
      &\functorF(\sfS\A_1)\ar[d]^{\functorF(\sfS\pi)}\ar[r]^-{\partial_1^{e_\mc}}&\\
      \ar[r]^-\cong_-{-\ident}
      &\functorF(\sfS\A_2)\ar[r]
      &0\ar[r]
      &\functorF(\A_2)\ar[r]^\cong_{-\beta_{\A_2}}
      &\functorF(\sfS\sfS\A_2)\ar[r]
      &0\ar[r]
      &\functorF(\sfS\A_2)\ar[r]^-\cong_-{-\ident}&}$$
  This takes care of the commutativity of two of the six squares. 

  Commutativity of 
  $$\xymatrix{
    \functorF(\A_0)\ar[d]^{\functorF(f_e)}_\cong\ar[r]^{\functorF(\iota)}
    &\functorF(\A_1)\ar@{=}[d]&
    \functorF(\sfS\A_0)\ar[d]^{\functorF(\sfS f_e)}_\cong\ar[r]^{\functorF(\sfS(\iota))}
    &\functorF(\sfS\A_1)\ar@{=}[d]\\
    \functorF(\sfC_\pi)\ar[r]^{\functorF(\pi_\mc)}
    &\functorF(\A_1)
    &\functorF(\sfS\sfC_\pi)\ar[r]^{\functorF(\sfS\pi_\mc)}
    &\functorF(\sfS\A_1)
    }$$
  follows directly from the $3\times 3$-diagram above.  Now we only need to check commutativity of 
  $$\xymatrix{
    \functorF(\sfS\A_2)\ar@{=}[d]\ar[r]^{\partial_1^e}&
    \functorF(\A_0)\ar[d]^{\functorF(f_e)}_\cong&
    \functorF(\A_2)\ar[d]^{\beta_{\A_2}}_\cong\ar[r]^{\partial_0^e}&
    \functorF(\sfS\A_0)\ar[d]^{\functorF(\sfS f_e)}_\cong\\
    \functorF(\sfS\A_2)\ar[r]^{\functorF(\iota_\mc)}&
    \functorF(\sfC_\pi)&
    \functorF(\sfS\sfS\A_2)\ar[r]^{\functorF(\sfS\iota_\mc)}&
    \functorF(\sfS\sfC_\pi)
    }$$
  Since $\sfC_\pi$ is the pullback, we get a canonical map
  $\sfC\A_1\rightarrow \sfC_\pi$ and commuting diagrams
  $$\def\objectstyle{\scriptstyle}
  \def\labelstyle{\scriptstyle}
  \xymatrix{
    \sfS\A_0\ar@{^(->}[r]\ar@{^(->}[d]&
    \sfC\A_0\ar@{->>}[r]\ar@{^(->}[d]&
    \A_0\ar@{^(->}[d]
    &
    \sfS\A_0\ar@{=}[r]\ar@{^(->}[d]&
    \sfS\A_0\ar@{->>}[r]\ar@{^(->}[d]&
    0\ar@{^(->}[d]
    &
    \sfS\A_0\ar@{=}[r]\ar@{^(->}[d]&
    \sfS\A_0\ar@{->>}[r]\ar@{^(->}[d]&
    0\ar@{^(->}[d] 
    \\
    \sfS\A_1\ar@{^(->}[r]\ar@{->>}[d]&
    \sfC\A_1\ar@{->>}[r]\ar@{->>}[d]&
    \A_1\ar@{->>}[d]
    &
    \sfC\A_0\ar@{^(->}[r]\ar@{->>}[d]&
    \sfC\A_1\ar@{->>}[r]\ar@{->>}[d]&
    \sfC\A_2\ar@{=}[d]
    &
    \sfS\A_1\ar@{^(->}[r]\ar@{->>}[d]&
    \sfC\A_1\ar@{->>}[r]\ar@{->>}[d]&
    \A_1\ar@{=}[d]
    \\
    \sfS\A_2\ar@{^(->}[r]&
    \sfC\A_2\ar@{->>}[r]&
    \A_2
    &
    \A_0\ar@{^(->}[r]^{f_e}&
    \sfC_\pi\ar@{->>}[r]&
    \sfC\A_2
    &
    \sfS\A_2\ar@{^(->}[r]&
    \sfC_\pi\ar@{->>}[r]&
    \A_1
    }$$
  with exact rows and columns.
  Using Lemma~\ref{lem-anticommute-kk} and
  Lemma~\ref{lem-mappingconesequence-new}, these diagrams give rise to 
  the following commutative diagrams 
  $$\def\objectstyle{\scriptstyle}
  \def\labelstyle{\scriptstyle}
  \xymatrix@C=20pt{
    \functorF(\A_2)\ar[r]^{-\beta_{\A_2}}\ar[d]_{\partial_0^e}&
    \functorF(\sfS\sfS\A_2)\ar[d]^{\partial_1^{\sfS e}}&
    \functorF(\A_0)\ar[r]^{\functorF(f_e)}\ar[d]_{-\beta_{\A_0}}&
    \functorF(\sfC_\pi)\ar[d]^{\partial_0^{e'}}&
    \functorF(\sfS\A_2)\ar[r]^{\functorF(\iota_\mc)}\ar[d]_{\partial_0^{\sfS e}}&
    \functorF(\sfC_\pi)\ar[d]^{\partial_0^{e'}}&\\
    \functorF(\sfS\A_0)\ar@{=}[r]&
    \functorF(\sfS\A_0)&
    \functorF(\sfS\sfS\A_0)\ar@{=}[r]&
    \functorF(\sfS\sfS\A_0)&
    \functorF(\sfS\sfS\A_0)\ar@{=}[r]&
    \functorF(\sfS\sfS\A_0)\\
    \functorF(\sfS\A_2)\ar@{=}[r]\ar[d]_{\partial_1^e}&
    \functorF(\sfS\A_2)\ar[d]^{\partial_0^{\sfS e}}&
    \functorF(\sfS\A_0)\ar[r]^{-\functorF(\sfS f_e)}\ar@{=}[d]&
    \functorF(\sfS\sfC_\pi)\ar[d]^{\partial_1^{e'}}&
    \functorF(\sfS\sfS\A_2)\ar[r]^{\functorF(\sfS\iota_\mc)}\ar[d]_{\partial_1^{\sfS e}}&
    \functorF(\sfS\sfC_\pi)\ar[d]^{\partial_1^{e'}}\\
    \functorF(\A_0)\ar[r]_{-\beta_{\A_0}}&
    \functorF(\sfS\sfS\A_0)&
    \functorF(\sfS\A_0)\ar@{=}[r]&
    \functorF(\sfS\A_0)&
    \functorF(\sfS\A_0)\ar@{=}[r]&
    \functorF(\sfS\A_0)
    }$$
  where $e'$ denotes the extension 
  $\extwithoutmaps{\sfS\A_0}{\sfC\A_1}{\sfC_\pi}$. 
  Consequently, 
  \begin{align*}
    \functorF(\iota_\mc)
    &=(\partial_0^{e'})^{-1}\circ\partial_0^{\sfS e}
    =-\functorF(f_e)\circ(\beta_{\A_0}^{-1})\circ\partial_0^{\sfS e}\\
    &=\functorF(f_e)\circ\beta_{\A_0}^{-1}\circ\beta_{\A_0}\circ\partial_1^{e}
    =\functorF(f_e)\circ\partial_1^{e},\\
    \functorF(\sfS\iota_\mc)
    &=(\partial_1^{e'})^{-1}\circ\partial_1^{\sfS e}
    =-\functorF(\sfS f_e)\circ\partial_1^{\sfS e}\\
    &=\functorF(\sfS f_e)\circ\partial_0^{e}\circ\beta_{\A_2}^{-1}.
  \end{align*}

  The proof in the contravariant case is dual. 
\end{proof}

\begin{corollary}
  Let $\functorF$ be an additive functor from $\mathgrc{S}$ 
  to the category $\mathgrc{Ab}$, which has the
  properties (H), (S), and (HX). 
  Let there be given a \starhomo 
  $$\phi:\A\rightarrow\B$$
  from a \ca \A to a \ca \B,  
  and let 
  $$e\colon\extwithmaps{\sfS\B}{\sfC_\phi}{\A}{\iota_\mc}{\pi_\mc}$$
  denote the corresponding mapping cone sequence.  Then there exist isomorphisms of cyclic six term exact sequences as follows: 
  $$\def\objectstyle{\scriptstyle}
  \def\labelstyle{\scriptstyle}
  \xymatrix{
    \ar[r]^-{-\functorF(\sfS\phi)}&
    \functorF(\sfS\B)\ar@{=}[d]\ar[r]^{\functorF(\iota_\mc)}&
    \functorF(\sfC_\phi)\ar@{=}[d]\ar[r]^{\functorF(\pi_\mc)}&
    \functorF(\A)\ar@{=}[d]\ar[r]^-{-\functorF(\phi)}&
    \functorF(\B)\ar[d]^{\beta_{\B}}_\cong\ar[r]^{\functorF(\sfS\iota_\mc)\circ\beta_\B}&
    \functorF(\sfS\sfC_\phi)\ar@{=}[d]\ar[r]^{\functorF(\sfS\pi_\mc)}&
    \functorF(\sfS\A)\ar@{=}[d]\ar[r]^-{-\functorF(\sfS\phi)}&\\
    \ar[r]^-{\partial_1^{e}}&
    \functorF(\sfS\B)\ar[r]^{\functorF(\iota_\mc)}&
    \functorF(\sfC_\phi)\ar[r]^{\functorF(\pi_\mc)}&
    \functorF(\A)\ar[r]^{\partial_0^{e}}&
    \functorF(\sfS\sfS\B)\ar[r]^{\functorF(\sfS\iota_\mc)}&
    \functorF(\sfS\sfC_\phi)\ar[r]^{\functorF(\sfS\pi_\mc)}&
    \functorF(\sfS\A)\ar[r]^-{\partial_1^{e}}&
    }$$
  in the covariant case, and 
  $$\def\objectstyle{\scriptstyle}
  \def\labelstyle{\scriptstyle}
  \xymatrix{
    \ar[r]^-{-\functorF(\phi)}&
    \functorF(\A)\ar@{=}[d]\ar[r]^{\functorF(\pi_\mc)}&
    \functorF(\sfC_\phi)\ar@{=}[d]\ar[r]^{\functorF(\iota_\mc)}&
    \functorF(\sfS\B)\ar@{=}[d]\ar[r]^-{-\functorF(\sfS\phi)}&
    \functorF(\sfS\A)\ar@{=}[d]\ar[r]^{\functorF(\sfS\pi_\mc)}&
    \functorF(\sfS\sfC_\phi)\ar@{=}[d]\ar[r]^{\functorF(\sfS\iota_\mc)}&
    \functorF(\B)\ar[d]^{\beta_{\B}}_\cong\ar[r]^-{-\functorF(\phi)}&\\
    \ar[r]^-{\partial^1_{e}}&
    \functorF(\A)\ar[r]^{\functorF(\pi_\mc)}&
    \functorF(\sfC_\phi)\ar[r]^{\functorF(\iota_\mc)}&
    \functorF(\sfS\B)\ar[r]^{\partial^0_{e}}&
    \functorF(\sfS\A)\ar[r]^{\functorF(\sfS\pi_\mc)}&
    \functorF(\sfS\sfC_\phi)\ar[r]^{\beta_{\B}\circ\functorF(\sfS\iota_\mc)}&
    \functorF(\sfS\sfS\B)\ar[r]^-{\partial^1_{e}}&
    }$$
  in the contravariant case.
\end{corollary}
\begin{proof}
  This follows from the first part of the proof of the previous proposition. 
\end{proof}

\section{Examples of concrete homology and cohomology theories}
\label{sec-examplesofconcretehomologytheories}

\begin{example}
  Let $\mathgrc{S}$ be the full subcategory of the category of all
  \cas, consisting of separable, nuclear algebras. 
  For each separable \ca \A, both $\kk(-,\A)$ and $\kk(\A,-)$ 
  are additive functors from $\mathgrc{S}$ to $\mathgrc{Ab}$, 
  which have the properties (H), (S), and (HX). 
  The first one is contravariant while the second is covariant. 
  So the above theory applies to these, and identifies the cyclic six
  term exact sequences 
  associated with extensions in these two cases (as defined in
  \cite{MR99g:46104}).  
\end{example}

\begin{example}\label{ex-K-theory}
  The functors $K_0$ and $K_1$ 
  are additive, covariant functors from the category of all
  separable \cas to the category $\mathgrc{Ab}$, which have the
  properties (H), (S), and (HX). 
  So the above theory applies to these two functors. 

  We have also a standard cyclic six term exact sequence in
  $K$-theory (as defined in \cite{MR2001g:46001}). 
  To avoid confusions, we write $\delta_0$ and $\delta_1$ for the
  exponential map and the index maps, respectively. 
  We will recall the definition here. 
  We have an isomorphism $\theta_{-}$ of functors from $K_1(-)$ to
  $K_0(\sfS(-))$, \ie, for each \ca \A we have an isomorphism
  $\theta_\A\colon K_1(\A)\rightarrow K_0(\sfS\A)$ and, moreover, for
  all \cas $\A$ and \B and all \starhomos
  $\varphi\colon\A\rightarrow\B$, the diagram
  $$\xymatrix{K_1(\A)\ar[r]^{K_1(\varphi)}\ar[d]_{\theta_\A}&
    K_1(\B)\ar[d]^{\theta_\B}\\
    K_0(\sfS\A)\ar[r]_{K_0(\sfS\varphi)}&K_0(\sfS\B)}$$
  commutes (\cf\ \cite[Theorem~10.1.3]{MR2001g:46001}). 

  The exponential map $\delta_0\colon K_0(\A_2)\rightarrow K_1(\A_0)$
  associated with a short exact sequence
  $\extwithoutmaps{\A_0}{\A_1}{\A_2}$ 
  is defined as the composition of the maps 
  $$\xymatrix{K_0(\A_2)\ar[r]^{\beta_{\A_2}}&
    K_1(\sfS\A_2)\ar[r]^{\overline{\delta_1}}&
    K_0(\sfS\A_0)\ar[r]^{\theta_{\A_0}^{-1}}&K_1(\A_0),}$$
  where $\overline{\delta_1}$ is the index map associated with the
  short exact sequence
  $$\extwithoutmaps{\sfS\A_0}{\sfS\A_1}{\sfS\A_2}.$$

  \begin{lemma}\label{lem-mappingconesequence-extranew}
    Let \A be a \ca. 
    The \scstes in $K$-theory associated with 
    $\extwithoutmaps{\sfS\A}{\sfC\A}{\A}$ 
    (as in \cite{MR2001g:46001}) is the
    sequence 
    $$\xymatrix{K_0(\sfS\A)\ar[r]&0\ar[r]&K_0(\A)\ar[d]_\cong^{-\beta_\A}\\
      K_1(\A)\ar[u]_\cong^{\theta_\A}&0\ar[l]&K_1(\sfS\A).\ar[l]}$$
  \end{lemma}
  \begin{proof}
    Since the cone, $\sfC\A$, of $\A$ is homotopy equivalent to the zero
    \ca, $K_0(\sfC\A)\cong K_1(\sfC\A)\cong 0$ 
    (\cf\ \cite[Example~4.1.5]{MR2001g:46001}).  That the index map is $\theta_\A$ follows directly from the
    definition of $\theta_\A$ 
    (\cf\ \cite[Proof of Theorem~10.1.3]{MR2001g:46001}).  The exponential map $\delta_0\colon K_0(\A)\rightarrow K_1(\sfS\A)$ is
    defined as the composition of the maps 
    $$\xymatrix{K_0(\A)\ar[r]^{\beta_\A}&
      K_1(\sfS\A)\ar[r]^{\overline{\delta_1}}&
      K_0(\sfS(\sfS\A))\ar[r]^{\theta_{\sfS\A}^{-1}}&K_1(\sfS\A),}$$
    where $\overline{\delta_1}$ is the index map associated with the
    short exact sequence
    $$\extwithoutmaps{\sfS(\sfS\A)}{\sfS(\sfC\A)}{\sfS\A}.$$
    
    We have a commuting diagram 
    $$\xymatrix{
      \sfS(\sfS\A)\ar@{^(->}[r]\ar@{^(->}[d]&\sfS(\sfC\A)\ar@{->>}[r]\ar@{^(->}[d]&\sfS(\A)\ar@{^(->}[d]\\
      \sfC(\sfS\A)\ar@{^(->}[r]\ar@{->>}[d]&\sfC(\sfC\A)\ar@{->>}[r]\ar@{->>}[d]&\sfC(\A)\ar@{->>}[d]\\
      \sfS\A\ar@{^(->}[r]&\sfC\A\ar@{->>}[r]&\A
      }$$
    with exact rows and columns. 
    This gives --- by Lemma~\ref{lem-anticommute-kk} and the above
    (applied to $\sfS\A$ instead of $\A$) --- rise to an anticommuting square
    $$\xymatrix{
      K_0(\A)\ar[r]^{\delta_0}_\cong\ar[d]^{\delta_0}_\cong&
      K_1(\sfS\A)\ar[d]^{\theta_{\sfS\A}}\\
      K_1(\sfS(\A))\ar[r]^{\overline{\delta_1}}&
      K_0(\sfS(\sfS\A))
      }$$
    Consequently, $\overline{\delta_1}=-\theta_{\sfS\A}$. 
    Now it follows that $\delta_0=-\beta_\A$. 
  \end{proof}
  
  Since the index and exponential maps are unique up to signs 
  (\cf\ \cite[Exercise~9.F]{MR95c:46116}), we have that the standard cyclic six term
  exact sequence in $K$-theory as defined here 
  differs from the cyclic six term exact
  sequence defined as above by change of sign of the index map 
  (under the identification $\theta_{-}$ of $K_1$ with $K_0\circ\sfS$). 
\end{example}

Thus we get the following corollaries: 

\begin{corollary}\label{cor:ksixmappingcone}
  Let there be given an extension 
  $$e:\extwithmaps{\A_0}{\A_1}{\A_2}{\iota}{\pi}.$$ 
  Then there exists an isomorphism of cyclic six term exact sequences as follows: 
  $$\def\objectstyle{\scriptstyle}
    \def\labelstyle{\scriptstyle}
    \xymatrix{
    \ar[r]&
    K_1(\A_2)\ar[d]^{\theta_{\A_2}}_\cong\ar[r]^-{-\delta_1^e}&
    K_0(\A_0)\ar[d]^{K_0(f_e)}_\cong\ar[r]^-{K_0(\iota)}&
    K_0(\A_1)\ar@{=}[d]\ar[r]^-{-K_0(\pi)}&
    K_0(\A_2)\ar[d]^{\beta_{\A_2}}_\cong\ar[r]^-{\delta_0^e}&
    K_1(\A_0)\ar[d]^{K_1(f_e)}_\cong\ar[r]^-{K_1(\iota)}&
    K_1(\A_1)\ar@{=}[d]\ar[r]^-{K_1(\pi)}&\\
    \ar[r]&K_0(\sfS\A_2)\ar[r]&K_0(\sfC_\pi)\ar[r]&K_0(\A_1)\ar[r]&    
    K_1(\sfS\A_2)\ar[r]&K_1(\sfC_\pi)\ar[r]&K_1(\A_1)\ar[r]&
    }$$
  {where} the second sequence
  is the \scstes 
  in  $K$-theory associated with $\mc(e)$. 
\end{corollary}

\begin{corollary}
  Let there be given a \starhomo 
  $$\phi:\A\rightarrow\B$$
  from a \ca \A to a \ca \B,  
  and let 
  $$e\colon\extwithmaps{\sfS\B}{\sfC_\phi}{\A}{\iota_{\mc}}{\pi_{\mc}}$$  
  denote the mapping cone sequence.  Then there exists an isomorphism of exact sequences as follows: 
  $$\def\objectstyle{\scriptstyle}
  \def\labelstyle{\scriptstyle}
    \xymatrix{
    \ar[r]&
    K_1(\B)\ar[d]^{\theta_{\B}}_\cong\ar[r]&
    K_0(\sfC_\phi)\ar@{=}[d]\ar[r]&
    K_0(\A)\ar@{=}[d]\ar[r]^-{-K_0(\phi)}&
    K_0(\B)\ar[d]^{\beta_{\B}}_\cong\ar[r]&
    K_1(\sfC_\phi)\ar@{=}[d]\ar[r]&
    K_1(\A)\ar@{=}[d]\ar[r]^-{K_1(\phi)}&\\
    \ar[r]&K_0(\sfS\B)\ar[r]&K_0(\sfC_\phi)\ar[r]&K_0(\A)\ar[r]&    
    K_1(\sfS\B)\ar[r]&K_1(\sfC_\phi)\ar[r]&K_1(\A)\ar[r]&
    }$$
  {where} 
  the second sequence
  is the \scstes 
  in $K$-theory associated with $e$. 
\end{corollary}

\begin{remark}
  Note that the way Bonkat associates cyclic six term exact
  sequences in ideal-related \kk-theory with short exact sequences 
  with completely positive contractive coherent splittings is
  completely analogous to the definitions of
  Section~\ref{sec-homologyandcohomologytheories} 
  (\cf\ \bonkat[Section~3.4]). 
\end{remark}

\begin{example}
  An instructive example to get a better understanding of 
  Lemma~\ref{lem-anticommute-kk} is to look at 
  $$\xymatrix{    
    \sfS\otimes\sfS\ar@{^(->}[d]\ar@{^(->}[r]&\sfS\otimes\sfC\ar@{^(->}[d]\ar@{->>}[r]&\sfS\otimes\C\ar@{^(->}[d]\\
    \sfC\otimes\sfS\ar@{->>}[d]\ar@{^(->}[r]&\sfC\otimes\sfC\ar@{->>}[d]\ar@{->>}[r]&\sfC\otimes\C\ar@{->>}[d]\\
    \C\otimes\sfS\ar@{^(->}[r]&\C\otimes\sfC\ar@{->>}[r]&\C\otimes\C 
    }$$
  where $\sfS=\sfS\C$ and $\sfC=\sfC\C$. 
  It is tempting to guess that the maps 
  $$K_0(\C\otimes\C)\rightarrow 
  K_1(\sfS\otimes\C)\rightarrow K_0(\sfS\otimes \sfS)$$
  $$K_0(\C\otimes\C)\rightarrow 
  K_1(\C\otimes \sfS)\rightarrow K_0(\sfS\otimes \sfS)$$ 
  are equal (after all,
  $\sfS\otimes\C$ is canonically isomorphic to $\C\otimes \sfS$) --- but
  this is not the case. 
  One map gives the Bott map while the other gives the anti-Bott map. 
  After some thought this seems reasonable after all, since the map 
  $\sfS\otimes \sfS\ni x\otimes y\mapsto y\otimes x\in \sfS\otimes \sfS$
  corresponds to the flip along the diagonal in
  $C_0((0,1)\times(0,1))$, which induces the automorphism $-\ident$ on
  $K_0$. 
\end{example}

\section{Ideal-related $K$-theory with coefficients}

In this section, we introduce a new invariant, \emph{Ideal-related
$K$-theory with coefficients}. 
  
\begin{definition}\label{def-dimensiondropinterval}
  Let $n\in\N_{\geq 2}$.  Denote the 
  (non-unital) \defemph{dimension-drop interval} by $\Io$, \ie, $\Io$ is the
  mapping cone of the unital \starhomo from $\C$ to $\Matn$.
\end{definition}

\begin{definition}
  Let $n\in\N_{\geq 2}$ and let \eo denote the mapping cone sequence 
  $$\eo\colon\extwithoutmaps{\sfS\Matn}{\Io}{\C}$$
  corresponding to the unital \starhomo from $\C$ to $\Matn$. 
  Moreover, let $\e{n}{i}=\mc^i(\eo)$, for all $i\in\N$ and write
  $$\ei\colon\extwithoutmaps{\sfS\C}{\Ii}{\Io},$$
  $$\e{n}{i}\colon\extwithoutmaps{\sfS\mathbb{I}_{n,i-2}}{\mathbb{I}_{n,i}}{\mathbb{I}_{n,i-1}}, 
  \text{ for }i\geq 2.$$
Similarly, 
  set $\Fo[1]\colon\extwithmaps{\C}{\C}{0}{\ident}{}$
  and 
  $\Fo\colon\extwithmaps{\Io}{\Io}{0}{\ident}{}$, 
  for all $n\in\N_{\geq 2}$. 
  Moreover, 
  set $\FFFF{n}{i}=\mc^i(\Fo)$ for all $n\in\N$
  and all $i\in\N$. 
\end{definition}

\begin{definition}
  Let $\ksix$ denote the functor, which to each extension of \cas
  associates the corresponding standard cyclic six term exact
  sequence (as defined in \cite{MR2001g:46001} 
  --- \cf\ Example~\ref{ex-K-theory}).  Let $\Homsix(\ksix(e_1),\ksix(e_2))$ denote the group of cyclic six term chain homomorphisms.

  As in \cite{MR1404333}, set 
  $K_i(-;\Z_n)=\kk^i(\Io,-)$. 
  Moreover, let $\underline{K}$ denote total $K$-theory as defined
  in \cite{MR1404333}. 
\end{definition}

\begin{remark}
  As is easily seen, the above cyclic six term exact sequence in
  $K$-theory differs from that defined by Bonkat in
  \cite[\S{}7.3]{bonkat} by the index and exponential maps having the
  opposite signs. 
  This makes no difference for the arguments and results in
  \cite{bonkat} (the important thing here is that we change the sign of
  either the index map or the exponential map compared with the
  definition of the connecting homomorphisms in \kk-theory). 

  By applying Lemma~\ref{lem-anticommute-kk} and 
  Lemma~\ref{lem-mappingconesequence-extranew} to 
  $$\xymatrix{
    \sfS\A_0\ar@{^(->}[r]^{\sfS\iota}\ar@{^(->}[d] &
    \sfS\A_1\ar@{->>}[r]^{\sfS\pi}\ar@{^(->}[d] &
    \sfS\A_2\ar@{^(->}[d] \\
    \sfC\A_0\ar@{^(->}[r]^{\sfC\iota}\ar@{->>}[d] &
    \sfC\A_1\ar@{->>}[r]^{\sfC\pi}\ar@{->>}[d] &
    \sfC\A_2\ar@{->>}[d] \\
    \A_0\ar@{^(->}[r]^{\iota} &
    \A_1\ar@{->>}[r]^\pi &
    \A_2 
    }$$
 there exists a commuting diagram 
  $$\xymatrix{
    \ar[r]^-{\delta_1^e} &
    K_0(\A_0)\ar[d]_\cong^{\beta_{\A_0}}\ar[r]^{K_0(\iota)} &
    K_0(\A_1)\ar[d]_\cong^{\beta_{\A_1}}\ar[r]^{K_0(\pi)} &
    K_0(\A_2)\ar[d]_\cong^{\beta_{\A_2}}\ar[r]^{\delta_0^e} &
    K_1(\A_0)\ar[d]_\cong^{\beta_{\A_0}}\ar[r]^{K_1(\iota)} &
    K_1(\A_1)\ar[d]_\cong^{\beta_{\A_1}}\ar[r]^{K_1(\pi)} &
    K_1(\A_2)\ar[d]_\cong^{\beta_{\A_2}}\ar[r]^-{\delta_1^e} & \\
    \ar[r]^-{\delta_0^{\sfS e}} &
    K_1(\sfS\A_0)\ar[r]^{K_1(\sfS\iota)} &
    K_1(\sfS\A_1)\ar[r]^{K_1(\sfS\pi)} &
    K_1(\sfS\A_2)\ar[r]^{\delta_1^{\sfS e}} &
    K_0(\sfS\A_0)\ar[r]^{K_0(\sfS\iota)} &
    K_0(\sfS\A_1)\ar[r]^{K_0(\sfS\pi)} &
    K_0(\sfS\A_2)\ar[r]^-{\delta_0^{\sfS e}} & 
     }$$
   Consequently, the definition of ``$(K_{*+1}A_i)$'' in \cite{bonkat}
   is just $\ksix(\sfS e)$ (up to canonical identification with our terminology). 
   The same argument works if we choose to work with the slightly
   different cyclic six term exact sequence defined in \cite{bonkat}.
   Note also that this is not true if we define the cyclic six term
   sequence using the abstract machinery of 
   Section~\ref{sec-homologyandcohomologytheories}. 
\end{remark}

\begin{definition}
{%
  For each extension $e$ of separable \cas, define the \defemph{ideal-related $K$-theory with coefficients},
  $\kE(e)$, of $e$ to be the (graded) group 
  $$\kE(e)=\bigoplus_{i=0}^5
  \left(\kkE(\FFFF{1}{i},e)\oplus\bigoplus_{n=2}^\infty\kkE(\e{n}{i},e)\oplus\kkE(\FFFF{n}{i},e)\right).$$
 A homomorphism $\alpha$ from $\kE(e_1)$ to $\kE(e_2)$ is a group homomorphism
  respecting the direct sum
  decomposition and the natural homomorphisms induced by the elements
  of $\kkE^j(e,e')$, for $j=0,1$, where $e$ and $e'$ are in
  $\sets{\e{n}{i},\FFFF{n}{i},\FFFF{1}{i}}{n\in\N_{\geq 2},i=0,1,2}$. 
  The set of homomorphisms from $\kE(e_1)$ to $\kE(e_2)$ will be denoted by $\Hom_\Lambda(\kE(e_1),\kE(e_2))$.  
  
  Let $x\in\kkE(e_1,e_2)$. Then $x$ induces an element of $\Hom_\Lambda(\kE(e_1),\kE(e_2))$ by 
  \begin{align*}
  &y\in\kkE(\FFFF{n}{i},e_1)\mapsto y\times x\in\kkE(\FFFF{n}{i},e_2),\quad n\in\N, \\
  &y\in\kkE(\e{n}{i},e_1)\mapsto y\times x\in\kkE(\e{n}{i},e_2),\quad n\in\N_{\geq 2}.
  \end{align*}
  Hence, if $\phi\colon e_1\rightarrow e_2$ is a homomorphism, then $\phi$ induces an element $\kE(\phi)\in\Hom_\Lambda(\kE(e_1),\kE(e_2))$. In this way $\kE$ becomes a functor on the category of extensions.}
\end{definition}

\begin{remark}
  For extensions $e_1:\extwithoutmaps{\A_0}{\A_1}{\A_2}$ and
  $e_2:\extwithoutmaps{\B_0}{\B_1}{\B_2}$ of separable \cas, there are natural
  homomorphisms $G_i\colon\kkE(e_1,e_2)\longrightarrow\kk(\A_i,\B_i)$,
  for $i=0,1,2$.  As in the proof of \cite[Satz~7.5.6]{bonkat}, the obvious diagram
  $$\def\objectstyle{\scriptstyle} \def\labelstyle{\scriptstyle}
  \xymatrix{
    \operatorname{Ext}_{\mathrm{six}}(\ksix(e_1),\ksix(\sfS e_2))\ar@{^(->}[r]\ar[d]
    & \kkE(e_1,e_2)\ar@{->>}[r]\ar[d]^{G_i}
    & \Homsix(\ksix(e_1),\ksix(e_2))\ar[d]\\
    \operatorname{Ext}(K_0(\A_i),K_1(\B_i))
    \oplus\operatorname{Ext}(K_1(\A_i),K_0(\B_i)) \ar@{^(->}[r] &
    \kk(\A_i,\B_i)\ar@{->>}[r] & \Hom(K_0(\A_i),K_0(\B_i))
    \oplus\Hom(K_1(\A_i),K_1(\B_i)) }$$
  commutes and is natural in
  $e_2$, for $i=0,1,2$ --- provided that $e_1$ belongs to the UCT
  class considered by Bonkat.
\end{remark}


Let  
$e\colon\extwithmaps{\A_0}{\A_1}{\A_2}{\iota}{\pi}$ 
and $e'\colon\extwithmaps{\B_0}{\B_1}{\B_2}{\iota'}{\pi'}$ 
be two given extensions. 
Then define 
\begin{equation*}
  \Lambda_{e,e'}\colon
  \Homsix(\ksix(e),\ksix(e'))  
  \longrightarrow
  \Homsix(\ksix(\mc(e)),\ksix(\mc(e'))) 
\end{equation*}
as follows:
Let $(\alpha_i)_{i=0}^5$ in 
$\Homsix(\ksix(e),\ksix(e'))$ be given. 
Then by Corollary~\ref{cor:ksixmappingcone} the diagram
$$\xymatrix{
  \ar[r]^{\delta_1^{\mc(e)}} 
  & K_0(\sfS\A_2) \ar[r]^{K_0(\iota_\mc)}\ar[d]^{\theta_{\A_2}^{-1}} 
  & K_0(\sfC_\pi) \ar[r]^{K_0(\pi_\mc)} \ar[d]^{K_0(f_e)^{-1}}  
  & K_0(\A_1) \ar[r]^{\delta_0^{\mc(e)}}\ar@{=}[d] 
  & K_1(\sfS\A_2) \ar[r]^{K_1(\iota_\mc)}\ar[d]^{\beta_{\A_2}^{-1}} 
  & K_1(\sfC_\pi) \ar[r]^{K_1(\pi_\mc)} \ar[d]^{K_1(f_e)^{-1}}   
  & K_1(\A_1) \ar[r]^-{\delta_1^{\mc(e)}}\ar@{=}[d] & \\
  \ar[r]^{K_1(\pi)} 
  & K_1(\A_2) \ar[r]^{-\delta_1^{e}}\ar[d]^{\alpha_5} 
  & K_0(\A_0) \ar[r]^{K_0(\iota)}\ar[d]^{\alpha_0} 
  & K_0(\A_1) \ar[r]^{-K_0(\pi)}\ar[d]^{\alpha_1} 
  & K_0(\A_2) \ar[r]^{\delta_0^{e}}\ar[d]^{\alpha_2} 
  & K_1(\A_0) \ar[r]^{K_1(\iota)}\ar[d]^{\alpha_3} 
  & K_1(\A_1) \ar[r]^-{K_1(\pi)}\ar[d]^{\alpha_4} &\\
  \ar[r]^{K_1(\pi')} 
  & K_1(\B_2) \ar[r]^{-\delta_1^{e'}}\ar[d]^{\theta_{\B_2}} 
  & K_0(\B_0) \ar[r]^{K_0(\iota')}\ar[d]^{K_0(f_{e'})} 
  & K_0(\B_1) \ar[r]^{-K_0(\pi')}\ar@{=}[d] 
  & K_0(\B_2) \ar[r]^{\delta_0^{e'}}\ar[d]^{\beta_{\B_2}} 
  & K_1(\B_0) \ar[r]^{K_1(\iota')}\ar[d]^{K_1(f_{e'})}  
  & K_1(\B_1) \ar[r]^-{K_1(\pi')}\ar@{=}[d] & \\
  \ar[r]^{\delta_1^{\mc(e')}} 
  & K_0(\sfS\B_2) \ar[r]^{K_0(\iota_\mc')} 
  & K_0(\sfC_{\pi'}) \ar[r]^{K_0(\pi_\mc')} 
  & K_0(\B_1) \ar[r]^{\delta_0^{\mc(e')}} 
  & K_1(\sfS\B_2) \ar[r]^{K_1(\iota_\mc')} 
  & K_1(\sfC_{\pi'}) \ar[r]^{K_1(\pi_\mc')}  
  & K_1(\B_1) \ar[r]^-{\delta_1^{\mc(e')}} & 
  }$$
commutes. 
Let $\Lambda_{e,e'}((\alpha_i)_{i=0}^5)$ denote the composition of
these maps. 
Clearly, 
$\Lambda_{e,e'}$ is an isomorphism.  
A computation shows that $\Lambda$ from $\Homsix(\ksix(e),\ksix(-))$ 
to $\Homsix(\ksix(\mc(e)),\ksix(\mc(-)))$ 
defined by $\Lambda(e')=\Lambda_{e,e'}$ 
is a natural transformation such that 
$\Lambda_{e,e}(\ksix(\ident_e))=\ksix(\ident_{\mc(e)})$.

Let $S\mathcal{E}$ be the subcategory of $\mathcal{E}$ consisting
only of extensions of separable \cas and morphism being triples of
\starhomos such that the obvious diagram commutes.  Consider the category $\catkkE$ whose objects are the objects of 
$S\mathcal{E}$ and the group of
morphisms is $\kkE(e_1,e_2)$.  
Consider the the composed functor $\kkE\circ\mc$ from $S\mathcal{E}$
to $\catkkE$, which sends an object $e$ of $S\mathcal{E}$ to
$\mc(e)$, and sends a morphism $(\phi_0,\phi_1,\phi_2)$ of
$S\mathcal{E}$ to $\kkE(\mc((\phi_0,\phi_1,\phi_2)))$. 
This is a stable, homotopy invariant, split exact functor, so by \bonkat[Satz~3.5.10 und~Satz~6.2.4], there exists a
unique functor $\widehat{\mc}$ from $\catkkE$ to $\catkkE$ such that the diagram 
$$\xymatrix{ 
  S\mathcal{E} \ar[r]^{\mc}\ar[d]_{\kkE} 
  & S\mathcal{E}\ar[d]^\kkE \\
  \catkkE\ar[r]_{\widehat{\mc}} & \catkkE  
  }$$
commutes.  By the universal property, the diagram 
\begin{equation*}
  \xymatrix{
    \kkE(e,e') \ar[r]\ar[d]_{\widehat{\mc}}\ar[d] 
    & \Homsix(\ksix(e),\ksix(e'))\ar[d]^{\Lambda(e')} \\
    \kkE(\mc(e),\mc(e')) \ar[r]   
    &  \Homsix(\ksix(\mc(e)),\ksix(\mc(e')))
    }
\end{equation*}
commutes, where the horizontal arrows are the natural maps in the UCT.

\begin{lemma}\label{lem:genBP1}
  Let $e$ and $e'$ be extensions of separable, nuclear \cas in \uct.  
  Then $\widehat{\mc}$ induces an isomorphism from $\kkE(e,e')$ to
  $\kkE(\mc(e),\mc(e'))$, 
  which is natural in both variables.
  \end{lemma}
\begin{proof}
  Let $\alpha_{e,e'}$ denote the map from $\kkE(e,e')$ to
  $\kkE(\mc(e),\mc(e'))$ induced by the functor $\widehat{\mc}$. 
  Since $\widehat{\mc}$ is a functor, clearly the map is going to be 
  natural (in both variables).  From Proposition~3.5.6 in \cite{bonkat} 
  (\cf\ also \cite[Lemma~3.2]{MR869779}), 
  it follows that $\widehat{\mc}$ is a group homomorphism.  Since $\Lambda_{e,e'}$ is an isomorphism, from the above diagram and
  the UCT of Bonkat \cite{bonkat}, we have that $\alpha_{e,e'}$ is an
  isomorphism whenever $\ksix(e')$ is injective.  

  When $e'$ is an arbitrary extension, then by 
  \cite[Proposition~7.4.3]{bonkat},  
  there exist an injective geometric resolution 
  $\extwithoutmaps{e_1}{e_2}{\sfS e'}$
  of $e'$, \ie, there exists a short exact sequence
  $\extwithoutmaps{e_1}{e_2}{\sfS e'}$ of extensions from
  $S\mathcal{E}$, with a completely positive  contractive coherent
  splitting, such that the induced six term exact $\ksix$-sequence
  degenerates to a short exact sequence 
  $\extwithoutmaps{\ksix(\sfS\sfS e')}{\ksix(e_1)}{\ksix(e_2)}$, which
  is an
  injective resolution of $\ksix(\sfS\sfS e')$.

  The cyclic six term exact sequences in $\kkE$-theory give a commuting diagram
  $$\def\objectstyle{\scriptstyle} \def\labelstyle{\scriptstyle}
  \xymatrix{
    \kkE(e,\sfS e_1)\ar[r]^-{\alpha_{e,\sfS e_1}}\ar[d] & 
    \kkE(\mc(e),\mc(\sfS e_1))\ar[d]\ar[r] &
    \kkE(\mc(e),\sfS\mc( e_1))\ar[d] \\
    \kkE(e,\sfS e_2)\ar[r]^-{\alpha_{e,\sfS e_2}}\ar[d] & 
    \kkE(\mc(e),\mc(\sfS e_2))\ar[d]\ar[r] &
    \kkE(\mc(e),\sfS\mc(e_2))\ar[d] \\
    \kkE(e,\sfS\sfS e')\ar[r]^-{\alpha_{e,\sfS\sfS e'}}\ar[d] & 
    \kkE(\mc(e),\mc(\sfS\sfS e'))\ar[r] &
    \kkE(\mc(e),\sfS\mc(\sfS e'))\ar[d] \\
    \kkE(e,e_1)\ar[r]^-{\alpha_{e,e_1}}\ar[d] & 
    \kkE(\mc(e),\mc(e_1))\ar[d]\ar@{=}[r] &
    \kkE(\mc(e),\mc(e_1))\ar[d] \\
    \kkE(e,e_2)\ar[r]^-{\alpha_{e,e_2}} &
    \kkE(\mc(e),\mc(e_2))\ar@{=}[r] & 
    \kkE(\mc(e),\mc(e_2)) 
    }$$
  with exact columns. 
  Naturality of $\alpha_{e,-}$ gives us commutativity of the squares on
  the left hand side, while naturality of the isomorphism from the
  functor $\mc\circ\sfS$ to the functor $\sfS\circ\mc$ 
  gives us commutativity of the squares on the right hand side 
  (\cf\ Lemma~\ref{lem-commutatorofMCwithSandC}). 
  The remaining rectangle is seen to commute by using the definition
  of the connecting homomorphisms and 
  Lemma~\ref{lem-bigbigbig-commutingdiagram}.
  By the Five Lemma, we have that $\alpha_{e,\sfS\sfS e'}$ is an
  isomorphism. 
  Therefore also $\alpha_{e,e'}$. 
\end{proof}

\begin{remark}
  Similarly, there exists a unique functor $\widehat{\sfS}$ from
  $\catkkE$ to $\catkkE$ such that the diagram 
  $$\xymatrix{
    S\mathcal{E}\ar[d]_{\kkE}\ar[r]^{\sfS} &
    S\mathcal{E}\ar[d]^{\kkE} \\
    \catkkE\ar[r]_{\widehat{\sfS}} & \catkkE
    }$$
  commutes. 
\end{remark}

\section{Some diagrams}
\label{sec-somediagrams}

In this section we construct 19 diagrams involving the groups of the
new invariant. 
These diagrams can in many cases be used to determine the new groups
introduced in the invariant (see Section~\ref{sec-examples} for some examples). 
They {are also used in the paper \cite{automorphism-se-gr-er}}, where the three authors
prove a Universal Multi-Coefficient Theorem (UMCT) for a certain class of \cas with one specified ideal,
which includes all the Cuntz-Krieger algebras of type (II) with one
specified ideal. 
The long proof of these diagrams is outlined in the next section.

\begin{assumption}
  Throughout this section,
  $e\colon\extwithmaps{\A_0}{\A_1}{\A_2}{\iota}{\pi}$ is a {(fixed)}
  extension of separable \cas. 
\end{assumption}

\begin{definition}
  Set $F_{1,i}=\kkE(\FFFF{1}{i},e)$, $F_{n,i}=\kkE(\FFFF{n}{i},e)$, and $H_{n,i}=\kkE(\e{n}{i},e)$, for
  all $n\in\N_{\geq 2}$ and all $i=0,1,2,3,4,5$. 
  For convenience, 
  identify indices modulo $6$, \ie,
  we write $F_{n,6}=F_{n,0}, F_{n,7}=F_{n,1}$ etc.
\end{definition}

\begin{remark}
  Let $e\colon\extwithmaps{\A_0}{\A_1}{\A_2}{\iota}{\pi}$ be a given
  extension of \cas. 
  Then consider the extensions
  $$\mc^3(e)\colon
  \extwithmaps{\sfS\sfC_\pi}{\sfC_{(\pi_{\mc})_{\mc}}}{\sfC_{\pi_\mc}}%
  {((\iota_{\mc})_{\mc})_{\mc}}{((\pi_{\mc})_{\mc})_{\mc}}$$
  and 
  $$\sfS(e)\colon\extwithmaps{\sfS\A_0}{\sfS\A_1}{\sfS\A_2}{\sfS\iota}{\sfS\pi}.$$
  There are canonical \starhomos $\sfS\A_0\rightarrow\sfS\sfC_\pi$,
  $\sfS\A_1\rightarrow\sfC_{(\pi_{\mc})_{\mc}}$, and 
  $\sfS\A_2\rightarrow\sfC_{\pi_\mc}$, which all induce isomorphisms on
  the level of $K$-theory. 
  But these do \emph{not}, in general, induce a morphism of
  extensions --- in fact
  not even of the corresponding cyclic six term exact sequences. 
  Using Corollary~\ref{cor:ksixmappingcone}, 
  we easily see, that the diagram 
  $$\def\objectstyle{\scriptstyle} \def\labelstyle{\scriptstyle}
  \xymatrix{
    \ar[r]^-{\delta_1^{\sfS e}} &
    K_0(\sfS\A_0)\ar[r]^{K_0(\sfS\iota)}\ar[d]^\cong_{\alpha_0} &
    K_0(\sfS\A_1)\ar[r]^{K_0(\sfS\pi)}\ar[d]^\cong_{-\alpha_1} &
    K_0(\sfS\A_2)\ar[r]^{\delta_0^{\sfS e}}\ar[d]^\cong_{\alpha_2} &
    K_1(\sfS\A_0)\ar[r]^{K_1(\sfS\iota)}\ar[d]^\cong_{\alpha_3} &
    K_1(\sfS\A_1)\ar[r]^{K_1(\sfS\pi)}\ar[d]^\cong_{-\alpha_4} &
    K_1(\sfS\A_2)\ar[r]^-{\delta_1^{\sfS e}}\ar[d]^\cong_{\alpha_5} & \\
    \ar[r]^-{\delta_1^{\mc^3(e)}} &
    K_0(\sfS\sfC_\pi)\ar[r]^{K_0(\iota')} & 
    K_0(\sfC_{(\pi_{\mc})_{\mc}})\ar[r]^{K_0(\pi')} &
    K_0(\sfC_{\pi_\mc})\ar[r]^{\delta_0^{\mc^3(e)}} &
    K_1(\sfS\sfC_\pi)\ar[r]^{K_1(\iota')} &
    K_1(\sfC_{(\pi_{\mc})_{\mc}})\ar[r]^{K_1(\pi')} &
    K_1(\sfC_{\pi_\mc})\ar[r]^-{\delta_1^{\mc^3(e)}} &
    }$$
  commutes, where $\alpha_i$ are the induced maps as mentioned 
  above, and $\iota'$ and $\pi'$ denote the maps 
  $((\iota_{\mc})_{\mc})_{\mc}$ and 
  $((\pi_{\mc})_{\mc})_{\mc}$, respectively\footnote{Here we also use that the canonical
    identifications $K_i(\A_j)\rightarrow K_{1-i}(\sfS\A_j)$ give an
    isomorphism of the corresponding cyclic six term exact sequences.}.
  We expect that it is possible to find a functorial way to implement
  the $\kkE$-equivalences between $\mc^3(e)$ and $\sfS e$, but can not
  see how to do this --- not even how to make a canonical choice of
  $\kkE$-equivalences.  
\end{remark}

\begin{definition}
  The previous remark showed that $\mc^3(e)$ and $\sfS e$ are
  $\kkE$-equivalent (assuming the UCT), but the remark did not give 
   a canonical way to choose a
  specific $\kkE$-equivalence (so we get a functorial identification of
  the two functors). 

  For our purposes, it is enough to have the following lemma. 
  Let $e\colon\extwithmaps{\A_0}{\A_1}{\A_2}{\iota}{\pi}$ be a given
  extension of separable, nuclear \cas in \uct. 
  Assume, moreover, that $\Extsix(\ksix(e),\ksix(\sfS e))$ is the
  trivial group. 
  For each such extension $e$, 
  define 
  $$\mathbf{x}_e\in\kkE(\sfS e,\mc^3(e))$$
  to be the unique element
  inducing
  $(\alpha_0,-\alpha_1,\alpha_2,\alpha_3,-\alpha_4,\alpha_5)$ in
  $\Homsix(\ksix(\sfS e),\ksix(\mc^3 (e)))$ (as defined in the
  preceding remark).
\end{definition}

\begin{lemma}\label{lem-functorialbetweenSandmc3}
  Let $e$ and $e'$ 
  be two given
  extensions of separable, nuclear \cas in \uct. 
  Assume, moreover, that $\Extsix(\ksix(e),\ksix(\sfS e))$,
  $\Extsix(\ksix(e'),\ksix(\sfS e'))$, and
  $\Extsix(\ksix(e),\ksix(\sfS e'))$ are trivial groups. 
  Let $\phi$ be a morphism from $e$ to $e'$, and set
  $\mathbf{x}=\kkE(\phi)$ in $\kkE(e,e')$. 
  Then 
  $$\kkE(\sfS\phi)\times\mathbf{x}_{e'}
  =\widehat{\sfS}\mathbf{x}\times\mathbf{x}_{e'}
  =\mathbf{x}_{e}\times\widehat{\mc}^3(\mathbf{x})
  =\mathbf{x}_e\times\kkE(\mc^3(\phi)).$$
\end{lemma}
\begin{proof}
  From the assumptions and the UCT of Bonkat, we see that the canonical
  homomorphisms 
  $$\kkE(e,e')\longrightarrow\Homsix(\ksix(e),\ksix(e')),$$
  $$\kkE(\sfS e,\sfS e')\longrightarrow\Homsix(\ksix(\sfS e),\ksix(\sfS e')),$$
  $$\kkE(\mc^3(e),\mc^3(e'))\longrightarrow\Homsix(\ksix(\mc^3(e)),\ksix(\mc^3(e'))),$$
  $$\kkE(\sfS e,\mc^3(e'))\longrightarrow\Homsix(\ksix(\sfS e),\ksix(\mc^3(e')))$$
  are functorial isomorphisms.  
  Consequently, it is enough to prove that the result holds for the
  induced maps in $K$-theory, \ie,
  $$\ksix(\mathbf{x}_{e'})\circ\ksix(\widehat{\sfS}\mathbf{x})
  =\ksix(\widehat{\mc}^3(\mathbf{x}))\circ\ksix(\mathbf{x}_{e}).$$
   Again to prove this, it is enough to show that 
  $$\psi_i'\circ\sfS \phi_i=(\mc^3(\phi))_i\circ\psi_i,$$
  for $i=0,1,2$, where $\psi_0$ ($\psi_1$, and $\psi_2$, respectively) is the canonical
  \starhomos from the ideal (the extension, and the quotient, respectively)
  of $\sfS e$ to the ideal (the extension, and the quotient, respectively) of
  $\mc^3(e)$ --- and correspondingly for $\psi_i'$. 
  This equation is straightforward to check. 
\end{proof}

\begin{remark}\label{rem-mcsixismczero}
  Let $e$ 
  be an extension of separable, nuclear \cas in \uct, and assume
  that $\Extsix(\ksix(e),\ksix(\sfS e))$ 
  is the trivial group. 
  Then there exists a $\kkE$-equivalence
  $\widehat{\sfS}\mathbf{x}_{e}\times\mathbf{x}_{\mc^3(e)}$ 
  from $\sfS\sfS e$ to $\mc^6(e)$. 
  Composed with the standard
  $\kkE$-equivalence from $e$ to $\sfS\sfS e$, this gives a canonical
  $\kkE$-equivalence from $e$ to $\mc^6(e)$.

  \label{rem-mcvsxe}
  It is also easy to show that 
  $$\mathbf{x}_{\mc(e)}=-\kkE(\theta_e)\times\widehat{\mc}(\mathbf{x}_e).$$
\end{remark}

\begin{definition}
  For an extension $e$, let $\mathbf{b}_e$ denote the element of
  $\kkE(e,\sfS\sfS e)$ induced by the Bott element --- this is a $\kkE$-equivalence. 
  Moreover, let $\mathbf{z}_n$ denote the $\kkE$-equivalence in
  $\kkE(\sfS\FFFF{1}{0},\mathfrak{i}(\eo))$ induced by the canonical
  embedding $\C\rightarrow\Matn$.  Let $\mathbf{w}_n$ denote the $\kkE$-equivalence from
  $\extwithmaps{0}{\sfS\Matn}{\sfS\Matn}{0}{\ident}$ to
  $\mathfrak{q}(\eii)$ induced by the canonical embedding
  $\sfS\Matn\rightarrow\Ii$. 
\end{definition}

For each $n\in\N_{\geq 2}$, we will, during the following three
definitions, define 36 homomorphisms, 
$$\xymatrix@R=12pt{
  F_{1,i+1}\ar[r]^{h_{n,i}^{1,1,in}}  &
  H_{n,i}  \ar[r]^{h_{n,i}^{1,1,out}} &
  F_{1,i+3} \\ 
  F_{n,i}  \ar[r]^{h_{n,i}^{n,1,in}}  &
  H_{n,i}  \ar[r]^{h_{n,i}^{n,1,out}} &
  F_{1,i+2} \\
  F_{1,i+2}\ar[r]^{h_{n,i}^{1,n,in}}  &
  H_{n,i}  \ar[r]^{h_{n,i}^{1,n,out}} &
  F_{n,i+1}, 
  }$$
where we identify indices modulo 6 
(so we write \eg\ $h_{n,6}^{*,*,*}=h_{n,0}^{*,*,*}$).  

\begin{definition}
  For each $n\in\N_{\geq 2}$,
  there exists a short exact sequence
  $\extwithmaps{\mathfrak{i}(\eo)}{\eo}{\mathfrak{q}(\eo)}{\mathfrak{i}_{\eo}}{\mathfrak{q}_{\eo}}$ of
  extensions.  We define $h_{n,0}^{1,1,in}$ and $h_{n,0}^{1,1,out}$ by 
  $$\xymatrix@C=48pt{
    F_{1,1}\ar[rr]^{h_{n,0}^{1,1,in}}\ar@{=}[d] &&
    H_{n,0}\ar[rr]^{h_{n,0}^{1,1,out}}\ar@{=}[d] &&
    F_{1,3}\ar@{=}[d] \\ 
    \kkE(\mathfrak{q}(\eo),e)\ar[rr]_{\kkE(\mathfrak{q}_{\eo})\times -} &&
    \kkE(\eo,e)\ar[rr]_{\mathbf{x}_{\FFFF{1}{0}}^{-1}\times\mathbf{z}_n\times\kkE(\mathfrak{i}_{\eo})\times -} &&
    \kkE(\FFFF{1}{3},e).}$$
  By applying the functor $\widehat{\mc}$,
  we define $h_{n,i}^{1,1,in}$
  and $h_{n,i}^{1,1,out}$, for $i=1,2,3,4,5$, \ie, 
  \begin{align*}
    h_{n,i}^{1,1,in}
    &=\kkE(\widehat{\mc}^i(\kkE(\mathfrak{q}_{\eo})),e),\\
    h_{n,i}^{1,1,out}
    &=\kkE(\widehat{\mc}^i(\mathbf{x}_{\FFFF{1}{0}}^{-1}\times\mathbf{z}_n\times\kkE(\mathfrak{i}_{\eo})),e),
  \end{align*}
  for all $i=0,1,2,3,4,5$ (of course we use the canonical
  $\kkE$-equivalences from Remark~\ref{rem-mcsixismczero} to identify $\kkE(\FFFF{1}{j+6},e)$
  with $\kkE(\FFFF{1}{j},e)$). 
\end{definition}

\begin{definition}
  For each $n\in\N_{\geq 2}$,
  there exists a short exact sequence
  $\extwithmaps{\mathfrak{i}(\ei)}{\ei}{\mathfrak{q}(\ei)}{\mathfrak{i}_{\ei}}{\mathfrak{q}_{\ei}}$ of
  extensions.  Define $h_{n,1}^{n,1,in}$ and $h_{n,1}^{n,1,out}$ by 
  $$\xymatrix@C=48pt{
    F_{n,1}\ar[rr]^{h_{n,1}^{n,1,in}}\ar@{=}[d] &&
    H_{n,1}\ar[rr]^{h_{n,1}^{n,1,out}}\ar@{=}[d] &&
    F_{1,3}\ar@{=}[d] \\ 
    \kkE(\mathfrak{q}(\ei),e)\ar[rr]_{\kkE(\mathfrak{q}_{\ei})\times -} &&
    \kkE(\ei,e)\ar[rr]_{\mathbf{x}_{\FFFF{1}{0}}^{-1}\times\kkE(\mathfrak{i}_{\ei})\times -} &&
    \kkE(\FFFF{1}{3},e).}$$
  By applying the functor $\widehat{\mc}$, we
  define $h_{n,i}^{n,1,in}$
  and $h_{n,i}^{n,1,out}$, for $i=0,2,3,4,5$, \ie, 
  \begin{align*}
    h_{n,i}^{n,1,in}
    &=\kkE(\widehat{\mc}^{i-1}(\kkE(\mathfrak{q}_{\ei})),e),\\
    h_{n,i}^{n,1,out}
    &=\kkE(\widehat{\mc}^{i-1}(\mathbf{x}_{\FFFF{1}{0}}^{-1}\times\kkE(\mathfrak{i}_{\ei})),e),
  \end{align*}
  for all $i=1,2,3,4,5,6$.
\end{definition}

\begin{definition}
  For each $n\in\N_{\geq 2}$, 
  there exists a short exact sequence
  $\extwithmaps{\mathfrak{i}(\eii)}{\eii}{\mathfrak{q}(\eii)}{\mathfrak{i}_{\eii}}{\mathfrak{q}_{\eii}}$ of
  extensions.  Define 
  $h_{n,2}^{1,n,in}$ and $h_{n,2}^{1,n,out}$ by 
  $$\xymatrix@C=48pt{
    F_{1,4}\ar[rr]^{h_{n,2}^{1,n,in}}\ar@{=}[d] &&
    H_{n,2}\ar[rr]^{h_{n,2}^{1,n,out}}\ar@{=}[d] &&
    F_{n,3}\ar@{=}[d] \\ 
    \kkE(\FFFF{1}{4},e)\ar[rr]_{\kkE(\mathfrak{q}_{\eii})\times\mathbf{w}_n^{-1}\times\widehat{\mc}(\mathbf{z}_n^{-1})\times\mathbf{x}_{\FFFF{1}{1}}\times -} &&
    \kkE(\eii,e)\ar[rr]_{\mathbf{x}_{\FFFF{n}{0}}^{-1}\times\kkE(\mathfrak{i}_{\eii})\times -} &&
    \kkE(\FFFF{n}{3},e).}$$
   By applying the functor $\widehat{\mc}$, we define $h_{n,i}^{1,n,in}$
  and $h_{n,i}^{1,n,out}$, for $i=0,1,3,4,5$, \ie, 
  \begin{align*}
    h_{n,i}^{1,n,in}
    &=\kkE(\widehat{\mc}^{i-2}(\kkE(\mathfrak{q}_{\eii})\times\mathbf{w}_n^{-1}\times\widehat{\mc}(\mathbf{z}_n^{-1})\times\mathbf{x}_{\FFFF{1}{1}}),e),\\
    h_{n,i}^{1,n,out}
    &=\kkE(\widehat{\mc}^{i-2}(\mathbf{x}_{\FFFF{n}{0}}^{-1}\times\kkE(\mathfrak{i}_{\eii})),e),
  \end{align*}
  for all $i=2,3,4,5,6,7$. 
\end{definition}

\begin{definition}
  Now, we define homomorphisms $f_{n,i}$ from $F_{n,i}$ to $F_{n,i+1}$,
  for all $n\in\N$ and $i=0,1,2,3,4,5$. 
  We set 
  $$\xymatrix{
    \kkE(\sfS\sfS\FFFF{n}{0},e)\ar[rr] & &
    \kkE(\FFFF{n}{1},e)\ar[rr] & &
    \kkE(\FFFF{n}{2},e)\ar[ddd] \\
    &F_{n,0}\ar[r]^{f_{n,0}}\ar[ul]^-\cong_-{\mathbf{b}_{\FFFF{n}{0}}^{-1}\times -} &
                                    F_{n,1}\ar[r]^{f_{n,1}}\ar@{=}[u] &
    F_{n,2}\ar[d]^{f_{n,2}}\ar@{=}[ur] \\
    &F_{n,5}\ar[u]^{f_{n,5}}\ar[dl]_-\cong^-{\mathbf{x}_{\FFFF{n}{2}}\times -} &
    F_{n,4}\ar[l]^{f_{n,4}}\ar[d]_-\cong^-{\mathbf{x}_{\FFFF{n}{1}}\times -} &
    F_{n,3}\ar[l]^{f_{n,3}}\ar[dr]_-\cong^-{\mathbf{x}_{\FFFF{n}{0}}\times -} \\
    \kkE(\sfS\FFFF{n}{2},e)\ar[uuu] & &
    \kkE(\sfS\FFFF{n}{1},e)\ar[ll]         & &
    \kkE(\sfS\FFFF{n}{0},e)\ar[ll]        
    }$$
  where the outer sequence is the
  cyclic six term exact sequence in $\kkE$-theory induced by the short
  exact sequence
  $\extwithoutmaps{\mathfrak{i}(\FFFF{n}{2})}{\FFFF{n}{2}}{\mathfrak{q}(\FFFF{n}{2})}$ 
  (which is exactly 
  $\extwithoutmaps{\sfS\FFFF{n}{0}}{\FFFF{n}{2}}{\FFFF{n}{1}}$). 
\end{definition}

\begin{definition}
  \label{def-Bocksteinoperations}
  Now, we will define the Bockstein operations, 
  $$\xymatrix{F_{1,i}\ar[r]^{\rho_{n,i}}&F_{n,i}\ar[r]^{\beta_{n,i}}&F_{1,i+3}},$$
  for all $n\in\N_{\geq 2}$ and $i=0,1,2,3,4,5$. 

  The extension $\eo\colon\extwithoutmaps{\sfS\Matn}{\Io}{\C}$ induces a
  short exact sequence
  $\extwithmaps{\mathfrak{i}(\eo)}{\Fo}{\FFFF{1}{0}}{x}{y}$. 
  We set 
  $$\xymatrix@C=48pt{
    F_{1,0}\ar[rr]^{\rho_{n,0}}\ar@{=}[d]  & &
    F_{n,0}\ar[rr]^{\beta_{n,0}}\ar@{=}[d] & &
    F_{1,3}\ar@{=}[d] \\
    \kkE(\FFFF{1}{0},e)\ar[rr]_-{\kkE(y)\times -} &&
    \kkE(\FFFF{n}{0},e)\ar[rr]_-{\mathbf{x}_{\FFFF{1}{0}}^{-1}\times\mathbf{z}_n\times\kkE(x)\times -} &&
    \kkE(\FFFF{1}{3},e).    
    }$$
By applying the functor $\widehat{\mc}$,
 we define $\rho_{n,i}$ and
  $\beta_{n,i}$, for $i=1,2,3,4,5$, \ie, 
  \begin{align*}
    \rho_{n,i}
    &=\kkE(\widehat{\mc}^{i}(\kkE(y)),e), \\
    \beta_{n,i}
    &=\kkE(\widehat{\mc}^{i}(\mathbf{x}_{\FFFF{1}{0}}^{-1}\times\mathbf{z}_n\times\kkE(x)),e),
  \end{align*}
  for all $i=0,1,2,3,4,5$ (of course we use the canonical
  $\kkE$-equivalences from Remark~\ref{rem-mcsixismczero} to make
  identifications modulo 6). 
\end{definition}

\begin{definition}
  For each $n\in\N$,
  set 
  $\tilde{f}_{n,i}=f_{n,i}$ for $i=1,2,4,5$ and 
  $\tilde{f}_{n,i}=-f_{n,i}$ for $i=0,3$. 
\end{definition} 

\begin{theorem}\label{thm-diagrams}
  For all $n\in\N$ and all $i=0,1,2,3,4,5$, 
  $$\xymatrix{
    F_{n,i-1}\ar[r]^{f_{n,i-1}} & 
    F_{n,i}\ar[r]^{f_{n,i}} & 
    F_{n,i+1} 
    }$$
  is exact. 
  For all $n\in\N_{\geq 2}$ and all $i=0,1,2,3,4,5$, 
  $$\xymatrix@C=36pt{
    F_{1,i+1}\ar[r]^{h_{n,i}^{1,1,in}} &
    H_{n,i}\ar[r]^{h_{n,i}^{1,1,out}} &
    F_{1,i+3}\ar[d]_{nf_{1,i+3}}\\
    F_{1,i}\ar[u]_{nf_{1,i}} &
    H_{n,i+3}\ar[l]^{h_{n,i+3}^{1,1,out}} &
    F_{1,i+4}\ar[l]^{h_{n,i+3}^{1,1,in}}, 
    }\qquad
  \xymatrix@C=36pt{
    F_{n,i}\ar[r]^{h_{n,i}^{n,1,in}} &
    H_{n,i}\ar[r]^{h_{n,i}^{n,1,out}} &
    F_{1,i+2}\ar[d]_{f_{n,i+2}\circ\rho_{n,i+2}}\\
    F_{1,i+5}\ar[u]_{f_{n,i+5}\circ\rho_{n,i+5}} &
    H_{n,i+3}\ar[l]^{h_{n,i+3}^{n,1,out}} &
    F_{n,i+3}\ar[l]^{h_{n,i+3}^{n,1,in}}, 
    }$$
  $$\vcenter{\xymatrix@C=36pt{
    F_{1,i+2}\ar[r]^{h_{n,i}^{1,n,in}} &
    H_{n,i}\ar[r]^{h_{n,i}^{1,n,out}} &
    F_{n,i+1}\ar[d]_{\beta_{n,i+2}\circ f_{n,i+1}}\\
    F_{n,i+4}\ar[u]_{\beta_{n,i+5}\circ f_{n,i+4}} &
    H_{n,i+3}\ar[l]^{h_{n,i+3}^{1,n,out}} &
    F_{1,i+5}\ar[l]^{h_{n,i+3}^{1,n,in}}, 
    }}\quad\text{and}\quad
  \vcenter{\xymatrix@C=36pt{
    F_{1,i}\ar[r]^{\rho_{n,i}} & 
    F_{n,i}\ar[r]^{\beta_{n,i}} & 
    F_{1,i+3}\ar[d]_{\times n} \\
    F_{1,i}\ar[u]_{\times n} &
    F_{n,i+3}\ar[l]^{\beta_{n,i+3}} & 
    F_{1,i+3}\ar[l]^{\rho_{n,i+3}}  
    }}$$
  are exact, and, moreover, all the three diagrams 
  \begin{equation}
    \label{eq:thmdiagramone}
    \vcenter{\xymatrix@C=32pt@R=32pt{
      F_{1,i}\ar[d]_{\rho_{n,i}}\ar[r]^{\tilde{f}_{1,i}} &
      F_{1,i+1}\ar[d]_{\tildeh_{n,i}^{1,1,in}}\ar[dr]^{\tilde{f}_{1,i+1}} \\
      F_{n,i}\ar[dr]_{\beta_{n,i}}\ar[r]^{\tildeh_{n,i}^{n,1,in}} &
      H_{n,i}\ar[d]^{\tildeh_{n,i}^{1,1,out}}\ar[r]_{\tildeh_{n,i}^{n,1,out}} &
      F_{1,i+2}\ar[d]^{\tilde{f}_{1,i+2}} \\
      & F_{1,i+3}\ar[r]_{\times n} &
      F_{1,i+3}
      }}
  \end{equation}
  \begin{equation}
    \label{eq:thmdiagramtwo}
    \vcenter{\xymatrix@C=32pt@R=32pt{
      F_{1,i+1}\ar[d]_{\tilde{f}_{1,i+1}}\ar[r]^{\times n} &
      F_{1,i+1}\ar[d]_{\tildeh_{n,i}^{1,1,in}}\ar[dr]^{\rho_{n,i+1}} \\
      F_{1,i+2}\ar[dr]_{\tilde{f}_{1,i+2}}\ar[r]^{\tildeh_{n,i}^{1,n,in}} &
      H_{n,i}\ar[d]^{\tildeh_{n,i}^{1,1,out}}\ar[r]_{\tildeh_{n,i}^{1,n,out}} &
      F_{n,i+1}\ar[d]^{-\beta_{n,i+1}} \\
      & F_{1,i+3}\ar[r]_{\tilde{f}_{1,i+3}} &
      F_{1,i+4}
      }}
  \end{equation}
  \begin{equation}
    \label{eq:thmdiagramthree}
    \vcenter{\xymatrix@C=32pt@R=32pt{
      F_{n,i+5}\ar[d]_{-\beta_{n,i+5}}\ar[r]^{\tilde{f}_{n,i+5}} &
      F_{n,i}\ar[d]_{\tildeh_{n,i}^{n,1,in}}\ar[dr]^{\tilde{f}_{n,i}} \\
      F_{1,i+2}\ar[dr]_{\times n}\ar[r]^{\tildeh_{n,i}^{1,n,in}} &
      H_{n,i}\ar[d]^{\tildeh_{n,i}^{n,1,out}}\ar[r]_{\tildeh_{n,i}^{1,n,out}} &
      F_{n,i+1}\ar[d]^{\tilde{f}_{n,i+1}} \\
      & F_{1,i+2}\ar[r]_{\rho_{n,i+2}} &
      F_{n,i+2}
      }}
  \end{equation}
  commute. 
\end{theorem}
\begin{proof}
  See next section. 
\end{proof}

\begin{corollary}
  For each $n\in\N_{\geq 2}$ and $i=0,1,2,3,4,5$, the two squares 
  $$\vcenter{\xymatrix{
      F_{1,i}\ar[d]_{\rho_{n,i}}\ar[r]^{\tilde{f}_{1,i}} &
      F_{1,i+1}\ar[d]^{\rho_{n,i+1}} \\
      F_{n,i}\ar[r]_{\tilde{f}_{n,i}} &
      F_{n,i+1}
      }}\quad\text{and}\quad
  \vcenter{\xymatrix{
      F_{n,i}\ar[d]_{-\beta_{n,i}}\ar[r]^{\tilde{f}_{n,i}} &
      F_{n,i+1}\ar[d]^{\beta_{n,i+1}} \\
      F_{1,i+3}\ar[r]_{\tilde{f}_{1,i+3}} &
      F_{1,i+4} }}$$
  commute.
\end{corollary}
\begin{proof}
  This follows directly from the previous theorem: 
  \begin{align*}
    \tilde{f}_{n,i}\circ\rho_{n,i}
    &=\tilde{f}_{n,i}\circ \tildeh_{n,i-1}^{1,n,out}\circ
    \tildeh_{n,i-1}^{1,1,in}&\text{by }\eqref{eq:thmdiagramtwo}\\
    &=\rho_{n,i+1}\circ \tildeh_{n,i-1}^{n,1,out}\circ
    \tildeh_{n,i-1}^{1,1,in}&\text{by }\eqref{eq:thmdiagramthree} \\
    &=\rho_{n,i+1}\circ \tilde{f}_{1,i}
    &\text{by }\eqref{eq:thmdiagramone} \\
    \beta_{n,i+1}\circ \tilde{f}_{n,i}
    &=\beta_{n,i+1}\circ \tildeh_{n,i}^{1,n,out}\circ \tildeh_{n,i}^{n,1,in}
    &\text{by }\eqref{eq:thmdiagramthree} \\
    &=-\tilde{f}_{1,i+3}\circ \tildeh_{n,i}^{1,1,out}\circ \tildeh_{n,i}^{n,1,in}
    &\text{by }\eqref{eq:thmdiagramtwo} \\
    &=-\tilde{f}_{1,i+3}\circ\beta_{n,i}
    &\text{by }\eqref{eq:thmdiagramone}
    &\qedhere
  \end{align*}
\end{proof}

\begin{remark}
  From the preceding theorem and corollary, it follows that,
  for each $n\in\N_{\geq 2}$ and $i=0,1,2$, we have the following 
  --- both horizontally and vertically six term cyclic ---
  commuting diagrams with exact rows and columns: 
  \begin{equation}
    \def\objectstyle{\scriptstyle} \def\labelstyle{\scriptstyle}
    \tag{$D_0$}\label{eq:diagram0} 
    \vcenter{\xymatrix@!=20pt@C=36pt@R=18pt{
        &\ar[d]^{\times n}&\ar[d]^{\times n}&\ar[d]^{\times n}
        &\ar[d]^{\times n}&\ar[d]^{\times n}&\ar[d]^{\times n} \\
        \ar[r]^{\tilde{f}_{1,5}} &
        F_{1,0}\ar[d]^{\rho_{n,0}}\ar[r]^{\tilde{f}_{1,0}} &
        F_{1,1}\ar[d]^{\rho_{n,1}}\ar[r]^{\tilde{f}_{1,1}} &
        F_{1,2}\ar[d]^{\rho_{n,2}}\ar[r]^{\tilde{f}_{1,2}} &
        F_{1,3}\ar[d]^{\rho_{n,3}}\ar[r]^{\tilde{f}_{1,3}} &
        F_{1,4}\ar[d]^{\rho_{n,4}}\ar[r]^{\tilde{f}_{1,4}} &
        F_{1,5}\ar[d]^{\rho_{n,5}}\ar[r]^{\tilde{f}_{1,5}} & \\
        \ar[r]^{\tilde{f}_{n,5}} &
        F_{n,0}\ar[d]^{\beta_{n,0}}\ar[r]^{\tilde{f}_{n,0}} &
        F_{n,1}\ar[d]^{-\beta_{n,1}}\ar[r]^{\tilde{f}_{n,1}} &
        F_{n,2}\ar[d]^{\beta_{n,2}}\ar[r]^{\tilde{f}_{n,2}} &
        F_{n,3}\ar[d]^{-\beta_{n,3}}\ar[r]^{\tilde{f}_{n,3}} &
        F_{n,4}\ar[d]^{\beta_{n,4}}\ar[r]^{\tilde{f}_{n,4}} &
        F_{n,5}\ar[d]^{-\beta_{n,5}}\ar[r]^{\tilde{f}_{n,5}} & \\
        \ar[r]^{\tilde{f}_{1,2}} &
        F_{1,3}\ar[d]^{\times n}\ar[r]^{\tilde{f}_{1,3}} &
        F_{1,4}\ar[d]^{\times n}\ar[r]^{\tilde{f}_{1,4}} &
        F_{1,5}\ar[d]^{\times n}\ar[r]^{\tilde{f}_{1,5}} &
        F_{1,0}\ar[d]^{\times n}\ar[r]^{\tilde{f}_{1,0}} &
        F_{1,1}\ar[d]^{\times n}\ar[r]^{\tilde{f}_{1,1}} &
        F_{1,2}\ar[d]^{\times n}\ar[r]^{\tilde{f}_{1,2}} & \\
        &&&&&&
        }}
  \end{equation}
  \begin{equation}
    \def\objectstyle{\scriptstyle} \def\labelstyle{\scriptstyle}
    \tag{$D_{1}$}\label{eq:diagram1-1} 
    \vcenter{\xymatrix@!=20pt@C=36pt@R=18pt{
        &\ar[d]^{\times n}&\ar[d]^{n\tilde{f}_{1,i}}&\ar[d]&\ar[d]^{\times n}&\ar[d]^{n\tilde{f}_{1,i+3}}&\ar[d] \\
        \ar[r]^{\tilde{f}_{1,i+5}} & 
        F_{1,i}  \ar[d]^{\rho_{n,i}}\ar[r]^{\tilde{f}_{1,i}}&
        F_{1,i+1}\ar[d]^{\tildeh_{n,i}^{1,1,in}}\ar[r]^{\tilde{f}_{1,i+1}} &
        F_{1,i+2}\ar@{=}[d]\ar[r]^{\tilde{f}_{1,i+2}} &
        F_{1,i+3}\ar[d]^{\rho_{n,i+3}}\ar[r]^{\tilde{f}_{1,i+3}} &
        F_{1,i+4}\ar[d]^{\tildeh_{n,i+3}^{1,1,in}}\ar[r]^{\tilde{f}_{1,i+4}} &
        F_{1,i+5}\ar@{=}[d]\ar[r]^{\tilde{f}_{1,i+5}} & \\
        \ar[r]^{\rho_{n,i}\circ \tilde{f}_{1,i+5}} &
        F_{n,i}\ar[d]^{\beta_{n,i}}\ar[r]^{\tildeh_{n,i}^{n,1,in}} &
        H_{n,i}\ar[d]^{\tildeh_{n,i}^{1,1,out}}\ar[r]^{\tildeh_{n,i}^{n,1,out}} &
        F_{1,i+2}\ar[d]\ar[r]^{\rho_{n,i+3}\circ \tilde{f}_{1,i+2}} &
        F_{n,i+3}\ar[d]^{\beta_{n,i+3}}\ar[r]^{\tildeh_{n,i+3}^{n,1,in}} &
        H_{n,i+3}\ar[d]^{\tildeh_{n,i+3}^{1,1,out}}\ar[r]^{\tildeh_{n,i+3}^{n,1,out}} &
        F_{1,i+5}\ar[d]\ar[r]^{\rho_{n,i}\circ \tilde{f}_{1,i+5}} & \\
        \ar[r] &
        F_{1,i+3}\ar[d]^{\times n}\ar@{=}[r] &
        F_{1,i+3}\ar[d]^{n\tilde{f}_{1,i+3}}\ar[r] &
        0 \ar[d]\ar[r] &
        F_{1,i}\ar[d]^{\times n}\ar@{=}[r] &
        F_{1,i}\ar[d]^{n\tilde{f}_{1,i}}\ar[r] &
        0 \ar[d]\ar[r] & \\
        &&&&&&
        }} 
    \end{equation}
    \begin{equation}
    \def\objectstyle{\scriptstyle} \def\labelstyle{\scriptstyle}
    \tag{$D_{1}^\star$}\label{eq:diagram1-2} 
    \vcenter{\xymatrix@!=20pt@C=36pt@R=18pt{
        &\ar[d]&\ar[d]^{n\tilde{f}_{1,i}}&\ar[d]^{\tilde{f}_{1,i}}&\ar[d]&\ar[d]^{n\tilde{f}_{1,i+3}}&\ar[d]^{\tilde{f}_{1,i+3}} \\
        \ar[r] & 
        0  \ar[d]\ar[r]&
        F_{1,i+1}\ar[d]^{\tildeh_{n,i}^{1,1,in}}\ar@{=}[r] &
        F_{1,i+1}\ar[d]^{\tilde{f}_{1,i+1}}\ar[r] &
        0\ar[d]\ar[r] &
        F_{1,i+4}\ar[d]^{\tildeh_{n,i+3}^{1,1,in}}\ar@{=}[r] &
        F_{1,i+4}\ar[d]^{\tilde{f}_{1,i+4}}\ar[r] & \\
        \ar[r]^{\rho_{n,i}\circ \tilde{f}_{1,i+5}} &
        F_{n,i}\ar@{=}[d]\ar[r]^{\tildeh_{n,i}^{n,1,in}} &
        H_{n,i}\ar[d]^{\tildeh_{n,i}^{1,1,out}}\ar[r]^{\tildeh_{n,i}^{n,1,out}} &
        F_{1,i+2}\ar[d]^{\tilde{f}_{1,i+2}}\ar[r]^{\rho_{n,i+3}\circ \tilde{f}_{1,i+2}} &
        F_{n,i+3}\ar@{=}[d]\ar[r]^{\tildeh_{n,i+3}^{n,1,in}} &
        H_{n,i+3}\ar[d]^{\tildeh_{n,i+3}^{1,1,out}}\ar[r]^{\tildeh_{n,i+3}^{n,1,out}} &
        F_{1,i+5}\ar[d]^{\tilde{f}_{1,i+5}}\ar[r]^{\rho_{n,i}\circ \tilde{f}_{1,i+5}} & \\
        \ar[r]^{\rho_{n,i}} &
        F_{n,i}\ar[d]\ar[r]^{\beta_{n,i}} &
        F_{1,i+3}\ar[d]^{n\tilde{f}_{1,i+3}}\ar[r]^{\times n} &
        F_{1,i+3}\ar[d]^{\tilde{f}_{1,i+3}}\ar[r]^{\rho_{n,i+3}} &
        F_{n,i+3}\ar[d]\ar[r]^{\beta_{n,i+3}} &
        F_{1,i}\ar[d]^{n\tilde{f}_{1,i}}\ar[r]^{\times n} &
        F_{1,i} \ar[d]^{\tilde{f}_{1,i}}\ar[r]^{\rho_{n,i}} & \\
        &&&&&&
        }} 
  \end{equation}
  \begin{equation}
    \def\objectstyle{\scriptstyle} \def\labelstyle{\scriptstyle}
    \tag{$D_{2}$}\label{eq:diagram2-1} 
    \vcenter{\xymatrix@!=20pt@C=36pt@R=18pt{
        &\ar[d]^{\tilde{f}_{1,i}}&\ar[d]^{n\tilde{f}_{1,i}}&\ar[d]&\ar[d]^{\tilde{f}_{1,i+3}}&\ar[d]^{n\tilde{f}_{1,i+3}}&\ar[d] \\
        \ar[r]^{\beta_{n,i+4}} & 
        F_{1,i+1}\ar[d]^{\tilde{f}_{1,i+1}}\ar[r]^{\times n} &
        F_{1,i+1}\ar[d]^{\tildeh_{n,i}^{1,1,in}}\ar[r]^{\rho_{n,i+1}} &
        F_{n,i+1}\ar@{=}[d]\ar[r]^{\beta_{n,i+1}} &
        F_{1,i+4}\ar[d]^{\tilde{f}_{1,i+4}}\ar[r]^{\times n} &
        F_{1,i+4}\ar[d]^{\tildeh_{n,i+3}^{1,1,in}}\ar[r]^{\rho_{n,i+4}} &
        F_{n,i+4}\ar@{=}[d]\ar[r]^{\beta_{n,i+4}} & \\
        \ar[r]^{\tilde{f}_{1,i+1}\circ\beta_{n,i+4}} &
        F_{1,i+2}\ar[d]^{\tilde{f}_{1,i+2}}\ar[r]^{\tildeh_{n,i}^{1,n,in}} &
        H_{n,i}\ar[d]^{\tildeh_{n,i}^{1,1,out}}\ar[r]^{\tildeh_{n,i}^{1,n,out}} &
        F_{n,i+1}\ar[d]\ar[r]^{\tilde{f}_{1,i+4}\circ\beta_{n,i+1}} &
        F_{1,i+5}\ar[d]^{\tilde{f}_{1,i+5}}\ar[r]^{\tildeh_{n,i+3}^{1,n,in}} &
        H_{n,i+3}\ar[d]^{\tildeh_{n,i+3}^{1,1,out}}\ar[r]^{\tildeh_{n,i+3}^{1,n,out}} &
        F_{n,i+4}\ar[d]\ar[r]^{\tilde{f}_{1,i+1}\circ\beta_{n,i+4}} & \\
        \ar[r] &
        F_{1,i+3}\ar[d]^{\tilde{f}_{1,i+3}}\ar@{=}[r] &
        F_{1,i+3}\ar[d]^{n\tilde{f}_{1,i+3}}\ar[r] &
        0\ar[d]\ar[r] &
        F_{1,i}\ar[d]^{\tilde{f}_{1,i}}\ar@{=}[r] &
        F_{1,i}\ar[d]^{n\tilde{f}_{1,i}}\ar[r] &
        0\ar[d]\ar[r] & \\
        &&&&&&
        }} 
  \end{equation}
  \begin{equation}
    \def\objectstyle{\scriptstyle} \def\labelstyle{\scriptstyle}
    \tag{$D_{2}^\star$}\label{eq:diagram2-2} 
    \vcenter{\xymatrix@!=20pt@C=36pt@R=18pt{
        &\ar[d]&\ar[d]^{n\tilde{f}_{1,i}}&\ar[d]^{\times n}&\ar[d]&\ar[d]^{n\tilde{f}_{1,i+3}}&\ar[d]^{\times n} \\
        \ar[r] & 
        0\ar[d]\ar[r] &
        F_{1,i+1}\ar[d]^{\tildeh_{n,i}^{1,1,in}}\ar@{=}[r] &
        F_{1,i+1}\ar[d]^{\rho_{n,i+1}}\ar[r] &
        0\ar[d]\ar[r] &
        F_{1,i+4}\ar[d]^{\tildeh_{n,i+3}^{1,1,in}}\ar@{=}[r] &
        F_{1,i+4}\ar[d]^{\rho_{n,i+4}}\ar[r] & \\
        \ar[r]^{-\tilde{f}_{1,i+1}\circ\beta_{n,i+4}} &
        F_{1,i+2}\ar@{=}[d]\ar[r]^{\tildeh_{n,i}^{1,n,in}} &
        H_{n,i}\ar[d]^{\tildeh_{n,i}^{1,1,out}}\ar[r]^{\tildeh_{n,i}^{1,n,out}} &
        F_{n,i+1}\ar[d]^{-\beta_{n,i+1}}\ar[r]^{-\tilde{f}_{1,i+4}\circ\beta_{n,i+1}} &
        F_{1,i+5}\ar@{=}[d]\ar[r]^{\tildeh_{n,i+3}^{1,n,in}} &
        H_{n,i+3}\ar[d]^{\tildeh_{n,i+3}^{1,1,out}}\ar[r]^{\tildeh_{n,i+3}^{1,n,out}} &
        F_{n,i+4}\ar[d]^{-\beta_{n,i+4}}\ar[r]^{-\tilde{f}_{1,i+1}\circ\beta_{n,i+4}} & \\
        \ar[r]^{\tilde{f}_{1,i+1}} &
        F_{1,i+2}\ar[d]\ar[r]^{\tilde{f}_{1,i+2}} &
        F_{1,i+3}\ar[d]^{n\tilde{f}_{1,i+3}}\ar[r]^{\tilde{f}_{1,i+3}} &
        F_{1,i+4}\ar[d]^{\times n}\ar[r]^{\tilde{f}_{1,i+4}} &
        F_{1,i+5}\ar[d]\ar[r]^{\tilde{f}_{1,i+5}} &
        F_{1,i}\ar[d]^{n\tilde{f}_{1,i}}\ar[r]^{\tilde{f}_{1,i}} &
        F_{1,i+1}\ar[d]^{\times n}\ar[r]^{\tilde{f}_{1,i+1}} & \\
        &&&&&&
        }} 
  \end{equation}
  \begin{equation}
    \def\objectstyle{\scriptstyle} \def\labelstyle{\scriptstyle}
    \tag{$D_{3}$}\label{eq:diagram3-1} 
    \vcenter{\xymatrix@!=20pt@C=36pt@R=18pt{
        &\ar[d]^{\rho_{n,i+5}}&\ar[d]^{\tilde{f}_{n,i+5}\circ\rho_{n,i+5}}&\ar[d]&
        \ar[d]^{\rho_{n,i+2}}&\ar[d]^{\tilde{f}_{n,i+2}\circ\rho_{n,i+2}}&\ar[d] \\
        \ar[r]^{\tilde{f}_{n,i+4}} &
        F_{n,i+5}\ar[d]^{-\beta_{n,i+5}}\ar[r]^{\tilde{f}_{n,i+5}} &
        F_{n,i}\ar[d]^{\tildeh_{n,i}^{n,1,in}}\ar[r]^{\tilde{f}_{n,i}} &
        F_{n,i+1}\ar@{=}[d]\ar[r]^{\tilde{f}_{n,i+1}} &
        F_{n,i+2}\ar[d]^{-\beta_{n,i+2}}\ar[r]^{\tilde{f}_{n,i+2}} &
        F_{n,i+3}\ar[d]^{\tildeh_{n,i+3}^{n,1,in}}\ar[r]^{\tilde{f}_{n,i+3}} &
        F_{n,i+4}\ar@{=}[d]\ar[r]^{\tilde{f}_{n,i+4}} & \\
        \ar[r] &
        F_{1,i+2}\ar[d]^{\times n}\ar[r]^{\tildeh_{n,i}^{1,n,in}} &
        H_{n,i}\ar[d]^{\tildeh_{n,i}^{n,1,out}}\ar[r]^{\tildeh_{n,i}^{1,n,out}} & 
        F_{n,i+1}\ar[d]\ar[r] &
        F_{1,i+5}\ar[d]^{\times n}\ar[r]^{\tildeh_{n,i+3}^{1,n,in}} &
        H_{n,i+3}\ar[d]^{\tildeh_{n,i+3}^{n,1,out}}\ar[r]^{\tildeh_{n,i+3}^{1,n,out}} & 
        F_{n,i+4}\ar[d]\ar[r] & \\
        \ar[r] &
        F_{1,i+2}\ar[d]^{\rho_{n,i+2}}\ar@{=}[r] & 
        F_{1,i+2}\ar[d]^{\tilde{f}_{n,i+2}\circ\rho_{n,i+2}}\ar[r] & 
        0\ar[d]\ar[r] & 
        F_{1,i+5}\ar[d]^{\rho_{n,i+5}}\ar@{=}[r] & 
        F_{1,i+5}\ar[d]^{\tilde{f}_{n,i+5}\circ\rho_{n,i+5}}\ar[r] & 
        0\ar[d]\ar[r] & \\
        &&&&&&
        }} 
  \end{equation}
  \begin{equation}
    \def\objectstyle{\scriptstyle} \def\labelstyle{\scriptstyle}
    \tag{$D_{3}^\star$}\label{eq:diagram3-2} 
    \vcenter{\xymatrix@!=20pt@C=36pt@R=18pt{
        &\ar[d]&\ar[d]^{\tilde{f}_{n,i+5}\circ\rho_{n,i+5}}&\ar[d]^{\tilde{f}_{n,i+5}}&
        \ar[d]&\ar[d]^{\tilde{f}_{n,i+2}\circ\rho_{n,i+2}}&\ar[d]^{\tilde{f}_{n,i+2}} \\
        \ar[r] &
        0\ar[d]\ar[r] &
        F_{n,i}\ar[d]^{\tildeh_{n,i}^{n,1,in}}\ar@{=}[r] &
        F_{n,i}\ar[d]^{\tilde{f}_{n,i}}\ar[r] &
        0\ar[d]\ar[r] &
        F_{n,i+3}\ar[d]^{\tildeh_{n,i+3}^{n,1,in}}\ar@{=}[r] &
        F_{n,i+3}\ar[d]^{\tilde{f}_{n,i+3}}\ar[r] & \\
        \ar[r] &
        F_{1,i+2}\ar@{=}[d]\ar[r]^{\tildeh_{n,i}^{1,n,in}} &
        H_{n,i}\ar[d]^{\tildeh_{n,i}^{n,1,out}}\ar[r]^{\tildeh_{n,i}^{1,n,out}} & 
        F_{n,i+1}\ar[d]^{\tilde{f}_{n,i+1}}\ar[r] &
        F_{1,i+5}\ar@{=}[d]\ar[r]^{\tildeh_{n,i+3}^{1,n,in}} &
        H_{n,i+3}\ar[d]^{\tildeh_{n,i+3}^{n,1,out}}\ar[r]^{\tildeh_{n,i+3}^{1,n,out}} & 
        F_{n,i+4}\ar[d]^{\tilde{f}_{n,i+4}}\ar[r] & \\
        \ar[r]^{\beta_{n,i+5}} &
        F_{1,i+2}\ar[d]\ar[r]^{\times n} & 
        F_{1,i+2}\ar[d]^{\tilde{f}_{n,i+2}\circ\rho_{n,i+2}}\ar[r]^{\rho_{n,i+2}} & 
        F_{n,i+2}\ar[d]^{\tilde{f}_{n,i+2}}\ar[r]^{\beta_{n,i+2}} & 
        F_{1,i+5}\ar[d]\ar[r]^{\times n} & 
        F_{1,i+5}\ar[d]^{\tilde{f}_{n,i+5}\circ\rho_{n,i+5}}\ar[r]^{\rho_{n,i+5}} & 
        F_{n,i+5}\ar[d]^{\tilde{f}_{n,i+5}}\ar[r]^{\beta_{n,i+5}} & \\
        &&&&&&
        }}
  \end{equation}
\end{remark}

\begin{remark}
  Note that Diagrams~$(D_i)$ and~$(D_i^\star)$ with two extra conditions each are
  equivalent, for $i=1,2,3$. 
  For instance, Diagram~$(D_1)$ with the extra condition 
  $$n h^{1,1,out}_{n,i}=\tilde{f}_{1,i+2}\circ h^{n,1,out}_{n,i},\qquad
  n h^{1,1,out}_{n,i+3}=\tilde{f}_{1,i+5}\circ h^{n,1,out}_{n,i+3}$$
  is equivalent to Diagram~$(D_1^\star)$ with the extra condition 
  $$h^{n,1,in}_{n,i}\circ \rho_{n,i}=h^{1,1,in}_{n,i}\circ
  \tilde{f}_{1,i},\qquad 
  h^{n,1,in}_{n,i+3}\circ \rho_{n,i+3}=h^{1,1,in}_{n,i+3}\circ
  \tilde{f}_{1,i+3}.$$
\end{remark}

\section{Proof of Theorem~\protect\protectreference}

The purpose of this section is to prove Theorem~\ref{thm-diagrams}. 
First we need some results, which will be useful in the proof.

\begin{remark}
  Let $\A$ be a separable, nuclear \ca in \uct. 
  Set $e_0\colon\extwithmaps{\A}{\A}{0}{\ident}{}$, and set
  $e_i=\mc^i(e_0)$. 
  As earlier we know that 
  $$e_0\colon\extwithmaps{\A}{\A}{0}{\ident}{},$$
  $$e_1\colon\extwithmaps{0}{\A}{\A}{}{\ident},$$
  $$e_2\colon\extwithmaps{\sfS\A}{\sfC\A}{\A}{\iota}{\ev_1},$$
  $$e_3\colon\extwithmaps{\sfS\A}{\sfC\A\oplus_{\ev_1,\ev_1}\sfC\A}{\sfC\A}{(0,\iota)}{\pi_1},$$
  where $\pi_1$ is the projection onto the first coordinate. 

  Note that there exists a canonical morphism,  $\phi=(\ident,(0,\iota),0)$, 
  from $\sfS e_0$ to $e_3$, which induces
  a $\kkE$-equivalence. 
  It is evident that $\kkE(\phi)$ is exactly
  $\mathbf{x}_{e_0}$ in the case that 
  $\Ext_\Z^1(K_i(\A),K_{1-i}(\A))=0$, for $i=0,1$. 
  Also note that in this case, $\kkE(\mc(\phi))=-\mathbf{x}_{e_1}$
  (according to Remark~\ref{rem-mcvsxe}). 

  Note that $\mathfrak{i}(e_2)=\sfS e_0$, $\mathfrak{q}(e_2)=e_1$, and
  $\mc(\mathfrak{i}(e_2))=\sfS e_1$. 
  So if applying $\mc^0$, $\mc^1$, and $\mc^2$ to the short exact sequence 
  $\extwithmaps{\mathfrak{i}(e_2)}{e_2}{\mathfrak{q}(e_2)}{\mathfrak{i}_{e_2}}{\mathfrak{q}_{e_2}}$, we get just 
  $\extwithmaps{\sfS e_0}{e_2}{e_1}{\mathfrak{i}_{e_2}}{\mathfrak{q}_{e_2}}$,
  $\extwithmaps{\sfS
    e_1}{e_3}{e_2}{\mc(\mathfrak{i}_{e_2})}{\mc(\mathfrak{q}_{e_2})}$,
  and $\extwithmaps{\mc\sfS
    e_1}{e_4}{e_3}{\mc^2(\mathfrak{i}_{e_2})}{\mc^2(\mathfrak{q}_{e_2})}$, respectively. 
\end{remark}

\begin{proposition}\label{prop-rotate}
  Let \A be a separable, nuclear \ca in \uct satisfying
  $\Ext^1_\Z(K_i(\A),K_{1-i}(\A))=0$, for $i=0,1$, and 
  let $e$ be an extension of separable \cas. 
  Set $e_0\colon\extwithmaps{\A}{\A}{0}{\ident}{}$, and set
  $e_i=\mc^i(e_0)$. 
  Then we have 
  $$\xymatrix{
    \kkE(\sfS\sfS e_1,e)\ar[dr]_\cong^{\mathbf{b}_{e_1}\times -}\ar[rr] & \ar@{}[d]^{\text{Anti-commutes}} &
    \kkE(e_2,e)\ar@{=}[d]\ar[rr] & \ar@{}[d]_{\text{Commutes}} &
    \kkE(e_3,e)\ar[dl]_\cong^{\mathbf{x}_{e_0}\times -}\ar[ddd]_{\text{Commutes}\quad} \\
    &\kkE(e_1,e)\ar[r]&
    \kkE(e_2,e)\ar[r]&
    \kkE(\sfS e_0,e)\ar[d]\\
    &\kkE(\sfS\sfS e_0,e)\ar[u]&
    \kkE(\sfS e_2,e)\ar[l]&
    \kkE(\sfS e_1,e)\ar[l]\\
    \kkE(\sfS e_3,e)
    \ar[ur]_\cong^{\widehat{\sfS}\mathbf{x}_{e_0}\times -}\ar[uuu]_{\quad\text{Commutes}} & \ar@{}[u]_{\text{Commutes}} &
    \kkE(\sfS e_2,e)\ar@{=}[u]\ar[ll] & \ar@{}[u]^{\text{Anti-commutes}} &
    \kkE(\sfS e_1,e)\ar@{=}[ul]\ar[ll]     
    }$$
  where the inner and outer sequences are the cyclic
  six term exact sequences in $\kkE$-theory induced by
  $\extwithmaps{\mathfrak{i}(e_2)}{e_2}{\mathfrak{q}(e_2)}{\mathfrak{i}_{e_2}}{\mathfrak{q}_{e_2}}$ and 
  $\extwithmaps{\mc(\mathfrak{i}(e_2))}{\mc(e_2)}{\mc(\mathfrak{q}(e_2))}{\mc(\mathfrak{i}_{e_2})}{\mc(\mathfrak{q}_{e_2})}$, 
  respectively.  Moreover, we have that
  $$\xymatrix{
    \kkE(\sfS\mc\sfS e_1,e)
    \ar[dr]_\cong^{\quad\mathbf{b}_{e_2}\times\kkE(\sfS\theta_{e_1})\times -}\ar[rr] & 
    \ar@{}[d]^{\text{Anti-commutes}} &
    \kkE(e_3,e)\ar@{=}[d]\ar[rr] & \ar@{}[d]_{\text{Anti-commutes}} &
    \kkE(e_4,e)\ar[dl]_\cong^{\mathbf{x}_{e_1}\times -}\ar[ddd]_{\text{Anti-commutes}\quad} \\
    &\kkE(e_2,e)\ar[r]&
    \kkE(e_3,e)\ar[r]&
    \kkE(\sfS e_1,e)\ar[d]\\
    &\kkE(\sfS\sfS e_1,e)\ar[u]&
    \kkE(\sfS e_3,e)\ar[l]&
    \kkE(\sfS e_2,e)\ar[l]\\
    \kkE(\sfS e_4,e)
    \ar[ur]_\cong^{\widehat{\sfS}\mathbf{x}_{e_1}\times -}\ar[uuu]_{\quad\text{Anti-commutes}} & \ar@{}[u]_{\text{Anti-commutes}} &
    \kkE(\sfS e_3,e)\ar@{=}[u]\ar[ll] & \ar@{}[u]^{\text{Anti-commutes}} &
    \kkE(\mc\sfS e_1,e)\ar[ul]_\cong^{\kkE(\theta_{e_1})}\ar[ll]     
    }$$
  where the inner and outer sequences are the cyclic
  six term exact sequences in $\kkE$-theory induced by
  $\extwithmaps{\mc(\mathfrak{i}(e_2))}{\mc(e_2)}{\mc(\mathfrak{q}(e_2))}{\mc(\mathfrak{i}_{e_2})}{\mc(\mathfrak{q}_{e_2})}$ and
  $\extwithmaps{\mc^2(\mathfrak{i}(e_2))}{\mc^2(e_2)}{\mc^2(\mathfrak{q}(e_2))}{\mc^2(\mathfrak{i}_{e_2})}{\mc^2(\mathfrak{q}_{e_2})}$, 
  respectively.
\end{proposition}
\begin{proof}
  First, writing out the short exact sequences   
  $\extwithmaps{\sfS e_0}{e_2}{e_1}{\mathfrak{i}_{e_2}}{\mathfrak{q}_{e_2}}$,
  $\extwithmaps{\sfS e_1}{e_3}{e_2}{\mc(\mathfrak{i}_{e_2})}{\mc(\mathfrak{q}_{e_2})}$, and 
  $\extwithmaps{\mc\sfS e_1}{e_4}{e_3}{\mc^2(\mathfrak{i}_{e_2})}{\mc^2(\mathfrak{q}_{e_2})}$ we get: 
  $$\xymatrix{
    \sfS\A\ar@{=}[d]\ar@{=}[r] &
    \sfS\A\ar@{^(->}[d]^{\iota}\ar[r] &
    0\ar@{^(->}[d] \\
    \sfS\A\ar[d]\ar@{^(->}[r]^\iota &
    \sfC\A\ar@{->>}[d]^\pi\ar@{->>}[r]^\pi &
    \A\ar@{=}[d] \\
    0\ar[r] &
    \A\ar@{=}[r] &
    \A
    }\qquad
  \xymatrix{
    0\ar[d]\ar[r] &
    \sfS\A\ar@{^(->}[d]^{\iota_2}\ar@{=}[r] &
    \sfS\A\ar@{^(->}[d]^\iota \\
    \sfS\A\ar@{=}[d]\ar@{^(->}[r]^-{\iota_1} &
    \sfC\A\oplus_{\pi,\pi}\sfC\A\ar@{->>}[d]^{\pi_1}\ar@{->>}[r]^-{\pi_2} &
    \sfC\A\ar@{->>}[d]^{\pi} \\
    \sfS\A\ar@{^(->}[r]^{\iota} &
    \sfC\A\ar@{->>}[r]^\pi & 
    \A
    }$$
  $$\xymatrix@C=60pt{
    \sfS\sfS\A
    \ar@{^(->}[d]^{\iota_{\sfS\A}}
    \ar@{^(->}[r]^{\sfS\iota} &
    \sfS\sfC\A
    \ar@{^(->}[d]^{(0,0,\iota_{\sfC\A})}
    \ar@{->>}[r]^{\sfS\pi} & 
    \sfS\A
    \ar@{^(->}[d]^{\iota_2} \\
    \sfC\sfS\A
    \ar@{->>}[d]^{\ev_1}
    \ar@{^(->}[r]^-{(\ev_1,0,\sfC\iota)} &
    (\sfC\A\oplus_{\pi,\pi}\sfC\A)\oplus_{\pi_1,\ev_1}\sfC\sfC\A
    \ar@{->>}[d]^{(f,g,h)\mapsto(f,g)}
    \ar@{->>}[r]^-{(f,g,h)\mapsto(g,h(-,1))} &
    \sfC\A\oplus_{\pi,\pi}\sfC\A
    \ar@{->>}[d]^{\pi_1} \\
    \sfS\A
    \ar@{^(->}[r]^{\iota_1} &
    \sfC\A\oplus_{\pi,\pi}\sfC\A
    \ar@{->>}[r]^{\pi_2} &
    \sfC\A
    }$$
  Now, we write out the cyclic six term exact sequences of cyclic six term
  exact sequences corresponding to these three short exact sequences--- where we horizontally use the \kk-boundary maps and vertically
  use the $\ksix$-boundary maps. 
  For convenience we will identify $K_1$ with $K_0\circ\sfS$. 
  Moreover, we let
  \begin{align*}
    \tilde{\A} &=\sfC\A\oplus_{\pi,\pi}\sfC\A \\
    \tilde{\tilde{\A}}&=\tilde{\A}\oplus_{\pi_1,\ev_1}\sfC\sfC\A .
  \end{align*}
  The diagrams are:    
  \begin{equation*}
    \def\objectstyle{\scriptstyle} 
    \def\labelstyle{\scriptstyle}
    \xymatrix@!=20pt@C=42pt@R=14pt{
      &\ar[d]&\ar[d]^{\ident}&\ar[d]&\ar[d]&\ar[d]^{-\ident}&\ar[d] \\
      \ar[r] & 
      K_0(\sfS\A)\ar@{=}[d]\ar@{=}[r] & 
      K_0(\sfS\A)\ar[d]^{K_0(\iota)}\ar[r] & 
      0 \ar[d]\ar[r] &
      K_0(\sfS\sfS\A)\ar@{=}[d]\ar@{=}[r] & 
      K_0(\sfS\sfS\A)\ar[d]^{K_0(\sfS\iota)}\ar[r] & 
      0 \ar[d]\ar[r] & \\
      \ar[r]^{-\ident} &
      K_0(\sfS\A)\ar[d]\ar[r]^{K_0(\iota)} &
      0\ar[d]^{K_0(\pi)}\ar[r]^{K_0(\pi)} &
      K_0(\A)\ar@{=}[d]\ar[r]^{-\beta_\A} &
      K_0(\sfS\sfS\A)\ar[d]\ar[r]^{K_0(\sfS\iota)} &
      0\ar[d]^{K_0(\sfS\pi)}\ar[r]^{K_0(\sfS\pi)} &
      K_0(\sfS\A)\ar@{=}[d]\ar[r]^{-\ident} & \\ 
      \ar[r] & 
      0 \ar[d]\ar[r] &
      K_0(\A)\ar[d]^{-\beta_\A}\ar@{=}[r] &
      K_0(\A)\ar[d]\ar[r] &
      0 \ar[d]\ar[r] &
      K_0(\sfS\A)\ar[d]^{\beta_{\sfS\A}}\ar@{=}[r] &
      K_0(\sfS\A)\ar[d]\ar[r] & \\
      \ar[r] &
      K_0(\sfS\sfS\A)\ar@{=}[d]\ar@{=}[r] & 
      K_0(\sfS\sfS\A)\ar[d]^{K_0(\sfS\iota)}\ar[r] & 
      0 \ar[d]\ar[r] &
      K_0(\sfS\sfS\sfS\A)\ar@{=}[d]\ar@{=}[r] & 
      K_0(\sfS\sfS\sfS\A)\ar[d]^{K_0(\sfS\sfS\iota)}\ar[r] & 
      0 \ar[d]\ar[r] & \\
      \ar[r]^{-\ident} &
      K_0(\sfS\sfS\A)\ar[d]\ar[r]^{K_0(\sfS\iota)} &
      0\ar[d]^{K_0(\sfS\pi)}\ar[r]^{K_0(\sfS\pi)} &
      K_0(\sfS\A)\ar@{=}[d]\ar[r]^{-\beta_{\sfS\A}} &
      K_0(\sfS\sfS\sfS\A)\ar[d]\ar[r]^{K_0(\sfS\sfS\iota)} &
      0\ar[d]^{K_0(\sfS\sfS\pi)}\ar[r]^{K_0(\sfS\sfS\pi)} &
      K_0(\sfS\sfS\A)\ar@{=}[d]\ar[r]^{-\ident} & \\ 
      \ar[r] & 
      0 \ar[d]\ar[r] &
      K_0(\sfS\A)\ar[d]^{\ident}\ar@{=}[r] &
      K_0(\sfS\A)\ar[d]\ar[r] &
      0 \ar[d]\ar[r] &
      K_0(\sfS\sfS\A)\ar[d]^{-\ident}\ar@{=}[r] &
      K_0(\sfS\sfS\A)\ar[d]\ar[r] & \\
      &&&&&&
      }
  \end{equation*}
  \begin{equation*}
    \def\objectstyle{\scriptstyle} 
    \def\labelstyle{\scriptstyle}
    \xymatrix{
      &\ar[d]&\ar[d]&\ar[d]^{\ident}&\ar[d]&\ar[d]&\ar[d]^{-\ident} \\
      \ar[r] &
      0\ar[d]\ar[r] &
      K_0(\sfS\A)\ar[d]^{K_0(\iota_2)}\ar@{=}[r] &
      K_0(\sfS\A)\ar[d]^{K_0(\iota)}\ar[r] &
      0\ar[d]\ar[r] &
      K_0(\sfS^2\A)\ar[d]^{K_0(\sfS\iota_2)}\ar@{=}[r] &
      K_0(\sfS^2\A)\ar[d]^{K_0(\sfS\iota)}\ar[r] & \\
      \ar[r] &
      K_0(\sfS\A)\ar@{=}[d]\ar[r]^{K_0(\iota_1)} &
      K_0(\tilde{\A})\ar[d]^{K_0(\pi_1)}\ar[r]^{K_0(\pi_2)} &
      0\ar[d]^{K_0(\pi)}\ar[r] &
      K_0(\sfS^2\A)\ar@{=}[d]\ar[r]^{K_0(\sfS\iota_1)} &
      K_0(\sfS\tilde{A})\ar[d]^{K_0(\sfS\pi_1)}\ar[r]^{K_0(\sfS\pi_2)} &
      0\ar[d]^{K_0(\sfS\pi)}\ar[r] & \\
      \ar[r]^{-\ident} &
      K_0(\sfS\A)\ar[d]\ar[r]^{K_0(\iota)} &
      0\ar[d]\ar[r]^{K_0(\pi)} &
      K_0(\A)\ar[d]^{-\beta_\A}\ar[r]^{-\beta_\A}&
      K_0(\sfS^2\A)\ar[d]\ar[r]^{K_0(\sfS\iota)} &
      0\ar[d]\ar[r]^{K_0(\sfS\pi)} &
      K_0(\sfS\A)\ar[d]^{\beta_{\sfS\A}}\ar[r]^{-\ident}& \\
      \ar[r] &
      0\ar[d]\ar[r] &
      K_0(\sfS^2\A)\ar[d]^{K_0(\sfS\iota_2)}\ar@{=}[r] &
      K_0(\sfS^2\A)\ar[d]^{K_0(\sfS\iota)}\ar[r] &
      0\ar[d]\ar[r] &
      K_0(\sfS^3\A)\ar[d]^{K_0(\sfS^2\iota_2)}\ar@{=}[r] &
      K_0(\sfS^3\A)\ar[d]^{K_0(\sfS^2\iota)}\ar[r] & \\
      \ar[r] &
      K_0(\sfS^2\A)\ar@{=}[d]\ar[r]^{K_0(\sfS\iota_1)} &
      K_0(\sfS\tilde{\A})\ar[d]^{K_0(\sfS\pi_1)}\ar[r]^{K_0(\sfS\pi_2)} &
      0\ar[d]^{K_0(\sfS\pi)}\ar[r] &
      K_0(\sfS^3\A)\ar@{=}[d]\ar[r]^{K_0(\sfS^2\iota_1)} &
      K_0(\sfS^2\tilde{\A})\ar[d]^{K_0(\sfS^2\pi_1)}\ar[r]^{K_0(\sfS^2\pi_2)} &
      0\ar[d]^{K_0(\sfS^2\pi)}\ar[r] & \\
      \ar[r]^{-\ident} &
      K_0(\sfS^2\A)\ar[d]\ar[r]^{K_0(\sfS\iota)} &
      0\ar[d]\ar[r]^{K_0(\sfS\pi)} &
      K_0(\sfS\A)\ar[d]^{\ident}\ar[r]^{-\beta_{\sfS\A}}&
      K_0(\sfS^3\A)\ar[d]\ar[r]^{K_0(\sfS^2\iota)} &
      0\ar[d]\ar[r]^{K_0(\sfS^2\pi)} &
      K_0(\sfS^2\A)\ar[d]^{-\ident}\ar[r]^{-\ident}& \\
      &&&&&&
    }
  \end{equation*}
  \begin{equation*}
    \def\objectstyle{\scriptstyle} 
    \def\labelstyle{\scriptstyle}
    \xymatrix{
      &\ar[d]^{\ident}&\ar[d]&\ar[d]&\ar[d]^{-\ident}&\ar[d]&\ar[d] \\
      \ar[r]^{\ident} &
      K_0(\sfS^2\A)\ar[d]\ar[r] &
      0\ar[d]\ar[r] &
      K_0(\sfS\A)\ar[d]^{K_0(\iota_2)}\ar[r]^{\beta_{\sfS\A}} &
      K_0(\sfS^3\A)\ar[d]\ar[r] &
      0\ar[d]\ar[r] &
      K_0(\sfS^2\A)\ar[d]^{K_0(\sfS\iota_2)}\ar[r]^{\ident} & \\
      \ar[r] &
      0\ar[d]\ar[r] &
      K_0(\tilde{\tilde{\A}})\ar[d]\ar[r] &
      K_0(\tilde{\A})\ar[d]\ar[r] &
      0\ar[d]\ar[r] &
      K_0(\sfS(\tilde{\tilde{\A}}))\ar[d]\ar[r] &
      K_0(\sfS\tilde{\A})\ar[d]\ar[r] & \\
      \ar[r] &
      K_0(\sfS\A)\ar[d]^{-\beta_{\sfS\A}}\ar[r]^{K_0(\iota_1)} &
      K_0(\tilde{\A})\ar[d]\ar[r] &
      0\ar[d]\ar[r] &
      K_0(\sfS^2\A)\ar[d]^{\beta_{\sfS^2\A}}\ar[r]^{K_0(\sfS\iota_1)} &
      K_0(\sfS\tilde{\A})\ar[d]\ar[r] &
      0\ar[d]\ar[r] & \\
      \ar[r]^{\ident} &
      K_0(\sfS^3\A)\ar[d]\ar[r] &
      0\ar[d]\ar[r] &
      K_0(\sfS^2\A)\ar[d]^{K_0(\sfS\iota_2)}\ar[r]^{\beta_{\sfS^2\A}} &
      K_0(\sfS^4\A)\ar[d]\ar[r] &
      0\ar[d]\ar[r] &
      K_0(\sfS^3\A)\ar[d]^{K_0(\sfS^2\iota_2)}\ar[r]^{\ident} & \\
      \ar[r] &
      0\ar[d]\ar[r] &
      K_0(\sfS(\tilde{\tilde{\A}}))\ar[d]\ar[r] &
      K_0(\sfS\tilde{\A})\ar[d]\ar[r] &
      0\ar[d]\ar[r] &
      K_0(\sfS^2(\tilde{\tilde{\A}}))\ar[d]\ar[r] &
      K_0(\sfS^2\tilde{\A})\ar[d]\ar[r] & \\
      \ar[r] &
      K_0(\sfS^2\A)\ar[d]^{\ident}\ar[r]^{K_0(\sfS\iota_1)} &
      K_0(\sfS\tilde{\A})\ar[d]\ar[r] &
      0\ar[d]\ar[r] &
      K_0(\sfS^3\A)\ar[d]^{-\ident}\ar[r]^{K_0(\sfS^2\iota_1)} &
      K_0(\sfS^2\tilde{\A})\ar[d]\ar[r] &
      0\ar[d]\ar[r] & \\
      &&&&&&}
  \end{equation*}

  Note that $\mathbf{x}_{e_0}$, $-\mathbf{x}_{e_1}$, and $\kkE(\theta_{e_1})$ are induced
  by the morphisms 
  $$\xymatrix{
    \sfS\A\ar@{=}[d]\ar@{=}[r] & 
    \sfS\A\ar@{^(->}[d]^{\iota_2} &
    0\ar[d]\ar[r] & \sfS\sfC\A\ar@{^(->}[d] & 
    \sfS\sfS\A\ar@{^(->}[d]\ar[r]^{\flip} & 
    \sfS\sfS\A\ar@{^(->}[d] \\
    \sfS\A\ar[d]\ar[r]^-{\iota_2} & 
    \sfC\A\oplus_{\pi,\pi}\sfC\A\ar[d]^{\pi_1} &
    \sfS\A\ar[d]\ar[r]^-{(0,\iota,0)} &
    (\sfC\A\oplus_{\pi,\pi}\sfC\A)\oplus_{\pi_1,\ev_1}\sfC\sfC\A\ar@{->>}[d] &
    \sfC\sfS\A\ar@{->>}[d]\ar[r]^{\flip} & 
    \sfS\sfC\A\ar@{->>}[d] \\
    0\ar[r] & \sfC\A &
    \sfS\A\ar[r]^-{\iota_2} & \sfC\A\oplus_{\pi,\pi}\sfC\A & 
    \sfS\A\ar@{=}[r] & \sfS\A }$$
  
  Using all these diagrams, a long, tedious, straightforward verification
  shows the proposition. 
\end{proof}

\begin{remark}
  What we actually showed in the proof of the preceding proposition,
  is that the corresponding diagrams of morphisms in the category
  $\catkkE$ (\ie, before we apply $\kkE(-,e)$) 
  commute respectively anti-commute. 
  This observation will be useful in the sequel. 
\end{remark}

\begin{proof}[Proof of the first part of Theorem~\ref{thm-diagrams}]
  By definition, 
  $\xymatrix{F_{n,i-1}\ar[r]^{f_{n,i-1}}&F_{n,i}\ar[r]^{f_{n,i}}&F_{n,i+1}}$ 
  is exact for all $n\in\N$ and all $i=0,1,2,3,4,5$. 

  Note that there exists a commuting square 
  $$\xymatrix{
    \C\ar@{=}[d]\ar@{=}[r] &
    \C\ar[d] \\
    \C\ar[r] & 
    \Matn,
    }$$
  where the maps $\C\rightarrow\Matn$ are the unital \starhomos. 
  By naturality of the mapping cone construction, this induces a
  morphism $\phi=(\phi_0,\phi_1,\phi_2)$ from $\FFFF{1}{2}$ to $\eo$. 
  This gives a commuting diagram 
  $$\xymatrix{
    \mathfrak{i}(\FFFF{1}{2}) \ar[d]^{(\phi_0,\phi_0,0)}\ar@{^(->}[r] &
    \FFFF{1}{2} \ar[d]^{\phi}\ar@{->>}[r] &
    \mathfrak{q}(\FFFF{1}{2}) \ar@{=}[d] \\
    \mathfrak{i}(\eo) \ar@{^(->}[r] &
    \eo \ar@{->>}[r] &
    \mathfrak{q}(\eo) 
    }$$
  of short exact sequences. 
  If we apply $\kkE(-,e)$ to this diagram we will get a morphism
  between two cyclic six term exact sequences in $\kkE$-theory. 
  Using the standard equivalences introduced so far, we arrive at the
  commuting diagrams
  \begin{equation*}
    \def\objectstyle{\scriptstyle} \def\labelstyle{\scriptstyle}
    \xymatrix@!=20pt@C=70pt@R=18pt{
      \ar[r]^-{f_{1,0}} & 
      \kkE(\FFFF{1}{1},e)\ar[r]^{f_{1,1}}\ar@{=}[d] & 
      \kkE(\FFFF{1}{2},e)\ar[r]^{f_{1,2}}\ar@{=}[d] & 
      \kkE(\FFFF{1}{3},e)\ar[r]^{f_{1,3}}\ar[d]^{\mathbf{x}_{\FFFF{1}{0}}\times -}_\cong & \\
      \ar[r] & 
      \kkE(\FFFF{1}{1},e)\ar[r]\ar@{=}[d] & 
      \kkE(\FFFF{1}{2},e)\ar[r]\ar@{<-}[d]^{\kkE(\phi)\times-} & 
      \kkE(\sfS\FFFF{1}{0},e)\ar[r]\ar@{<-}[d]^{\kkE((\phi_0,\phi_0,0))\times-} & \\
      \ar[r]&\kkE(\FFFF{1}{1},e)\ar[r]\ar@{=}[d] & 
      \kkE(\eo,e)\ar[r]\ar@{=}[d] &
      \kkE(\mathfrak{i}(\eo),e)\ar[r]\ar[d]^{\mathbf{x}_{\FFFF{1}{0}}^{-1}\times\mathbf{z}_n\times -}_\cong & \\
      &\kkE(\FFFF{1}{1},e)\ar[r]\ar@{=}[d] & 
      \kkE(\eo,e)\ar[r]\ar@{=}[d] &
      \kkE(\FFFF{1}{3},e)\ar@{=}[d] &  \\
      & 
      F_{1,1}\ar[r]^{h_{n,0}^{1,1,in}} & 
      H_{n,0}\ar[r]^{h_{n,0}^{1,1,out}} &
      F_{1,3} & \\
      \ar[r]^{f_{1,3}} &
      \kkE(\FFFF{1}{4},e)\ar[r]^{f_{1,4}}\ar[d]^{\mathbf{x}_{\FFFF{1}{1}}\times -}_\cong & 
      \kkE(\FFFF{1}{5},e)\ar[r]^{f_{1,5}}\ar[d]^{\mathbf{x}_{\FFFF{1}{2}}\times -}_\cong & 
      \kkE(\FFFF{1}{0},e)\ar[r]^-{f_{1,0}}\ar[d]^{\mathbf{b}^{-1}_{\FFFF{1}{0}}\times -}_\cong & \\
      \ar[r] &
      \kkE(\sfS\FFFF{1}{1},e)\ar[r]\ar@{=}[d] & 
      \kkE(\sfS\FFFF{1}{2},e)\ar[r]\ar@{<-}[d]^{\kkE(\sfS\phi)\times-} & 
      \kkE(\sfS\sfS\FFFF{1}{0},e)\ar[r]\ar@{<-}[d]^{\kkE(\sfS(\phi_0,\phi_0,0))\times-} & \\
      \ar[r] &
      \kkE(\sfS\FFFF{1}{1},e)\ar[r]\ar[d]^{\mathbf{x}_{\FFFF{1}{1}}^{-1}\times -}_\cong &
      \kkE(\sfS\eo,e)\ar[r]\ar[d]^{\mathbf{x}_{\e{n}{0}}^{-1}\times -}_\cong &
      \kkE(\sfS\mathfrak{i}(\eo),e)\ar[r]
      \ar[d]^{\mathbf{b}_{\FFFF{1}{0}}\times\widehat{\sfS}\mathbf{z}_n\times -}_\cong & \\
      &\kkE(\FFFF{1}{4},e)\ar[r]\ar@{=}[d] &
      \kkE(\e{n}{3},e)\ar[r]\ar@{=}[d] &
      \kkE(\FFFF{1}{0},e)\ar@{=}[d] &  \\
      &
      F_{1,4}\ar[r]^{h_{n,3}^{1,1,in}} & 
      H_{n,3}\ar[r]^{h_{n,3}^{1,1,out}} &
      F_{1,0} &
    }
  \end{equation*}
  with exact rows. 
  We use Lemma~\ref{lem-functorialbetweenSandmc3} for commutativity of
  the two squares between row three and four in the lower part of the
  diagram --- and 
 we use that 
  $$\Extsix(\ksix(\e{n}{0}),\ksix(\sfS \FFFF{1}{1}))=0,$$
  $$\Extsix(\ksix(\sfS\FFFF{1}{3}),\ksix(\sfS\e{n}{3}))=0.$$
  This is easily verified using projective resolutions. 
  
  It is easy to verify that, up to a sign, we have 
  $$\kkE(\mathbf{z}_n^{-1}\times\kkE((\phi_0,\phi_0,0)),e)=n\ident
  \quad\text{and}\quad
  \kkE(\widehat{\sfS}\mathbf{z}_n^{-1}\times\widehat{\sfS}\kkE((\phi_0,\phi_0,0)),e)=n\ident.$$
  Consequently, $nf_{1,0}$ and $nf_{1,3}$ are exactly the connecting
  homomorphisms of the cyclic six term exact sequence in the bottom. 

  This proves exactness of the first of the four cyclic sequences in
  the theorem, for $i=0,3$. 

  This same result also works for $i=1,2,4,5$, by invoking
  Proposition~\ref{prop-rotate} (recall that we do not care about
  the signs, because that does not change exactness). 

  Note that there exists a commuting square 
  $$\xymatrix{
    \Io\ar@{=}[d]\ar@{=}[r] &
    \Io\ar[d] \\
    \Io\ar[r] & 
    \C,
    }$$
  where the maps $\Io\rightarrow\C$ are the canonical surjective \starhomos. 
  By naturality of the mapping cone construction, this induces a
  morphism $\phi=(\phi_0,\phi_1,\phi_2)$ from $\FFFF{n}{2}$ to $\ei$. 
  This gives a commuting diagram 
  $$\xymatrix{
    \mathfrak{i}(\FFFF{n}{2}) \ar[d]^{(\phi_0,\phi_0,0)}\ar@{^(->}[r] &
    \FFFF{n}{2} \ar[d]^{\phi}\ar@{->>}[r] &
    \mathfrak{q}(\FFFF{n}{2}) \ar@{=}[d] \\
    \mathfrak{i}(\ei) \ar@{^(->}[r] &
    \ei \ar@{->>}[r] &
    \mathfrak{q}(\ei) 
    }$$
  of short exact sequences. 
  If we apply $\kkE(-,e)$ to this diagram we will get a morphism
  between two cyclic six term exact sequences in $\kkE$-theory. 
  Using the standard equivalences introduced so far, we arrive at the
  commuting diagram
  \begin{equation*}
    \def\objectstyle{\scriptstyle} \def\labelstyle{\scriptstyle}
    \xymatrix@!=20pt@C=70pt@R=18pt{
      \ar[r]^-{f_{n,0}} & 
      \kkE(\FFFF{n}{1},e)\ar[r]^{f_{n,1}}\ar@{=}[d] & 
      \kkE(\FFFF{n}{2},e)\ar[r]^{f_{n,2}}\ar@{=}[d] & 
      \kkE(\FFFF{n}{3},e)\ar[r]^{f_{n,3}}\ar[d]^{\mathbf{x}_{\FFFF{n}{0}}\times -}_\cong & \\
      \ar[r] & 
      \kkE(\FFFF{n}{1},e)\ar[r]\ar@{=}[d] & 
      \kkE(\FFFF{n}{2},e)\ar[r]\ar@{<-}[d]^{\kkE(\phi)\times-} & 
      \kkE(\sfS\FFFF{n}{0},e)\ar[r]\ar@{<-}[d]^{\kkE((\phi_0,\phi_0,0))\times-} & \\
      \ar[r]&\kkE(\FFFF{n}{1},e)\ar[r]\ar@{=}[d] & 
      \kkE(\ei,e)\ar[r]\ar@{=}[d] &
      \kkE(\sfS\FFFF{1}{0},e)\ar[r]\ar[d]^{\mathbf{x}_{\FFFF{1}{0}}^{-1}\times -}_\cong & \\
      &\kkE(\FFFF{n}{1},e)\ar[r]\ar@{=}[d] & 
      \kkE(\ei,e)\ar[r]\ar@{=}[d] &
      \kkE(\FFFF{1}{3},e)\ar@{=}[d] & \\
      & 
      F_{n,1}\ar[r]^{h_{n,1}^{n,1,in}} & 
      H_{n,1}\ar[r]^{h_{n,1}^{n,1,out}} &
      F_{1,3} & \\
      \ar[r]^{f_{n,3}} & 
      \kkE(\FFFF{n}{4},e)\ar[r]^{f_{n,4}}\ar[d]^{\mathbf{x}_{\FFFF{n}{1}}\times -}_\cong & 
      \kkE(\FFFF{n}{5},e)\ar[r]^{f_{n,5}}\ar[d]^{\mathbf{x}_{\FFFF{n}{2}}\times -}_\cong & 
      \kkE(\FFFF{n}{0},e)\ar[r]^-{f_{n,0}}\ar[d]^{\mathbf{b}^{-1}_{\FFFF{n}{0}}\times -}_\cong & \\
      \ar[r] & 
      \kkE(\sfS\FFFF{n}{1},e)\ar[r]\ar@{=}[d] & 
      \kkE(\sfS\FFFF{n}{2},e)\ar[r]\ar@{<-}[d]^{\kkE(\sfS\phi)\times-} & 
      \kkE(\sfS\sfS\FFFF{n}{0},e)\ar[r]\ar@{<-}[d]^{\kkE(\sfS(\phi_0,\phi_0,0))\times-} & \\
      \ar[r] &
      \kkE(\sfS\FFFF{n}{1},e)\ar[r]\ar[d]^{\mathbf{x}_{\FFFF{n}{1}}^{-1}\times -}_\cong &
      \kkE(\sfS\ei,e)\ar[r]\ar[d]^{\mathbf{x}_{\e{n}{1}}^{-1}\times -}_\cong &
      \kkE(\sfS\sfS\FFFF{1}{0},e)\ar[r]
      \ar[d]^{\mathbf{b}_{\FFFF{1}{0}}\times -}_\cong & \\
      & \kkE(\FFFF{n}{4},e)\ar[r]\ar@{=}[d] &
      \kkE(\e{n}{4},e)\ar[r]\ar@{=}[d] &
      \kkE(\FFFF{1}{0},e)\ar@{=}[d] &  \\
      & F_{n,4}\ar[r]^{h_{n,4}^{n,1,in}} & 
      H_{n,4}\ar[r]^{h_{n,4}^{n,1,out}} &
      F_{1,0} &
    }
  \end{equation*}
    with exact rows. 
  We use Lemma~\ref{lem-functorialbetweenSandmc3} for commutativity of
  the two squares between row three and four in the lower part of the diagram --- and we use that 
  $$\Extsix(\ksix(\sfS\e{n}{1}),\ksix(\sfS\FFFF{n}{4}))=0,$$
  $$\Extsix(\ksix(\sfS\FFFF{1}{3}),\ksix(\sfS\e{n}{4}))=0.$$
  This is easily verified using projective resolutions. 

  Using naturality of $\mathbf{b}_-$ and
  Lemma~\ref{lem-functorialbetweenSandmc3}, it is easy to see that 
  $$\kkE(\mathbf{x}^{-1}_{\FFFF{n}{0}}\times\kkE((\phi_0,\phi_0,0))\times\mathbf{x}_{\FFFF{1}{0}},e)
  =\rho_{n,3},\quad\text{and}$$
  $$\kkE(\mathbf{b}_{\FFFF{n}{0}}\times\kkE(\sfS(\phi_0,\phi_0,0))\times\mathbf{b}^{-1}_{\FFFF{1}{0}},e)
  =\rho_{n,0}.$$
  Consequently, $f_{n,3}\circ\rho_{n,3}$ and $f_{n,0}\circ\rho_{n,0}$ are exactly the connecting
  homomorphisms of the cyclic six term exact sequence in the bottom. 

  This proves exactness of the second of the four cyclic sequences in
  the theorem, for $i=1,4$. 

  This same result also works for $i=0,2,3,5$, by invoking
  Proposition~\ref{prop-rotate}. 

  Note that there exists a commuting square 
  $$\xymatrix{
    \Ii\ar[d]\ar[r] &
    \Io\ar@{=}[d] \\
    \Io\ar@{=}[r] & 
    \Io,
    }$$
  where the maps $\Ii\rightarrow\Io$ are the canonical surjective \starhomos. 
  By naturality of the mapping cone construction, this induces a
  morphism $\phi=(\phi_0,\phi_1,\phi_2)$ from $\eii$ to $\FFFF{n}{2}$. 
  This gives a commuting diagram 
  $$\xymatrix{
    \mathfrak{i}(\eii) \ar@{=}[d]\ar@{^(->}[r] &
    \eii\ar[d]^{\phi}\ar@{->>}[r] &
    \mathfrak{q}(\eii) \ar[d]^{(0,\phi_2,\phi_2)} \\
    \mathfrak{i}(\FFFF{n}{2}) \ar@{^(->}[r] &
    \FFFF{n}{2} \ar@{->>}[r] &
    \mathfrak{q}(\FFFF{n}{2}) 
    }$$
  of short exact sequences. 
  If we apply $\kkE(-,e)$ to this diagram we will get a morphism
  between two cyclic six term exact sequences in $\kkE$-theory. 
  Using the standard equivalences introduced so far, we arrive at the
  commuting diagram
  \begin{equation*}
    \def\objectstyle{\scriptstyle} 
    \def\labelstyle{\scriptstyle}
    \xymatrix@!=20pt@C=70pt@R=18pt{
      & 
      F_{1,4}\ar@{=}[d]\ar[r]^{h_{n,2}^{1,n,in}} &
      H_{n,2}\ar@{=}[d]\ar[r]^{h_{n,2}^{1,n,out}} &
      F_{n,3}\ar@{=}[d] & \\
      & 
      \kkE(\FFFF{1}{4},e)
      \ar[r]
      \ar[d]_{\mathbf{w}_n^{-1}\times\widehat{\mc}(\mathbf{z}_n^{-1})\times\mathbf{x}_{\FFFF{1}{1}}\times -} & 
      \kkE(\eii,e)\ar[r]\ar@{=}[d] & 
      \kkE(\FFFF{n}{3},e)\ar[d]^{\mathbf{x}_{\FFFF{n}{0}}\times-}_\cong & \\
      \ar[r] & 
      \kkE(\mathfrak{q}(\eii),e)\ar[r]\ar@{<-}[d]^{\kkE((0,\phi_2,\phi_2))\times-} & 
      \kkE(\eii,e)\ar[r]\ar@{<-}[d]^{\kkE(\phi)\times-} & 
      \kkE(\sfS\FFFF{n}{0},e)\ar[r]\ar@{=}[d] & \\
      \ar[r]&\kkE(\FFFF{n}{1},e)\ar[r]\ar@{=}[d] & 
      \kkE(\FFFF{n}{2},e)\ar[r]\ar@{=}[d] &
      \kkE(\sfS\FFFF{n}{0},e)\ar[r]\ar@{<-}[d]^{\mathbf{x}_{\FFFF{n}{0}}\times -}_\cong & \\
      \ar[r]^-{f_{n,0}}&\kkE(\FFFF{n}{1},e)\ar[r]^{f_{n,1}} &
      \kkE(\FFFF{n}{2},e)\ar[r]^{f_{n,2}} &
      \kkE(\FFFF{n}{3},e)\ar[r]^{f_{n,3}} &  \\
      &
      F_{1,1}\ar@{=}[d]\ar[r]^{h_{n,5}^{1,n,in}} &
      H_{n,5}\ar@{=}[d]\ar[r]^{h_{n,5}^{1,n,out}} &
      F_{n,0}\ar@{=}[d] &  \\
      &
      \kkE(\FFFF{1}{1},e)
      \ar[r]
      \ar[d]_{\widehat{\sfS}(\mathbf{w}_n^{-1}\times\widehat{\mc}(\mathbf{z}_n^{-1}))\times\mathbf{b}^{-1}_{\FFFF{1}{1}}\times -} & 
      \kkE(\e{n}{5},e)\ar[r]\ar[d]^{\mathbf{x}_{\eii}\times -}_\cong & 
      \kkE(\FFFF{n}{0},e)\ar[d]^{\mathbf{b}^{-1}_{\FFFF{n}{0}}\times -}_\cong & \\
      \ar[r] &
      \kkE(\sfS\mathfrak{q}(\eii),e)\ar[r]\ar@{<-}[d]^{\kkE(\sfS(0,\phi_2,\phi_2))\times-} & 
      \kkE(\sfS\eii,e)\ar[r]\ar@{<-}[d]^{\kkE(\sfS\phi)\times-} & 
      \kkE(\sfS\sfS\FFFF{n}{0},e)\ar[r]\ar@{=}[d] & \\
      \ar[r] &
      \kkE(\sfS\FFFF{n}{1},e)\ar[r]\ar@{<-}[d]^{\mathbf{x}_{\FFFF{n}{1}}\times -}_\cong &
      \kkE(\sfS\FFFF{n}{2},e)\ar[r]\ar@{<-}[d]^{\mathbf{x}_{\FFFF{n}{2}}\times -}_\cong &
      \kkE(\sfS\sfS\FFFF{n}{0},e)\ar[r]
      \ar@{<-}[d]^{\mathbf{b}_{\FFFF{n}{0}}^{-1}\times -}_\cong & \\
      \ar[r]^{f_{n,3}} &
      \kkE(\FFFF{n}{4},e)\ar[r]^{f_{n,4}} &
      \kkE(\FFFF{n}{5},e)\ar[r]^{f_{n,5}} &
      \kkE(\FFFF{n}{0},e)\ar[r]^-{f_{n,0}} &  
    }
  \end{equation*}
  with exact rows. 
  We use Lemma~\ref{lem-functorialbetweenSandmc3} for commutativity the two squares on the
  right hand side between row two and three --- and we use that 
  $$\Extsix(\ksix(\sfS\e{n}{2}),\ksix(\sfS\FFFF{1}{7}))=0,$$
  $$\Extsix(\ksix(\FFFF{n}{0}),\ksix(\sfS\sfS\e{n}{2}))=0.$$
  This is easily verified using projective resolutions. 

  Using naturality of $\mathbf{b}_-$ and
  Lemma~\ref{lem-functorialbetweenSandmc3}, it is easy to see that 
  $$\kkE(\mathbf{x}_{\FFFF{1}{1}}\times\widehat{\mc}(\mathbf{z}_n)\times\mathbf{w}_n\times\kkE((0,\phi_2,\phi_2)),e)
  =-\beta_{n,1},\quad\text{and}$$
  $$\kkE(\mathbf{b}_{\FFFF{1}{1}}\times\widehat{\sfS}\widehat{\mc}(\mathbf{z}_n)\times\widehat{\sfS}\mathbf{w}_n\times\kkE(\sfS(0,\phi_2,\phi_2)),e)
  =-\beta_{n,3}.$$
  Consequently, $\beta_{n,4}\circ f_{n,3}$ and $\beta_{n,1}\circ f_{n,0}$ are exactly the connecting
  homomorphisms of the cyclic six term exact sequence in the top (up
  to a sign, of course). 

  This proves exactness of the third of the four cyclic sequences in
  the theorem, for $i=2,5$. 

  This same result also works for $i=0,1,3,4$, by invoking
  Proposition~\ref{prop-rotate}. 

  That the last one of the sequences is exact for all $i=0,1,2$ is
  straightforward to check.
\end{proof}

\begin{proof}[Proof of the second part of Theorem~\ref{thm-diagrams}]
  Diagram~\eqref{eq:thmdiagramone}.
  First we prove it for $i=1$. 
  We have a commuting diagram 
  $$\xymatrix{
    0\ar@{^(->}[d]\ar@{^(->}[r] &
    \mc(\mathfrak{i}(\eo))\ar@{^(->}[d]^{\mc(\mathfrak{i}_{\eo})}\ar@{=}[r] &
    \mc(\mathfrak{i}(\eo))\ar@{^(->}[d] \\
    \mathfrak{i}(\ei)\ar@{=}[d]\ar@{^(->}[r]^{\mathfrak{i}_{\ei}} &
    \ei\ar@{->>}[d]^{\mc(\mathfrak{q}_{\eo})}\ar@{->>}[r]^{\mathfrak{q}_{\ei}} &
    \mathfrak{q}(\ei)\ar@{->>}[d] \\
    \mathfrak{i}(\ei)\ar@{^(->}[r] &
    \mc(\mathfrak{q}(\eo))\ar@{->>}[r] &
    \FFFF{1}{1}
    }$$
  of objects from $S\mathcal{E}$ with short exact rows and short exact
  columns. 
  Note that the short exact sequences
  $\extwithoutmaps{\mathfrak{i}(\ei)}{\mc(\mathfrak{q}(\eo))}{\FFFF{1}{1}}$
  and 
  $\extwithoutmaps{\mc(\mathfrak{i}(\eo))}{\mathfrak{q}(\ei)}{\FFFF{1}{1}}$ 
  are exactly 
  the short exact sequences 
  $\extwithmaps{\mathfrak{i}(\FFFF{1}{2})}{\FFFF{1}{2}}{\mathfrak{q}(\FFFF{1}{2})}{\mathfrak{i}_{\FFFF{1}{2}}}{\mathfrak{q}_{\FFFF{1}{2}}}$
  and 
  $\mc$ applied to the sequence
  $\extwithmaps{\mathfrak{i}(\eo))}{\Fo}{\FFFF{1}{0}}{x}{y}$ from
  Definition~\ref{def-Bocksteinoperations}, respectively. 
  Now apply $\kkE(-,e)$, then one easily shows
  the commutativity of the diagram (using the definitions of the
  different maps)
  $$\xymatrix@C=32pt@R=32pt{
    F_{1,1}\ar[d]_{\rho_{n,1}}\ar[r]^{\tilde{f}_{1,1}} &
    F_{1,2}\ar[d]_{\tildeh_{n,1}^{1,1,in}}\ar[dr]^{\tilde{f}_{1,2}} \\
    F_{n,1}\ar[dr]_{\beta_{n,1}}\ar[r]^{\tildeh_{n,1}^{n,1,in}} &
    H_{n,1}\ar[d]^{\tildeh_{n,1}^{1,1,out}}\ar[r]_{\tildeh_{n,1}^{n,1,out}} &
    F_{1,3} \\ 
    & F_{1,4}    }$$
  If we apply $\mc^k$ to the diagram, for $k=1,2,3,4,5$, 
  we obtain commutativity of the corresponding part of
  Diagram~\eqref{eq:thmdiagramone},
  for $i=2,3,4,5,0$, respectively  --- 
  this is, indeed, a very long and tedious verification
  using the identifications and results above. 

  Now we prove commutativity of the remaining square in
  Diagram~\eqref{eq:thmdiagramone} for $i=1$. 
  Note that there exists a commuting diagram
  $$\xymatrix{
    \sfS\FFFF{1}{1}\ar[d]\ar[r] &
    \sfS\FFFF{1}{0}\ar[d]^{\mathfrak{i}_{\ei}} \\
    \mathfrak{q}(\eii)\ar[r]_\phi &
    \ei,
    }$$
  where the horizontal morphisms are the unique morphism which are the
  identity on the extension algebra, and the vertical morphism from
  $\sfS\FFFF{1}{1}$ to $\mathfrak{q}(\eii)$ is the morphism induced by
  the \starhomo $\sfS\C\rightarrow\Ii$ in the extension $\ei$. 
  It is easy to see that
  $\mc(\mathfrak{i}_{\eo})$ is exactly 
  $\phi\circ w_n$, where $w_n$ is the morphism inducing $\mathbf{w}_n$. 
  Now we get commutativity of 
  $$\xymatrix{
    H_{n,1}\ar[d]^{\tildeh_{n,1}^{1,1,out}}\ar[r]^{\tildeh_{n,1}^{n,1,out}} &
    F_{1,3}\ar[d]^{\tilde{f}_{1,3}} \\
    F_{1,4}\ar[r]_{\times n} &
    F_{1,4}
    }$$
  by applying $\kkE(-,e)$ to the above diagram. 
  If we first apply $\mc^k$ to the diagram, for $k=1,2,3,4,5$, 
  we obtain commutativity of the corresponding square of 
  Diagram~\eqref{eq:thmdiagramone},
  for $i=2,3,4,5,0$, respectively.

  Diagram~\eqref{eq:thmdiagramtwo}. 
  We first prove it for $i=2$. 
  We have a commutative diagram
  $$\xymatrix{
    \mathfrak{i}(\mc^2(\sfS\FFFF{1}{0}))
    \ar[r]^{\mathfrak{i}_{\mc^2(\sfS\FFFF{1}{0})}}
    \ar[d]^{\sfS(x\circ z_n)} &
    \mc^2(\sfS\FFFF{1}{0})
    \ar[r]^{\mathfrak{q}_{\mc^2(\sfS\FFFF{1}{0})}}
    \ar[d]^{\mc^2(\mathfrak{i}_{\eo}\circ z_n)} &
    \mathfrak{q}(\mc^2(\sfS\FFFF{1}{0}))
    \ar[d]^{w_n\circ\mc(z_n)} \\
    \mathfrak{i}(\eii)
    \ar[r]^{\mathfrak{i}_{\eii}}
    \ar[d]^{\sfS y} &
    \eii
    \ar[r]^{\mathfrak{q}_{\eii}}
    \ar[d]^{\mc^2(\mathfrak{q}_{\eo})} &
    \mathfrak{q}(\eii) \\
    \sfS\FFFF{1}{0}
    \ar[r]^{x_{\FFFF{1}{0}}} &
    \FFFF{1}{3}
    }$$
  where $x_{\FFFF{1}{0}}$, $z_n$ and $w_n$ denote the morphisms
  inducing $\mathbf{x}_{\FFFF{1}{0}}$, $\mathbf{z}_n$
  and $\mathbf{w}_n$, respectively, and the first column is the suspension of
  the short exact sequence introduced in
  Definition~\ref{def-Bocksteinoperations} 
  (note that we do not claim the columns and rows to be exact). 

  A computation shows that this gives rise to a commutative diagram
  (by applying $\kkE(-,e)$) 
  $$\xymatrix@C=32pt@R=32pt{
    & F_{1,3}\ar[d]_{\tildeh_{n,2}^{1,1,in}}\ar[dr]^{\rho_{n,3}} \\
    F_{1,4}\ar[dr]_{\tilde{f}_{1,4}}\ar[r]^{\tildeh_{n,2}^{1,n,in}} &
    H_{n,2}\ar[d]^{\tildeh_{n,2}^{1,1,out}}\ar[r]_{\tildeh_{n,2}^{1,n,out}} &
    F_{n,3}\ar[d]^{-\beta_{n,3}} \\
    & F_{1,5}\ar[r]_{\tilde{f}_{1,5}} &
    F_{1,0}.
    }$$
  If we apply $\mc^k$ to the diagram, for $k=1,2,3,4,5$, 
  we obtain commutativity of the corresponding part of 
  Diagram~\eqref{eq:thmdiagramthree},
  for $i=3,4,5,0,1$, respectively  --- 
  this is, indeed, a very long and tedious verification
  using the identifications and results above. 

  Now we prove commutativity of the remaining square in
  Diagram~\eqref{eq:thmdiagramtwo} for $i=2$. 
  Note that there exists a commuting diagram
  $$\xymatrix{
    \eii\ar[r]^{\mathfrak{q}_{\eii}}\ar[d]^{(*,\ident,0)} &
    \mathfrak{q}(\eii)\ar[d]^{(0,\ident,0)} \\
    [\extwithoutmaps{\Iii}{\Iii}{0}]\ar[r] &
    [\extwithoutmaps{\Ii}{\Ii}{0}],
    }$$
  where the bottom horizontal morphism is the morphism induced by the
  \starhomo from $\Iii$ to $\Ii$ in the extension $\eii$. 
  It is easy to see that 
  there exists a commuting square
  $$\xymatrix{
    \sfS\FFFF{1}{1}
    \ar[d]_{\kkE(\sfS(0,\ident_{\C},0))}
    \ar[r]^{\mathbf{x}_{\FFFF{1}{1}}} &
    \FFFF{1}{4}\ar[d]^{\kkE(\mc^3(0,\ident_\C,0))=\kkE(\mc^2(\mathfrak{q}_{\FFFF{1}{2}}))} \\
    \sfS\FFFF{1}{0}\ar[r]_{\mathbf{x}_{\FFFF{1}{0}}} &
    \FFFF{1}{3}
    }$$
  in $\catkkE$. 
  Using that 
  $\kkE(\eii,\sfS\FFFF{1}{0})$ is naturally isomorphic to 
  $\Homsix(\ksix(\eii),\ksix(\sfS\FFFF{1}{0}))$ (since 
  $\Extsix(\ksix(\eii),\ksix(\sfS\sfS\FFFF{1}{0}))=0$), 
  we can show that the square
  $$\xymatrix{
    \eii\ar[r]^{\kkE(\mc^2(\mathfrak{q}_{\eo}))}\ar[d]_{\kkE((*,\ident,0))} &
    \mc^2(\mathfrak{q}(\eo))\ar[d]^{\mathbf{x}^{-1}_{\FFFF{1}{0}}} \\
    [\extwithoutmaps{\Iii}{\Iii}{0}]\ar[r] &
    \sfS\FFFF{1}{0}
    }$$
  anti-commutes in $\catkkE$, where the bottom horizontal map is the
  canonical identification. 
  Using all this, we can show that we have a commuting diagram
  $$\xymatrix{
    F_{1,3}\ar[r]^{\times n}\ar[d]_{\tilde{f}_{1,3}} &
    F_{1,3}\ar[d]^{\tildeh_{n,2}^{1,1,in}} \\
    F_{1,4}\ar[r]_{h_{n,2}^{1,n,in}} &
    H_{n,2}.
    }$$
  If we first apply $\mc^k$ to the diagrams, for $k=1,2,3,4,5$, 
  we obtain commutativity of the corresponding square of 
  Diagram~\eqref{eq:thmdiagramtwo},
  for $i=3,4,5,0,1$, respectively.
  
  Diagram~\eqref{eq:thmdiagramthree}.  
  First we prove it for $i=2$. 
  We have a commuting diagram of objects from $S\mathcal{E}$ 
  $$\xymatrix{
    0\ar@{^(->}[d]\ar@{^(->}[r] &
    \mc(\mathfrak{i}(\ei))\ar@{^(->}[d]^{\mc(\mathfrak{i}_{\ei})}\ar@{=}[r] &
    \mc(\mathfrak{i}(\ei))\ar@{^(->}[d] \\
    \mathfrak{i}(\eii)\ar@{=}[d]\ar@{^(->}[r]^{\mathfrak{i}_{\eii}} &
    \eii\ar@{->>}[d]^{\mc(\mathfrak{q}_{\ei})}\ar@{->>}[r]^{\mathfrak{q}_{\eii}} &
    \mathfrak{q}(\eii)\ar@{->>}[d] \\
    \mathfrak{i}(\eii)\ar@{^(->}[r] &
    \mc(\mathfrak{q}(\ei))\ar@{->>}[r] &
    \FFFF{n}{1}
    }$$
  with short exact rows and columns. 
  The short exact sequence
  $\extwithoutmaps{\mathfrak{i}(\eii)}{\mc(\mathfrak{q}(\ei))}{\FFFF{n}{1}}$
  is exactly 
  the short exact sequence
  $\extwithmaps{\mathfrak{i}(\Fii)}{\Fii}{\mathfrak{q}(\Fii)}{\mathfrak{i}_{\Fii}}{\mathfrak{q}_{\Fii}}$. 
  Moreover, the short exact sequence 
  $\extwithoutmaps{\mc(\mathfrak{i}(\ei))}{\mathfrak{q}(\eii)}{\FFFF{n}{1}}$ 
  is exactly 
  the short exact sequence 
  $\extwithoutmaps{\sfS\FFFF{1}{1}}{e}{\FFFF{n}{1}}$ induced by the
  extension $\ei\colon\extwithoutmaps{\sfS\C}{\Ii}{\Io}$, where 
  $e$ is $\extwithoutmaps{0}{\Ii}{\Ii}$.

  A computation shows that this gives rise to a commutative diagram
  (by applying $\kkE(-,e)$) 
  $$\xymatrix@C=32pt@R=32pt{
    F_{n,1}\ar[d]_{-\beta_{n,1}}\ar[r]^{\tilde{f}_{n,1}} &
    F_{n,2}\ar[d]_{\tildeh_{n,2}^{n,1,in}}\ar[dr]^{\tilde{f}_{n,2}} \\
    F_{1,4}\ar[dr]_{\times n}\ar[r]^{\tildeh_{n,2}^{1,n,in}} &
    H_{n,2}\ar[d]^{\tildeh_{n,2}^{n,1,out}}\ar[r]_{\tildeh_{n,2}^{1,n,out}} &
    F_{n,3} \\ 
        & F_{1,4}
    }$$
  If we apply $\mc^k$ to the diagram, for $k=1,2,3,4,5$, 
  we obtain commutativity of the corresponding part of 
  Diagram~\eqref{eq:thmdiagramthree},
  for $i=3,4,5,0,1$, respectively  --- 
  this is, indeed, a very long and tedious verification
  using the identifications and results above. 

  Now we prove commutativity of the remaining square in
  Diagram~\eqref{eq:thmdiagramthree} for $i=2$. 
  We have a commuting diagram
  $$\xymatrix{
    \sfS\FFFF{n}{1}\ar[d]\ar[r] &
    e\ar[d] \\
    \mathfrak{i}(\eii)\ar[r]^{\mathfrak{i}_{\eii}} &
    \eii,
    }$$
  where $e$ is the extension $\extwithmaps{0}{\Iii}{\Iii}{}{\ident}$, the
  map from $\sfS\FFFF{n}{1}$ to $e$ is the one induced by the map
  $\sfS\Io\rightarrow\Iii$ in $\eii$, the map from $\sfS\FFFF{n}{1}$
  to $\sfS\FFFF{n}{0}=\mathfrak{i}(\eii)$ is the unique morphism which is the identity on
  the extension algebra, and the morphism from $e$ to 
  $\eii$ is the unique morphism which is  the identity on
  the extension algebra. 
  It is elementary to see that if we compose the morphism from $e$ to
  $\eii$ with the canonical identification of $e$ with
  $\sfS\FFFF{1}{1}$, we get exactly the morphism
  $\mc(\mathfrak{i}_{\ei})$. 
  If $\phi$ denotes the obvious morphism from $\FFFF{n}{1}$ to
  $\FFFF{n}{0}$, it is elementary to show that $\mc^3(\phi)$ is
  $\mc^2(\mathfrak{q}_{\FFFF{n}{2}})$. 
  Using all this, we see that this gives rise to a commuting square 
  $$\xymatrix{
    H_{n,2}\ar[d]^{h_{n,2}^{n,1,out}}\ar[r]_{h_{n,2}^{1,n,out}} &
    F_{n,3} \ar[d]^{\tilde{f}_{n,3}} \\
    F_{1,4}\ar[r]_{\rho_{n,4}} &
    F_{n,4}.
    }$$
  If we first apply $\mc^k$ to the diagrams, for $k=1,2,3,4,5$, 
  we obtain commutativity of the corresponding square of 
  Diagram~\eqref{eq:thmdiagramthree},
  for $i=3,4,5,0,1$, respectively.
\end{proof}

\section{Examples}
\label{sec-examples}

In the article \cite{cmb}, the authors examine the invariant
$\underline{K}_{\mathrm{six}}$ of extensions. 
For an extension \extwithoutmaps{\A_0}{\A_1}{\A_2}, the invariant consists of the six term exact
sequences: 
\begin{equation*}
  \xymatrix{
    K_0(\A_0)\ar[r] &
    K_0(\A_1)\ar[r] &
    K_0(\A_2)\ar[d] \\
    K_0(\A_2)\ar[u] &
    K_0(\A_1)\ar[l] &
    K_0(\A_0)\ar[l] \\
    K_0(\A_0;\Z_n)\ar[r] &
    K_0(\A_1;\Z_n)\ar[r] &
    K_0(\A_2;\Z_n)\ar[d] \\
    K_0(\A_2;\Z_n)\ar[u] &
    K_0(\A_1;\Z_n)\ar[l] &
    K_0(\A_0;\Z_n)\ar[l] 
  }
\end{equation*}
together with all the Bockstein operations. 
A homomorphism between the invariants is a family of group
homomorphisms respecting all the above maps as well as all the
individual Bockstein operations.  
We let
$\operatorname{Hom}_\Gamma(\underline{K}_{\mathrm{six}}(e),\underline{K}_{\mathrm{six}}(e'))$
denote the group of such homomorphisms. 

In the article \cite{cmb}, 
the authors prove that  
there is a natural
homomorphism from $\kkE(e,e')$ to 
$\operatorname{Hom}_\Gamma(\underline{K}_{\mathrm{six}}(e),\underline{K}_{\mathrm{six}}(e'))$.
Moreover, 
the authors prove through a series of examples that this homomorphism is neither surjective nor
injective. 
We take a closer look at {this series of} examples here. 

\begin{example}
  We will compute the invariant $\kE(\e{n}{i})$ for $n\in\N_{\geq 2}$
  and $i=0,1,\ldots,5$.  The groups are as in the table below

  \noindent
  \begin{center}
    \begin{tabular}{|c|c|c|c|c|c|c|c|}
      \hline
      & \e{n}{0} & \e{n}{1} & \e{n}{2} & \e{n}{3} & \e{n}{4} & \e{n}{5} \\
      \hline
      $F_{1,0}$ & 0        & 0        & $\Z_n$   & \Z       & \Z       & 0        \\
      \hline
      $F_{1,1}$ & 0        & 0        & 0        & $\Z_n$   & \Z       & \Z       \\
      \hline
      $F_{1,2}$ & \Z       & 0        & 0        & 0        & $\Z_n$   & \Z       \\
      \hline
      $F_{1,3}$ & \Z       & \Z       & 0        & 0        & 0        & $\Z_n$   \\
      \hline
      $F_{1,4}$ & $\Z_n$   & \Z       & \Z       & 0        & 0        & 0        \\
      \hline
      $F_{1,5}$ & 0        & $\Z_n$   & \Z       & \Z       & 0        & 0        \\
      \hline
      $F_{k,0}$ & 0              & 0              & $\Z_{(n,k)}$ & $\Z_k$         & $\Z_k$         & $\Z_{(n,k)}$   \\
      \hline
      $F_{k,1}$ & $\Z_{(n,k)}$   & 0              & 0              & $\Z_{(n,k)}$ & $\Z_k$         & $\Z_k$         \\
      \hline
      $F_{k,2}$ & $\Z_k$         & $\Z_{(n,k)}$   & 0              & 0              & $\Z_{(n,k)}$ & $\Z_k$         \\
      \hline
      $F_{k,3}$ & $\Z_k$         & $\Z_k$         & $\Z_{(n,k)}$   & 0              & 0              & $\Z_{(n,k)}$ \\
      \hline
      $F_{k,4}$ & $\Z_{(n,k)}$ & $\Z_k$         & $\Z_k$         & $\Z_{(n,k)}$   & 0              & 0            \\
      \hline
      $F_{k,5}$ & 0              & $\Z_{(n,k)}$ & $\Z_k$         & $\Z_k$         & $\Z_{(n,k)}$   & 0 \\
      \hline
      $H_{k,0}$ & $\Z$           & 0              & 0              & $\Z_{(n,k)}$ & $\Z_{nk}$      & $\Z\oplus\Z_{(n,k)}$ \\
      \hline
      $H_{k,1}$ & $\Z\oplus\Z_{(n,k)}$ & $\Z$           & 0              & 0              & $\Z_{(n,k)}$ & $\Z_{nk}$      \\
      \hline
      $H_{k,2}$ & $\Z_{nk}$      & $\Z\oplus\Z_{(n,k)}$ & $\Z$           & 0              & 0              & $\Z_{(n,k)}$ \\
      \hline
      $H_{k,3}$ & $\Z_{(n,k)}$ & $\Z_{nk}$      & $\Z\oplus\Z_{(n,k)}$ & $\Z$           & 0              & 0              \\
      \hline
      $H_{k,4}$ & 0              & $\Z_{(n,k)}$ & $\Z_{nk}$      & $\Z\oplus\Z_{(n,k)}$ & $\Z$           & 0              \\
      \hline
      $H_{k,5}$ & 0              & 0              & $\Z_{(n,k)}$ & $\Z_{nk}$      & $\Z\oplus\Z_{(n,k)}$ & $\Z$           \\
      \hline
    \end{tabular}
  \end{center}
   
  First case ($\e{n}{0}$):
  Using the first sequence from Theorem~\ref{thm-diagrams} with $i=0$, we get that $H_{k,0}=\Z$.
  Using the first sequence from Theorem~\ref{thm-diagrams} with $i=1$,
  we get that $H_{k,4}=0$. 
  Using the first sequence from Theorem~\ref{thm-diagrams} with $i=2$,
  we get that $H_{k,2}$ is $\Z_{nk}$ and that $H_{k,5}=0$. 
  Using the third sequence from Theorem~\ref{thm-diagrams} with $i=0$, we get that $H_{k,3}=\Z_{(n,k)}$. 
  If we put $i=1$, we see from 
  the second sequence that {$H_{k,1}$} fits into a
  short exact sequence
  $$\xymatrix{0\ar[r]&\Z_{(n,k)}\ar[r]&H_{k,1}\ar[r]&\Z\ar[r]&0}.$$
  Since $\Z$ is projective, $H_{k,1}=\Z\oplus\Z_{(n,k)}$.

  The other cases follow by symmetry. 

  The $H_{k,i}$'s above tell us what $\kkE(\e{k}{i},\e{n}{j})$ is for
  all $k,n\in\N_{\geq 2}$ and $i,j=0,1,2,3,4,5$. 
  A lengthy, quite straightforward computation shows that the group of
  homomorphisms from $\kE(\e{k}{i})$ to $\kE(\e{n}{j})$ agrees with
  these groups for   all $k,n\in\N_{\geq 2}$ and $i,j=0,1,2,3,4,5$,
  and thus all of the counterexamples from \cite{cmb} become examples when we augment the invariant this way. 

  A more systematic proof of the above claims (and some generalizations) {are given} in the paper
  \cite{automorphism-se-gr-er}. 
\end{example}

\section*{Acknowledgements}
This work was carried out over a long period of time and has been supported by the Faroese Research Council, the NordForsk Research Network ``Operator
Algebras and Dynamics'' (grant \#11580) and the 
Danish National Research Foundation (DNRF) through the Centre for
Symmetry and Deformation.

\bibliographystyle{amsalpha}

\providecommand{\bysame}{\leavevmode\hbox to3em{\hrulefill}\thinspace}
\providecommand{\MR}{\relax\ifhmode\unskip\space\fi MR }
\providecommand{\MRhref}[2]{%
  \href{http://www.ams.org/mathscinet-getitem?mr=#1}{#2}
}
\providecommand{\href}[2]{#2}

\end{document}